\documentclass[12 pt]{report}
\usepackage[margin=1 in]{geometry}
\usepackage{graphicx}
\usepackage{amsmath}
\usepackage{booktabs}
\usepackage{multirow}
\usepackage{bigdelim}
\usepackage{float}
\usepackage{bm}
\usepackage{mathtools}
\usepackage{tgtermes}
\usepackage{lipsum}
\usepackage{pythonhighlight}
\usepackage{tabto}
\usepackage{subcaption}
\usepackage{setspace}
\usepackage{amsfonts}
\usepackage{setspace}
\usepackage[round, sort]{natbib}
\usepackage{titlesec}
\usepackage{algorithm}
\usepackage{algpseudocode}
\usepackage{titletoc}
\usepackage{tocloft}
\usepackage{apptools}
\usepackage{etoolbox}
\usepackage{indentfirst}
\usepackage{datetime}
\usepackage{caption}
\usepackage{subcaption}
\usepackage{amssymb,amsthm}
\usepackage{url}
\usepackage{ifthen}
\usepackage{pgfplots}
\usepackage{mathrsfs}
\usepackage{blindtext}
\usepackage{mdframed}
\usepackage{subfig}
\usepackage{fancyvrb}
\usepackage{pgfplots}

\usepackage{pgf,tikz}
\usetikzlibrary{arrows}

\usepackage{listings}
\usepackage{xcolor}

\definecolor{codegreen}{rgb}{0,0.6,0}
\definecolor{codegray}{rgb}{0.5,0.5,0.5}
\definecolor{codepurple}{rgb}{0.58,0,0.82}
\definecolor{backcolour}{rgb}{0.95,0.95,0.92}

\lstdefinestyle{mystyle}{
    backgroundcolor=\color{backcolour},   
    commentstyle=\color{codegreen},
    keywordstyle=\color{magenta},
    numberstyle=\tiny\color{codegray},
    stringstyle=\color{codepurple},
    basicstyle=\ttfamily\footnotesize,
    breakatwhitespace=false,         
    breaklines=true,                 
    captionpos=b,                    
    keepspaces=true,                 
    numbers=left,                    
    numbersep=5pt,                  
    showspaces=false,                
    showstringspaces=false,
    showtabs=false,                  
    tabsize=2
}

\newtheorem{theo}{Theorem}

\newenvironment{ftheo}[1][]
  {\begin{mdframed}\begin{theo}\ifthenelse{\equal{#1}{}}{}{\ (#1)}}
  {\end{theo}\end{mdframed}}

\newtheorem{defn}{Definition}

\newenvironment{fdefn}[1][]
  {\begin{mdframed}\begin{defn}\ifthenelse{\equal{#1}{}}{}{\ (#1)}\normalfont}
  {\end{defn}\end{mdframed}}

\newdateformat{monthdate}{%
  \monthname[\THEMONTH], \THEYEAR}
  
\apptocmd{\appendix}
  {\addtocontents{toc}{%
   \protect\renewcommand{\protect\cftchappresnum}{Appendix~}%
   \protect\settowidth\mylength{\bfseries\protect\cftchappresnum
   \protect\cftchapaftersnum}%
   \setlength{\cftchapnumwidth}{6em}}%
  }{}{}
\newlength\mylength

\renewcommand{\cftchapaftersnum}{:}

\renewcommand\cftchappresnum{Chapter~}
\newcommand{\RR}{\mathbb{R}}
\newcommand{\norm}[2]{\left \lVert #1 \right \rVert_{#2}}
\newcommand{\abs}[1]{\left| #1 \right|}

\titlespacing*{\chapter}{0pt}{-19pt}{12pt}  
\titlespacing*{\section}{0pt}{\baselineskip}{0pt}
\titlespacing*{\subsection}{0pt}{\baselineskip}{0pt}
\titlespacing*{\subsubsection}{0pt}{\baselineskip}{0pt}

\titleformat{\chapter}[hang]
  {\normalfont\normalsize\bfseries\centering}
  {Chapter \thechapter:}{3pt}{\normalsize}
  
\titleformat*{\section}{\bfseries\fontsize{12}{0}}
\titleformat*{\subsection}{\bfseries\fontsize{12}{0}}
\titleformat*{\subsubsection}{\bfseries\fontsize{12}{0}}

\cftsetrmarg{0pt}%
 
\begin{document}
\spacing{1.5}


\thispagestyle{empty}
    \begin{center}
      \uppercase{Adaptive Mesh Refinement for Variational Inequalities}\\
      \hfill \break by  \\
      Stefano Fochesatto\\
      B.S. University of Alaska Fairbanks\\
      \hfill \break
      A Project Submitted in Partial Fulfillment of the Requirements\\ 
       for the Degree of\\
      \hfill \break
      Master of Science\\
       in\\
      Mathematics\\    
      \  
         
      University of Alaska Fairbanks\\
      \monthdate\today \\
      
    \end{center}
  \hspace{6cm} APPROVED:
\begin{tabbing}
\hspace{7cm}\=\hspace{2cm}\=\kill
  \> Ed Bueler, Committee Chair \\
  \> Jill Faudree, Committee Member \\
  \> David Maxwell, Committee Member \\
  \> Leah Wrenn Berman, Chair \\
  \> \> Department of Mathematics and Statistics\\
  \> Karsten Hueffer, Dean\\
  \> \> College of Natural Science and Mathematics\\
  \> Richard Collins, Director\\
  \> \> Graduate School
\end{tabbing}
  \clearpage

\newpage\null\thispagestyle{empty}\newpage

\pagenumbering{roman}
\setcounter{page}{3}

\addcontentsline{toc}{chapter}{Copyright}
  \vspace*{10\baselineskip}
  \begin{center}
    Copyright by Stefano Fochesatto\\

    All Rights Reserved
  \end{center}
  \clearpage



\chapter*{Abstract}
\addcontentsline{toc}{chapter}{Abstract}
\noindent
Variational inequalities play a pivotal role in a wide array of scientific and engineering
applications. This project presents two techniques for adaptive mesh refinement (AMR) in the context of variational inequalities, with a specific focus on the classical obstacle problem. 

We propose two distinct AMR strategies: Variable Coefficient Elliptic Smoothing (VCES) and Unstructured Dilation Operator (UDO). VCES uses a nodal active set indicator function as the initial iterate to a time-dependent heat equation problem. Solving a single step of this problem has the effect of smoothing the indicator about the free boundary. We threshold this smoothed indicator function to identify elements near the free boundary. Key parameters such as timestep and threshold values significantly influence the efficacy of this method.

The second strategy, UDO, focuses on the discrete identification of elements adjacent to the free boundary, employing a graph-based approach to mark neighboring elements for refinement. This technique resembles the dilation morphological operation in image processing, but tailored for unstructured meshes.

We also examine the theory of variational inequalities, the convergence behavior of finite element solutions, and implementation in the Firedrake finite element library. Convergence analysis reveals that accurate free boundary estimation is pivotal for solver performance. Numerical experiments demonstrate the effectiveness of the proposed methods in dynamically enhancing mesh resolution around free boundaries, thereby improving the convergence rates and computational efficiency of variational inequality solvers. Our approach integrates seamlessly with existing Firedrake numerical solvers, and it is promising for solving more complex free boundary problems.

\newpage

\chapter*{Acknowledgements}
\addcontentsline{toc}{chapter}{Acknowledgements}	
I would like to express my deepest gratitude to my committee members, Dr. Ed Bueler, Dr. Jill Faudree, and Dr. David Maxwell, for their unwavering support and invaluable guidance throughout my entire mathematical education. Their expertise and encouragement have been instrumental in the completion of this thesis and my development as a mathematician.

I am also sincerely thankful to my cohort of fellow DMS graduate students. Your camaraderie, insightful discussions, and constant support have made this challenging journey both rewarding and memorable. Sharing this experience with all of you has been an inspiration. My heartfelt thanks go to the Firedrake community for their support and collaboration. 

Finally, I am profoundly grateful to my family for their unconditional love and encouragement. Your steadfast support has been the foundation of my perseverance and success. I could not have reached this milestone without you.
\noindent

\newpage

\renewcommand{\contentsname}{Table of Contents}
\tableofcontents
\newpage

{%
\renewcommand*{\addvspace}[1]{}
\let\oldnumberline\numberline
\renewcommand{\numberline}[1]{\figurename~\oldnumberline{#1:}}
\listoffigures%
}
\addcontentsline{toc}{chapter}{List of Figures}
\newpage

{%
\renewcommand*{\addvspace}[1]{}
\let\oldnumberline\numberline
\renewcommand{\numberline}[1]{\tablename~\oldnumberline{#1:}}
\listoftables%
}
\addcontentsline{toc}{chapter}{List of Tables}
\newpage

\pagenumbering{arabic}

\chapter{Introduction}

The goal of this project is to introduce two techniques for adaptive mesh refinement for free boundary problems (variational inequalities). We will consider only the classical obstacle problem as an example of a variational inequality. 

In the context of the finite element method (FEM), the discretization of a partial differential equation (PDE) is described by a partition of its domain into a finite number of elements (i.e, a mesh) and a finite dimensional function space (i.e., a finite element space). In a two dimensional domain, these elements are usually triangles or rectangles and the basis of the finite element space is composed of hat functions over vertices in the mesh, with support over neighboring elements. 

The convergence of finite element solutions is most commonly achieved via approximating on increasingly refined meshes by the use of a finite basis of piecewise polynomials of low degree; this method of convergence is often referred to as $h$-refinement. We also see convergence achieved by approximating via increasingly higher degree basis of piecewise polynomials over a coarse mesh; this is referred to as $p$-refinement. Several schema have been explored that take advantage of both methods of convergence, so-called $hp$-refinement finite elements, further discussion of such methods can be found in \citet{demkowicz_computing_2007}.

The classical obstacle problem consists of finding the equilibrium position of an elastic membrane with a fixed boundary after some force is applied, where the membrane is constrained to lie above some obstacle. Let $\Omega \subset \RR$ be a bounded domain. The problem is then formulated as a constrained minimization problem where we seek to find the position of the membrane $u(x)$, with fixed valued $u(x) = g_D(x)$ on $\partial \Omega$ with a load $f$ applied, and where $u(x)$ is constrained to be above an obstacle $\psi(x)$, 
\begin{align}
    \underset{u}{\text{ minimize: }}  &I(u) = \int_\Omega \frac{1}{2} \abs{\nabla u}^2 - fu \\
  \text{ subject to: } v &\geq \psi, \\
    u |_{\partial\Omega} &= g_D. 
  \end{align}
  The admissable set for such a problem can be described by 
  \begin{equation}
    K_\psi = \{u \in X| u \geq \psi \}
  \end{equation}
  where $X$ is a Sobolev space with boundary conditions $g_D$ enforced. 
  
  We can describe another, equivalent, variational inequality (VI) formulation of the obstacle problem 
  \begin{equation}
    \int_\Omega \nabla u \cdot \nabla(v - u) \geq \int_\Omega f(v - u), \quad \text{ for all } v \in K_\psi.
  \end{equation} 
  Equivalence of these formulations is proven in Chapter 3. From a solution to (1.1)-(1.3) or (1.5) we may identify the inactive and active sets $I_u$ and $A_u$,
\begin{equation}
  I_u = \{x \in \Omega | u(x) > \psi(x)\} \quad A_u = \Omega \setminus I_u.
\end{equation} 
With these definitions of the inactive and active set we can also define the free boundary $\Gamma_u = \partial I_u \cap \Omega$.

The same problem can also be defined by its strong form formulation
\begin{align}
  -\nabla^2 u &= f \quad \text{ on } I_u,\\
  u &= \psi \quad \text{ on } A_u,\\
  u &= g_D \quad \text{ on } \partial \Omega.
\end{align}
 Observe that this `naive' strong form is not sufficient to describe the problem as neither $A_u$ or $I_u$ are known a priori. A better strong formulation is as a complementarity problem (CP). In this form each statement holds over the entire domain $\Omega$. A solution $u$ satisfies the following, 
\begin{align}
  -\nabla^2 u - f \geq 0\\
  u - \psi \geq 0\\
  (-\nabla^2u - f)(u - \psi) = 0
\end{align}
Each of these formulations are instrumental in our understanding of the problem. The constrained energy minimization formulation can be helpful for understanding the physics of the problem. The class of problems described by the VI and CP formulations is a superset of the minimization formulation \citep[page 319]{bueler_petsc_2021}. As we will see, the CP formulation is also instrumental in the implementation of numerical solvers.

In terms of numerics, the constraint of $u \geq \psi$ makes this problem nonlinear so we are required to use an iterative solver. In this project our main solver will be VI-adapted Newton with Reduced-Space Line Search (VINEWTONRSLS). In Chapter 3 we will show that numerical methods will not converge until the active and inactive sets stabilize and the free boundary is identified. For solvers like VINEWTONRSLS which employ a Newton iteration, we find that they can only adjust the approximated free boundary by one-cell per iteration \citep{graser_multigrid_2009} and therefore convergence is tied to mesh resolution and is proportional to the number of grid spaces between initial iterate free boundary and discrete solution free boundary \citep[page 324]{bueler_petsc_2021}.

Adaptation involves altering the discretization to achieve a more desirable solution, whether that be with the goal of reducing $L_2$ error or reducing the error in some post-computation quantity like drag or lift. For example, consider the flow of some incompressable fluid through a pipe with an obstructive obstacle. Computing its drag coefficient would be an example of such a desired quantity \citep[Chapter 1.1]{bangerth_adaptive_2003}.

Further ``-adaptive" methods have also been explored. These methods are designed to increase mesh resolution or polynomial degree locally by means of a local error estimator, for some goal term often referred to as the Quantity of Interest (QoI), usually denoted as $J(\cdot)$. The most common of these employ the Dual Weighted Residual (DWR) method, derived in \citet{becker_feed-back_1996}. The DWR method, like most adaptive refinement techniques begins with a rigorous, a posteriori analysis of the quantity $J(u) - J(u_n)$. Estimation of this error quantity is then decomposed into local element-wise error estimators which can be used as a heuristic to tag elements for refinement. Finally a refinement strategy is employed to refine the mesh. Choices of whether to refine elements or increase the degree of polynomial basis functions, and by how much, must be made. 

The techniques outlined above are referred to as tagging methods, in which the error indicators are used to ``tag" elements for refinement. However, there are ``metric-based" methods which control the size, shape and orientation of the elements instead \citep{alauzet_metric-based_2010}.  For this project we will primarily focus on $h$-adaptive tagging methods.

As we will see in Chapter 3, convergence of VI problems is dominated by the error in approximating the free boundary. An adaptive refinement scheme that is able to concentrate effort around the solution free boundary will both enhance convergence properties and reduce unnecessary computation in the active set.

In this project we will introduce two adaptive refinement schemes which can identify the free boundary to a high degree of accuracy and permit the user to vary the spread of the refinement area. In the first strategy we compute a node-wise indicator function for the active set for use as an initial iterate in a time-dependent heat equation problem. Solving a single step of this problem has the effect of smoothing the indicator about the free boundary. The result of this smoothing is then averaged over each element and thresholded for refinement. The second strategy focuses on the discrete identification of elements adjacent to the free boundary, employing a graph-based approach to mark neighboring elements for refinement.

Our implementations of these methods are written using the Firedrake finite element library and produce conforming meshes with no hanging nodes. These high quality meshes are suitable for various VI solvers, including as coarse grids for the FASCD multilevel solver in \citet{bueler_full_2023}.

The remaining chapters are organized as follows:
\begin{enumerate}
  \item[2] Briefly describe the finite element method, with relevant results from FEM error analysis, and introduce the Firedrake library. 
  \item[3] Provide background for solving VIs numerically. Explain why the error arising from free boundary approximation dominates. 
  \item[4] Introduce the two new adaptive refinement schemes for VIs, VCES and UDO, and provide a detailed description of their implementations.
  \item[5] Present numerical results demonstrating the effectiveness of the proposed methods in enhancing mesh resolution around free boundaries. 
\end{enumerate}

\chapter{A Brief Introduction to the Finite Element Method}
\section{Reference Problem: Poisson Equation on Square Domain}
For our brief introduction to the finite element method we aim to construct a solution for the Poisson equation on a square domain. The following is a formulation of such a problem, with mixed boundary conditions.
\begin{mdframed}
    Given $\Omega = (-1, 1) \times (-1, 1)$, find $u$ such that, 
    \begin{align}  \label{eq:poisson}
      -\nabla^2u &= f \quad\text{in}\quad \Omega\\
       u = g_D \quad \text{on}\quad \partial \Omega_D \qquad \text{and} &\qquad u' = g_N \quad \text{on} \quad \partial\Omega_N
    \end{align}
\end{mdframed}

Our reference problem, as described in \citet{elman_finite_2005}, will have only Dirichlet boundary conditions. It will become relevant later that we have an analytic solution to such a problem in order to verify our error analysis. To do so, we choose a solution $u^*$ which is smooth over $\Omega$. Let
\begin{equation}
  u^*(x, y) = \dfrac{2(1 + y)}{(3 + x)^2 + (1 + y)^2}. 
\end{equation}
We find that $-\nabla^2u^* = 0$ so $f = 0$. Let $g_D = u^*|_{\partial\Omega}$.

\begin{figure}[H]
  \begin{center}
    
    \begin{subfigure}[b]{0.45\textwidth}
      \centering
      \includegraphics[width=\textwidth]{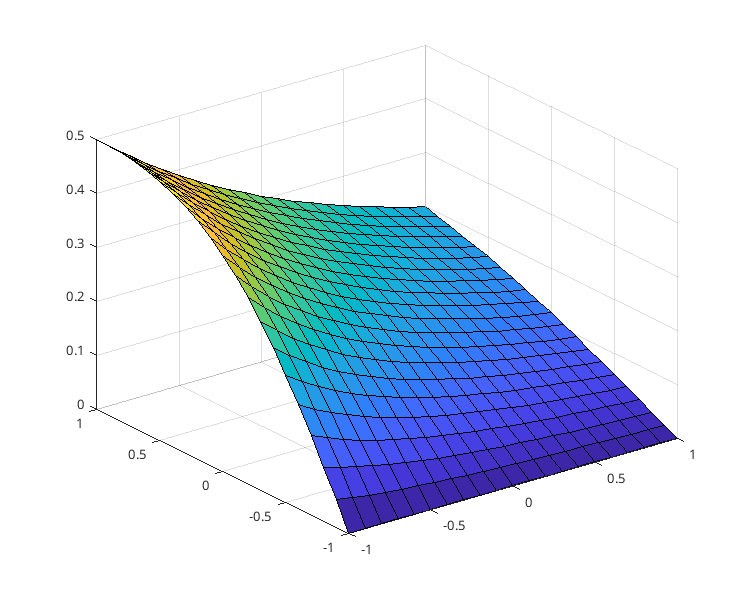}
      \label{fig:Reference Problem 3D}
    \end{subfigure}
    \hfill
    \begin{subfigure}[b]{0.45\textwidth}
      \centering
      \includegraphics[width=\textwidth]{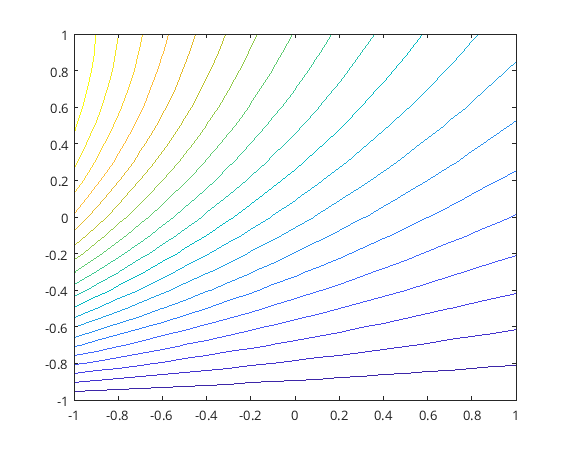}
      \label{fig:Reference Problem Contour}
    \end{subfigure}
    \hfill
    \caption{Plot of Reference Problem, surface plot (left), contour plot (right).}
    \label{fig: Reference Problem FEM Plots}
  \end{center}
\end{figure}

\section{Weak Formulation}
We say a function is a \emph{classical solution} to a boundary value problem if it is sufficiently differentiable over $\Omega$ to satisfy the strong form of the differential equation, for example (2.1). For our reference problem, a solution $u$ would need to be twice differentiable over $\Omega$ to be a classical solution. For problems defined over non-smooth domains or with discontinuous source functions there may be no solution satisfying such a differentiability requirement. This leads to an alternative description of the problem which expands the admissable set of solutions by "weakening" the differentiability requirements, hence we call this description the weak formulation. 

For our reference problem such a weak form is derived by the following functional. Let $v$ be a test function and multiply:
\begin{equation} \label{eq:weakformfunctional}
  \int_\Omega (\nabla^2u + f)v = 0.
\end{equation}
Distributing $v$, applying integration by parts and the divergence theorem we balance the differentiability requirements between $u$ and $v$:
\begin{align} \label{eq: weakformderivation}
  -\int_\Omega v\nabla^2u &= \int_\Omega \nabla u \cdot \nabla v - \int_\Omega \nabla \cdot (v \nabla u)\\ 
  &= \int_\Omega \nabla u \cdot \nabla v - \int_{\partial\Omega} v (\nabla u \cdot n)
\end{align}
Note that (2.5) follows when $\nabla \cdot (v \nabla u)$ is expanded via product rule. Substitution into (2.4) results in an integral form of the Poisson equation, 
\begin{equation}\label{eq: weakformpossonequation}
  \int_{\Omega} \nabla u \cdot \nabla v = \int_{\Omega} vf + \int_{\partial\Omega} v (\nabla u \cdot n). 
\end{equation}

Now we must define the function spaces which contain our solution and test functions. Furthermore we must incorporate the boundary conditions of our reference problem into the weak formulation. For expression (2.7) to be well defined, it is necessary for $f \in L_2(\Omega)$, $(\nabla u \cdot n) \in L^2(\partial \Omega)$ and $v, u \in H^1(\Omega)$ where 
\begin{equation} \label{eq: definitionofSobelevSpace}
  {H}^1(\Omega) \coloneqq \left\{u: \Omega \to \RR: u, \frac{\partial u}{\partial x}, \frac{\partial u}{\partial y} \in L_2(\Omega)\right\}. 
\end{equation}
\begin{equation}
  \norm{u}{H^1(\Omega)} = \left(\int_{\Omega} \abs{u}^2 + \abs{\nabla u}^2\right)^{1/2}.
\end{equation}

Incorporating the Dirichlet boundary conditions we define the following test and solution subspaces, 
\begin{equation} \label{eq test and solution subspaces}
  {H}_{g_D}^1(\Omega) \coloneqq \{u \in {H}^1(\Omega): u = g_D \text{ on } \partial \Omega_D\}
\end{equation}
\begin{equation} \label{eq test and solution subspaces}
  {H}_{0}^1(\Omega) \coloneqq \{v \in {H}^1(\Omega):  v = 0 \text{ on } \partial \Omega_D\}
\end{equation}
We incorporate Neumann boundary conditions by enforcing $\nabla u \cdot n = g_N$ in the weak form. Let $\Omega_N$ be the portion of $\partial\Omega$ where Neumann conditions are enforced. In our reference problem $\Omega_N = \emptyset$. So we arrive at the following weak formulation, 
\begin{mdframed}
  Find $u \in  {H}_{g_D}^1(\Omega) $ such that, 
  \begin{equation}  \label{eq:poisson}
    \int_{\Omega} \nabla u \cdot \nabla v = \int_{\Omega} vf + \int_{\partial\Omega_N} vg_N \text{ for all } v \in {H}_{0}^1(\Omega)
  \end{equation}
\end{mdframed}

\section{Galerkin Finite Element Method}
We will be approximating $u$ by considering finite-dimensional subspaces of the solution space $H_{g_D}(\Omega)$, and test space $H_{0}^1(\Omega)$, notated as $S^h_{g_D}$ and $S^h_{0}$ respectively. These subspaces are defined by a suitable discretization of $\Omega$ and some choice of basis functions. We then expand the weak formulation (2.12) via a finite basis of functions to get a system of equations which we can then solve for $u_h \in S^h_{g_D}$. Both $p$ and $h$ refinement will result in a better approximation of the solution space $H_{g_D}^1(\Omega)$ and therefore a better approximation $u_h \approx u$. This gives the finite-dimensional weak formulation, 
\begin{mdframed}
  Find $u_h \in  {S}^H_{g_D}(\Omega) $ such that, 
  \begin{equation}  \label{eq:poisson finite dim}
    \int_{\Omega} \nabla u_h \cdot \nabla v_h = \int_{\Omega} v_hf + \int_{\partial\Omega_N} v_hg_N \text{ for all } v_h \in {S}_{0}^h(\Omega)
  \end{equation}
\end{mdframed}

Assume that $S^h_{g_D}$ and $S^h_{0}$ are spaces of piecewise-linear and continuous functions of a triangulation of $\Omega$, with $g_D$ and zero boundary values, respectively. This restricts $g_D$ to be piecewise-linear on $\partial \Omega_D$. Let $\{\phi_1, \phi_2, \dots ,\phi_n\}$ be a basis for $S^h_{0}$. In order to incorporate Dirichlet conditions we extend this basis with $\phi_{n + 1}, \dots,\phi_{n + n_\partial}$ and solve for coefficients $u_j$ with $j \in \{n + 1, ..., n + n_\partial\}$. We enforce $u_j = (g_D)_j$ for $\{n + 1, ... n + n\partial \}$ where,  
\begin{equation}
  g_D = \sum_{j = n + 1}^{n + n_\partial} (g_D)_j \phi_j, \text{ on } \partial \Omega
\end{equation}
To illustrate the differences between interior and boundary basis functions consider Figure 2.2.

\definecolor{zzttqq}{rgb}{0.6,0.2,0}
\definecolor{ccqqqq}{rgb}{0.8,0,0}
\definecolor{zzqqqq}{rgb}{0.6,0,0}
\definecolor{qqqqff}{rgb}{0,0,1}

\begin{figure}[H]
    \begin{subfigure}{0.50\textwidth}
      \begin{tikzpicture}[scale = .5][line cap=round,line join=round,>=triangle 45,x=1cm,y=1cm]
        \clip(-3,-1) rectangle (13,8);
        \fill[color=zzttqq,fill=zzttqq,fill opacity=0.1] (4,4) -- (5.33,2.67) -- (1.33,1.33) -- cycle;
        \draw [line width=2.4pt,color=zzqqqq] (4,4)-- (0,0);
        \draw (0,0)-- (8,0);
        \draw (1.33,1.33)-- (5.33,2.67);
        \draw (4,4)-- (5.33,2.67);
        \draw (4,4)-- (12,4);
        \draw (12,4)-- (8,0);
        \draw (2.67,2.67)-- (6.67,4);
        \draw (2.67,2.67)-- (4,1.33);
        \draw (1.33,1.33)-- (2.67,0);
        \draw (0,0)-- (4,1.33);
        \draw (6.67,4)-- (8,2.67);
        \draw (5.33,2.67)-- (9.33,4);
        \draw (5.33,2.67)-- (6.67,1.33);
        \draw (4,1.33)-- (8,2.67);
        \draw (4,1.33)-- (5.33,0);
        \draw (2.67,0)-- (6.67,1.33);
        \draw (9.33,4)-- (10.67,2.67);
        \draw (8,2.67)-- (12,4);
        \draw (8,2.67)-- (9.33,1.33);
        \draw (6.67,1.33)-- (10.67,2.67);
        \draw (6.67,1.33)-- (8,0);
        \draw (5.33,0)-- (9.33,1.33);
        \draw (2.67,0)-- (6.67,4);
        \draw (5.33,0)-- (9.33,4);
        \draw (1.33,1.33)-- (9.33,1.33);
        \draw (2.67,2.67)-- (10.67,2.67);
        \draw [dash pattern=on 1pt off 1pt] (2.67,7.1)-- (4,4);
        \draw [dash pattern=on 1pt off 1pt] (2.67,7.1)-- (4.67,3.33);
        \draw [dash pattern=on 1pt off 1pt] (2.67,7.1)-- (5.33,2.67);
        \draw [dash pattern=on 1pt off 1pt] (2.67,7.1)-- (3.33,2);
        \draw [dash pattern=on 1pt off 1pt] (2.67,7.1)-- (1.33,1.33);
        \draw [color=zzttqq] (4,4)-- (5.33,2.67);
        \draw [color=zzttqq] (5.33,2.67)-- (1.33,1.33);
        \draw [color=zzttqq] (1.33,1.33)-- (4,4);
        \draw (2.67,7.1)-- (2.67,2.67);
        \draw (0,0)-- (4,4);
        \begin{scriptsize}
        \fill [color=qqqqff] (4,4) circle (1.5pt);
        \fill [color=qqqqff] (4,4) circle (1.5pt);
        \draw[color=zzqqqq] (-.75,2) node {Dirichlet Boundary};
        \fill [color=qqqqff] (1.33,1.33) circle (1.5pt);
        \fill [color=ccqqqq] (2.67,2.67) circle (1.5pt);
        \fill [color=qqqqff] (2.67,7.1) circle (1.5pt);
        \fill [color=qqqqff] (4.67,3.33) circle (1.5pt);
        \fill [color=qqqqff] (5.33,2.67) circle (1.5pt);
        \fill [color=qqqqff] (3.33,2) circle (1.5pt);
        \end{scriptsize}
      \end{tikzpicture}
      \label{fig:Reference Problem 3D}
    \end{subfigure}
    \hspace{.3cm}
    \begin{subfigure}{0.50\textwidth}
      \begin{tikzpicture}[scale = .50][line cap=round,line join=round,>=triangle 45,x=1.0cm,y=1.0cm]
        \clip(-1,-1) rectangle (13,8);
        \fill[color=zzttqq,fill=zzttqq,fill opacity=0.1] (6.67,4) -- (4.67,3.33) -- (2.67,2.67) -- (3.33,2) -- (4,1.33) -- (6,2) -- (8,2.67) -- (7.33,3.33) -- cycle;
        \draw (4,4)-- (0,0);
        \draw (0,0)-- (8,0);
        \draw (1.33,1.33)-- (5.33,2.67);
        \draw (4,4)-- (5.33,2.67);
        \draw (4,4)-- (12,4);
        \draw (12,4)-- (8,0);
        \draw (2.67,2.67)-- (6.67,4);
        \draw (2.67,2.67)-- (4,1.33);
        \draw (1.33,1.33)-- (2.67,0);
        \draw (0,0)-- (4,1.33);
        \draw (6.67,4)-- (8,2.67);
        \draw (5.33,2.67)-- (9.33,4);
        \draw (5.33,2.67)-- (6.67,1.33);
        \draw (4,1.33)-- (8,2.67);
        \draw (4,1.33)-- (5.33,0);
        \draw (2.67,0)-- (6.67,1.33);
        \draw (9.33,4)-- (10.67,2.67);
        \draw (8,2.67)-- (12,4);
        \draw (8,2.67)-- (9.33,1.33);
        \draw (6.67,1.33)-- (10.67,2.67);
        \draw (6.67,1.33)-- (8,0);
        \draw (5.33,0)-- (9.33,1.33);
        \draw (2.67,0)-- (6.67,4);
        \draw (5.33,0)-- (9.33,4);
        \draw (1.33,1.33)-- (9.33,1.33);
        \draw (2.67,2.67)-- (10.67,2.67);
        \draw [color=zzttqq] (6.67,4)-- (4.67,3.33);
        \draw [color=zzttqq] (4.67,3.33)-- (2.67,2.67);
        \draw [color=zzttqq] (2.67,2.67)-- (3.33,2);
        \draw [color=zzttqq] (3.33,2)-- (4,1.33);
        \draw [color=zzttqq] (4,1.33)-- (6,2);
        \draw [color=zzttqq] (6,2)-- (8,2.67);
        \draw [color=zzttqq] (8,2.67)-- (7.33,3.33);
        \draw [color=zzttqq] (7.33,3.33)-- (6.67,4);
        \draw [dash pattern=on 1pt off 1pt] (5.33,7.02)-- (6.67,4);
        \draw [dash pattern=on 1pt off 1pt] (5.33,7.02)-- (4.67,3.33);
        \draw [dash pattern=on 1pt off 1pt] (5.33,7.02)-- (2.67,2.67);
        \draw [dash pattern=on 1pt off 1pt] (5.33,7.02)-- (3.33,2);
        \draw [dash pattern=on 1pt off 1pt] (5.33,7.02)-- (4,1.33);
        \draw [dash pattern=on 1pt off 1pt] (5.33,7.02)-- (6,2);
        \draw [dash pattern=on 1pt off 1pt] (5.33,7.02)-- (8,2.67);
        \draw [dash pattern=on 1pt off 1pt] (5.33,7.02)-- (7.33,3.33);
        \draw (5.33,7.02)-- (5.33,2.67);
        \begin{scriptsize}
        \fill [color=qqqqff] (2.67,2.67) circle (1.5pt);
        \fill [color=ccqqqq] (5.33,2.67) circle (1.5pt);
        \fill [color=qqqqff] (6.67,4) circle (1.5pt);
        \fill [color=qqqqff] (4.67,3.33) circle (1.5pt);
        \fill [color=qqqqff] (3.33,2) circle (1.5pt);
        \fill [color=qqqqff] (4,1.33) circle (1.5pt);
        \fill [color=qqqqff] (6,2) circle (1.5pt);
        \fill [color=qqqqff] (8,2.67) circle (1.5pt);
        \fill [color=qqqqff] (7.33,3.33) circle (1.5pt);
        \fill [color=qqqqff] (5.33,7.02) circle (1.5pt);
        \end{scriptsize}
        \end{tikzpicture}
      \label{fig:Reference Problem Contour}
    \end{subfigure}
    \caption{ Basis functions on Dirichlet boundary node (left) and interior node (right).}
\end{figure}
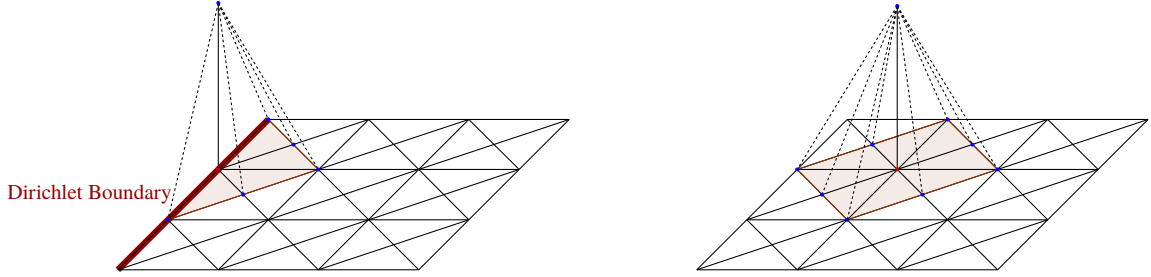

Then $u_h \in S^h_{g_D}$ is uniquely determined by $\{u_1, ..., u_n\} \in \RR$ in the expansion, 
\begin{equation}
  u_h = \underbrace{\sum_{j = 1}^{n} u_j \phi_j}_{\text{trial functions}} + \underbrace{\sum_{j = n + 1}^{n + n_\partial} (g_D)_j \phi_{j}.}_{\text{enforcing Dirichlet conditions}}
\end{equation}
Note that functions defined by the first sum are often referred to as trial functions. The space of trial functions now coincides with $S^h_{0}$, the same as the test function space. Expanding $u_h$ via (2.15) we get 
\begin{equation}
  \sum_{j = i}^n u_j \int_\Omega \nabla \phi_j \cdot \nabla v_h = \int_\Omega v_h f + \int_{\Omega_N} v_h g_N - \sum_{j = n + 1}^{n + n_\partial} u_j \int_\Omega \nabla \phi_{j + n} \cdot \nabla v_h \quad \text{ for all } \quad v_h \in {S}_{0}^h(\Omega)
\end{equation}
Then we set $v_h = \theta_i$. That is, we enforce the weak form on each basis function. 
\begin{equation}
  \sum_{j = i}^n u_j \int_\Omega \nabla \phi_j \cdot \nabla \phi_i = \int_\Omega \phi_i f + \int_{\Omega_N} \phi_i g_N - \sum_{j = n + 1}^{n + n_\partial} u_j \int_\Omega \nabla \phi_{j + n} \cdot \nabla \phi_i \quad \text{ for all }\quad \phi_i, i \in \{1, ..., n\}
\end{equation}
This system of equations is often referred to as a \textit{Galerkin system}. Let  
\begin{equation}
  a_{ij} = \int_\Omega \nabla \phi_j \cdot \nabla \phi_i
\end{equation}
\begin{equation}
  f_i = \int_\Omega \phi_i f + \int_\Omega \phi_i g_N - \sum_{j = n + 1}^{n + n_\partial} u_j \int_\Omega \nabla \phi_{j} \cdot \nabla \phi_i. 
\end{equation}
Then (2.17) becomes 
\begin{equation}
  Au = f
\end{equation}
The matrix $A$ formed in (2.20) is called the \textit{stiffness matrix}.

In summary, solving a PDE on a given mesh using FEM follows these general steps. 
\begin{enumerate}
  \item Derive the weak form via multiplication of a test function $v$ and integration by parts. Recall that $v$ should be chosen from a function which satisfies Dirichlet boundary conditions.
  \item Choose suitable basis $\{\phi_j\}_{j = 1}^{n}$ to expand the test and trial functions, minding the Dirichlet conditions on the boundary basis functions $\{\phi_j\}_{j = n + 1}^{n + n_\partial}$.
  \item Assemble and solve the resulting system of equations.  
\end{enumerate}

There are several implementation considerations, such as the representation and discretization of $\Omega$, quadrature methods for computing integrals in (2.18) and (2.19), and the subsequent solving step. Luckily the Firedrake Library \citep{ham_firedrake_2023} makes this process as high level as the steps described above, as we see next.

\section{Firedrake Implementation: Poisson Equation on Square Domain}

Firedrake \citep{ham_firedrake_2023} is an open source Python library for finite elements built on top of PETSc \citep{balay_petsc_2015} providing a high-level interface for the solution of partial differential equations and variational inequalities. As an introduction we will consider the Firedrake code to the reference problem defined in Section 2.1.

We begin by importing the Firedrake library and using a utility mesh constructor function.
\begin{Verbatim}[frame = single,numbers = left,baselinestretch=0.75]
  from firedrake import *
  mesh = RectangleMesh(10, 10, 1, 1, -1, -1)
\end{Verbatim}
Next we define the function space for our problem. In this example we are using piecewise linear functions. The ``CG" below stands for continuous Galerkin, meaning the basis functions are continuous over $\Omega$. The ``1" determines the degree of the basis functions. 
\begin{Verbatim}[frame = single,numbers = left,baselinestretch=0.75, firstnumber=last]
  V = FunctionSpace(mesh, "CG", 1)
  v = TestFunction(V)
\end{Verbatim}
Now we define the Dirichlet boundary conditions and the source term for our problem. These definitions use the Unified Form Language (UFL; \citep{UnifiedFormLanguage}) which is a domain specific language for defining variational forms. 
\begin{Verbatim}[frame = single,numbers = left,baselinestretch=0.75,firstnumber=last]
  x, y = SpatialCoordinate(mesh)
  f = Constant(0.0)
  gbdry = Function(V).interpolate((2*(1 + y))/((3 + x)**2 + (1 + y)**2))
\end{Verbatim}
Boundary indices are defined by the mesh constructor, which we match here.  
\begin{Verbatim}[frame = single,numbers = left,baselinestretch=0.75,  firstnumber=last]
  bdry_ids = (1, 2, 3, 4) 
  bcs = DirichletBC(V, gbdry, bdry_ids)
\end{Verbatim}
Next define the weak form of the problem using UFL. 
\begin{Verbatim}[frame = single,numbers = left,baselinestretch=0.75, firstnumber=last]
  u = Function(V)
  a = (inner(grad(u), grad(v))) * dx
  L = f*v * dx
\end{Verbatim}
The solve command will assemble the Galerkin system and solve it using the defined solver parameters. 
\begin{Verbatim}[frame = single,numbers = left,baselinestretch=0.75, firstnumber=last]
  solve(a == L, u, bcs=bcs, solver_parameters={'ksp_type': 'preonly', 
                                               'pc_type': 'lu'})
\end{Verbatim}
In this example the Galerkin system is solved directly using LU factorization. 

One could solve nonlinear PDEs using finite elements. For example consider the Lioville-Bratu equation, 
\begin{equation}
  \nabla^2u + \lambda e^u = 0.
\end{equation} 
The $\lambda e^u$ term makes this problem nonlinear. The resulting system of equations must be solved iteratively. To invoke the nonlinear solver in Firedrake, we put all the terms in the weak form on one side of the equation, passing ``F == 0" to the solver, and we specify a Newton solver.
\begin{Verbatim}[frame = single,numbers = left,baselinestretch=0.75, firstnumber=last]
  F = (dot(grad(u), grad(v)) - lambda * exp(u) * v) * dx
  
  solve(F == 0, u, bcs=bc, solver_parameters={'snes_type': 'newtonls',
                                              'ksp_type': 'preonly',
                                              'pc_type': 'lu'})
\end{Verbatim}
This combination of options will solve each linear Newton step directly using LU. The options regarding solver types and parameters in Firedrake are inherited from PETSc and are therefore extensive allowing for the user to specify the most appropriate solver for their problem \citep{bueler_petsc_2021}.

\section{Theory of Errors}
There are several results that we must understand with regards to finite element error analysis. These results motivate several ideas in adaptive mesh refinement. To begin, consider the following generic description of a weak formulation and its approximate solution constructed in the manner outlined in Section 2.3. 
\begin{fdefn}[Galerkin Approximation]
  Consider the linear variational problem of the form.
  Find $u \in H^1_{g_D}$ such that, 
  \begin{equation}
    a(u, v) = F(v) \text{ for all } v \in H^1_{0}.
  \end{equation}
  Where $a(\cdot, \cdot): H^1(\Omega) \times H^1(\Omega) \rightarrow \mathbb{R}$ is a bilinear form, and $F(\cdot): H^1(\Omega) \to \mathbb{R}$ is a bounded linear form.
  For a finite element space $S^h_{g_D} \subset H^1_{g_D}$ and $S^h_0 \subset S^h \cap H^1_{0}$, the \textit{Galerkin approximation} to the linear variational problem is $u_h \in S^h_{g_D}$ such that, 
  \begin{equation}
    a(u_h, v) = F(v) \quad \text{ for all }\quad v \in S^h_{0}
  \end{equation}
\end{fdefn}

\begin{ftheo}[Galerkin orthogonality]
  Given a linear variational problem (2.22), with solution $u$, and Galerkin approximation $u_h$ solving (2.23), then for all $v_h \in S^h_0$ we have that, 
  \begin{equation}
    a(u - u_h, v_h) = 0 \quad \text{ for all } \quad v_h \in S^h_0
  \end{equation}
\end{ftheo}
\begin{proof} Choose a test function $v_h \in S^h_0$. Since $S^h_0 \subset H^1_{0}$ we can apply both the weak form continuum problem and the finite dimensional problem, 
  \begin{equation}
    a(u, v_h) = F(v_h)
  \end{equation}
  \begin{equation}
    a(u_h, v_h) = F(v_h)
  \end{equation}
  Since $a$ is linear in the first argument, subtraction yields
  \begin{equation}
    a(u - u_h, v_h) = 0 \text{ for all } v_h \in S^h_0.
  \end{equation}
\end{proof}

Next we demonstrate the equivalence of norms $\norm{\nabla(\cdot)}{L_2}$ and $\norm{\cdot}{H^1}$, which will be used in the following analysis. Recall that the $H^1$ norm is defined in (2.9). To prove the equivalence of norms we will require the Poincare-Freidrichs inequality which we will quote from \citet[Chapter 2]{braess_finite_2007} without proof.
\begin{ftheo}[Poincare-Freidrichs Inequality] If $\Omega$ is a bounded domain and $u \in H_{0}^1(\Omega)$t then there exists a $C_1 > 0$ depending only on $\Omega$ such that,
  \begin{equation}
    \int_\Omega \abs{u}^2 \leq C_1 \int_\Omega \abs{\nabla u}^2
  \end{equation}
\end{ftheo}
To show the equivalence of norms, let $u \in H_{0}^1$ and note the by definition of the $L_2$ and $H^1$ norms we have, 
  \begin{equation}
    \norm{\nabla u}{L_2} \leq \norm{u}{H^1}.
  \end{equation}
  By Poincare-Freidrichs we get, 
  \begin{align}
    \norm{u}{H^1}^2 &= \int_\Omega \abs{u}^2 + \abs{\nabla u}^2\\
    &\leq \int_\Omega C_1\abs{\nabla u}^2 + \abs{\nabla u}^2\\
    &= (1 + C_1) \int_\Omega \abs{\nabla u}^2\\
    &= (1 + C_1)\norm{\nabla u}{L_2}^2
  \end{align}
Hence, 
\begin{equation}
 \norm{\nabla u}{L_2} \leq  \norm{u}{H^1} \leq \sqrt{1 + C_1} \norm{\nabla u}{L_2}.
\end{equation}

The following theorem applies to linear variational problems with symmetric bilinear forms; this case is covered in \citet[Chapter 6.2]{wait_finite_1985}. We will consider the Poisson case,
\begin{equation}
  a(u, v) = \int_{\Omega} \nabla u \cdot \nabla v.
\end{equation}
\begin{ftheo}[Best Approximation Property]
  \begin{equation}
    \norm{\nabla u - \nabla u_h}{L_2} = \min\{ \norm{\nabla u - \nabla v_h}{L_2} : v_h \in S^h_{g_D}\}
  \end{equation}
\end{ftheo}
\begin{proof}
  Let $v_h \in S^h_{g_D}$. Note that $a(\cdot, \cdot)$ defines an inner product on $H^1_{0} \times H^1_{0}$ with associated norm $\norm{\nabla u}{L_2} = \sqrt{a(u, u)}$. Note that since $u - u_h \in H^1_{0}$ (again it is assumed that we have a conforming method i.e. $S^h_0 \subset H^1_{0}$) by definition we get, 
  \begin{align*}
    \norm{\nabla(u - u_h)}{L_2}^2 &= a(u - u_h, u - u_h)\\
    &= a(u - u_h, u - v_h + v_h - u_h)\\
    &= a(u - u_h, u - v_h) + a(u - u_h,v_h - u_h)\quad \text{ bilinearity}\\
    &= a(u - u_h, u - v_h) \quad \text{Galerkin orthogonality}\\
    &\leq \norm{\nabla(u - u_h)}{L_2} \norm{\nabla(u - v_h)}{L_2} \quad \text{Cauchy-Schwarz}
  \end{align*}
  So it follows that, 
  \begin{equation}
    \norm{\nabla(u - u_h)}{L_2} \leq \norm{\nabla(u - v_h)}{L_2} \text{ for all } v_h \in S^h_{g_D}
  \end{equation}
\end{proof}

The following theorem generalizes the result to linear variational problems with asymmetric bilinear forms and is generally referred to as quasi-optimality.
\begin{ftheo}[Cea's Lemma; e.g. \citep{elman_finite_2005}]
  Let $u$ be the solution to a linear variational problem on $H^1$ and $u_h$ be the Galerkin approximation on $S^h$. Recall the following, 
  \begin{enumerate}
    \item[(a)] $F(\cdot)$ is a bounded linear form, so there exists a $C > 0$ such that, 
    \begin{equation}
      \abs{F(v)} \leq C \norm{v}{H^1}, \forall v \in H^1.
    \end{equation}
    \item[(b)] $a(\cdot, \cdot)$ is continuous bilinear form, so there exists a $\gamma > 0$ such that,
    \begin{equation}
      \abs{a(u, v)} \leq \gamma \norm{u}{H^1} \norm{v}{H^1}, \forall u, v \in H^1.
    \end{equation}
    
    \item[(c)] Assume that $a(\cdot, \cdot)$ is coercive, so there exists an $\alpha > 0$ such that, 
    \begin{equation}
      a(u, u) \geq \alpha \norm{u}{H^1}^2, \forall u \in H^1.
    \end{equation}

  \end{enumerate}

  Then, 
  \begin{equation}
    \norm{u - u_h}{H^1} \leq \frac{\gamma}{\alpha} \min_{v \in S^h}\norm{u - v}{H^1}
  \end{equation}
\end{ftheo}

Choosing the test function as the linear interpolant of $u$ in $S^h$ we find that Cea's Lemma implies that the approximation error is bounded by a multiple of the error in interpolating $u$ in $S^h$. Let $\pi_h(u)$ be this linear interpolant of $u$ in $S^h$. Then 
\begin{equation}
  \norm{u - u_h}{H^1} \leq \frac{\gamma}{\alpha} \norm{u -\pi_h(u)}{H^1}.
\end{equation}

Next we relate $\norm{u - \pi_h(u)}{H^1}$ to properties of the triangulation. This requires several definitions and arguments regarding a reference element, a mapping function to the reference element, its derivatives, and the aspect ratios of the elements in the discretization. For brevity we will cite the result from \citet{elman_finite_2005} and briefly outline a few key points. 

\begin{fdefn}[shape regular; e.g. \citep{elman_finite_2005}]
  A family of triangulations $\{T_h\}$ is shape regular if there exists a minimum angle $\theta_* \neq 0$ such that every element in $T_h$ satisfies $\theta_T \geq \theta_*$.
\end{fdefn}

The following is a key result which bounds the $H^1$ error of the interpolant in terms of the second derivatives of $u$ and $h$. 

\begin{ftheo}[e.g. \cite{elman_finite_2005}]
  Let $u \in H^2$. Let $T_h$ be a triangulation, define
  \begin{equation}
    \norm{D^2u}{}^2 := \int_\Omega \left(\left(\frac{\partial^2 u}{\partial x^2}\right)^2 + \left(\frac{\partial^2 u}{\partial x\partial y}\right)^2+\left(\frac{\partial^2 u}{\partial y^2}\right)^2\right)
  \end{equation}
  and let $h_k$ be the largest length of $\triangle_k \in T_h$, then there exists some constant $C_2$ such that
  \begin{equation}
    \norm{\nabla(u - \pi_h(u))}{L_2}^2 \leq C_2\sum_{{\triangle_k}\in {T_h}} h_k^2 \norm{D^2u}{\triangle_k}^2.
   \end{equation}
\end{ftheo}

Here are the key arguments summarized from \citep{elman_finite_2005} which are used in proving the result above. First let $\triangle_*$ be a reference element and $\overline{(\cdot)}$ denote the mapping between any given element in $T_h$ to the reference element. The first key result is bounding the element-wise interpolation error by the element-wise interpolation error over the reference element, 
\begin{equation}
  \norm{\nabla(u - \pi_h(u))}{\triangle_k}^2 \leq \underbrace{\frac{2}{\abs{\triangle_k}}}_{\text{aspect ratio term}}h_k^2\norm{\nabla(\overline{u} - \pi_h(\overline{u}))}{\triangle_*}
\end{equation}
The ratio between the area of the element and reference element is the Jacobian determinant of the  change of variables. The Bramble Hilbert lemma \citep{elman_finite_2005} is applied to bound the interpolation error over the reference element with the norm of second derivative of $\overline{u}$. 
\begin{equation}
  \norm{\nabla(\overline{u} - \pi_h(\overline{u}))}{\triangle_*}^2 \leq C_2 \norm{D^2\overline{u}}{\triangle_*}^2
\end{equation}
We apply another change of variables to $\norm{D^2\overline{u}}{\triangle_*}^2$ to get a bound in terms of $\norm{D^2u}{\triangle_k}^2$,
\begin{equation}
  \norm{D^2\overline{u}}{\triangle_*} \leq 18h_k^2 \frac{h_k^2}{\abs{\triangle_k}}\norm{D^2u}{\triangle_k}
\end{equation}

A geometric argument converts the aspect ratio term to one based on minimum angle over an element,
\begin{equation}
  \norm{\nabla(u - \pi_h(u))}{L_2}^2 \leq C_2\sum_{{\triangle_k}\in {T_h}} \frac{1}{\sin^2 \theta_k} h_k^2 \norm{D^2u}{\triangle_k}^2.
 \end{equation}
Shape regularity is applied to get a global bound on the interpolation error, 
\begin{equation}
  \norm{\nabla(u - \pi_h(u))}{L_2}^2 \leq C_2\sum_{{\triangle_k}\in {T_h}} h_k^2 \norm{D^2u}{\triangle_k}^2 \leq C_2h^2 \norm{D^2u}{}^2.
\end{equation}
Further details of the proof can be found in \citet{elman_finite_2005}. 

Now Cea's Lemma,  the Poincare-Friedrichs inequality, and Theorem 5 combine to show that for a problem like the reference problem introduced in Section 2.1 we can expect a convergence rate of $O(h)$ in the $H^1$ norm,
\begin{align*}
  \norm{u - u_h}{H^1} &\leq \frac{\gamma}{\alpha}\norm{u - \pi_h(u)}{H^1},\\
  &\leq \frac{\gamma}{\alpha} \sqrt{1 + C_1} \norm{\nabla(u - \pi_h(u))}{L_2},\\
  &\leq \frac{\gamma}{\alpha} \sqrt{1 + C_1} C_2 h \norm{D^2u}{},\\
  &= O(h).
\end{align*}

A further analysis in \citet{brennerscott}, using a different argument demonstrates an expected convergence rate of $O(h^2)$ in the $L_2$ norm. 

Using our Firedrake implementation on our reference problem described Section 2.1 we can easily verify these results, Figure 2.3 shows the result. 
\begin{figure}[H]
  \begin{center}
      \includegraphics[width=.7\textwidth]{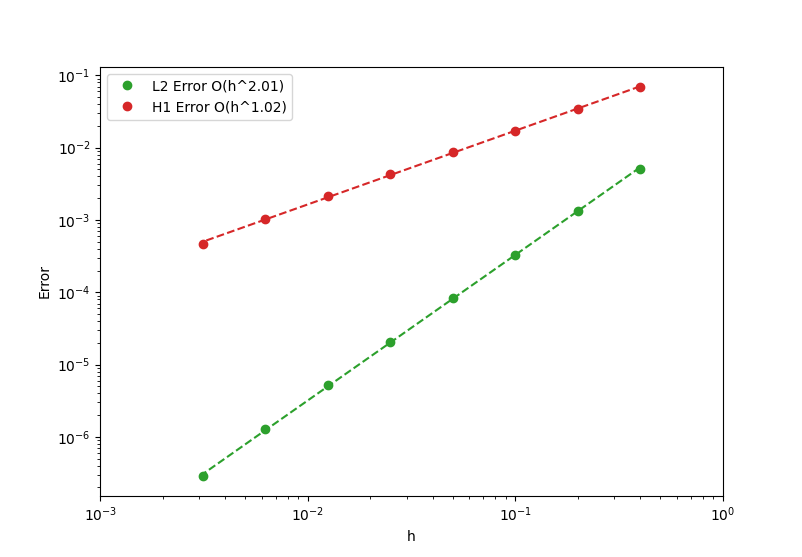}
      \label{fig:Reference Problem Convergence}
    \caption{Convergence of reference problem}
  \end{center}
\end{figure}

As we explore adaptive mesh refinement starting in Chapter 4, maintaining shape regularity will be key to preserving the convergence rates shown in Figure 2.3. One could easily consider a mark and refine strategy which fails to keep the minimum angle bounded away from zero. Consider the refinement strategies for an equilateral triangular element in Figure 2.4. Note the refinement strategy used in the sequence of triangulation on the top preserves the minimum angle, while the strategy used on the bottom sequence does not.

\begin{figure}[H]
  \begin{center}
  \includegraphics[width = .8\textwidth]{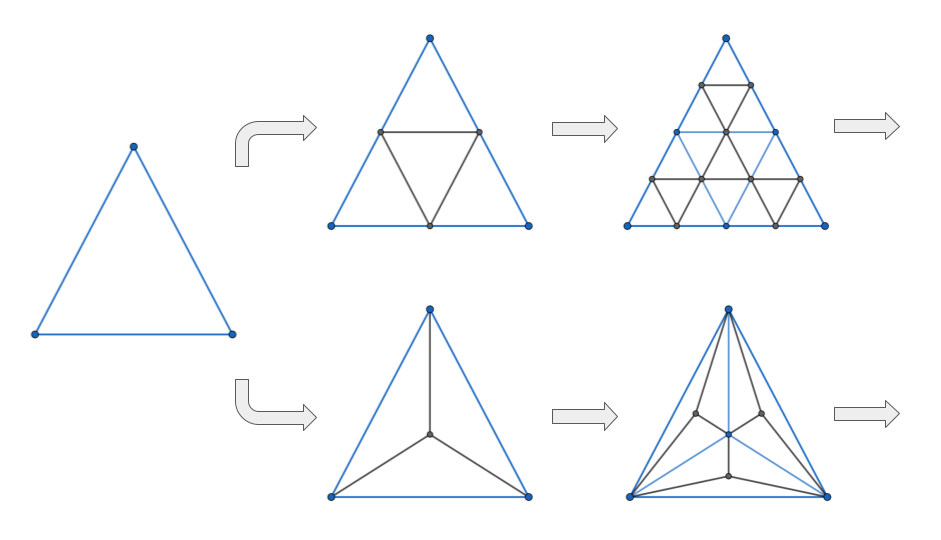}
\caption{Refinement schemes.}
\end{center}
\end{figure}

Another straightforward observation is that we can balance the $h_k$ and $\norm{D^2u}{\triangle_k}$ terms in Theorem 5 by refining the mesh in regions where the second derivative of the solution is large. This observation is central to many implementations of AMR for PDEs \citep[page 48, Figure 1.19]{elman_finite_2005}.

\chapter{Introduction to Variational Inequalities}
As mentioned in the Introduction, obstacle problems have several different formulations which we will prove are equivalent, or nearly so, in the following theorem. 
\begin{ftheo} Fix $f \in L^2(\Omega)$, $\psi \in C(\overline{\Omega})$, $g_D \in C(\overline{\Omega})$, $g_D \geq \psi$, $X = H^1(\Omega)$, and
  \begin{equation}
    K_\psi = \{v \in X| v \geq \psi\}
  \end{equation}
   where $X$ is a Sobolev space with boundary conditions equal to $g$ as before. Then the following statements are equivalent for a solution $u \in K_\psi$:
  \begin{enumerate}
    \item[(a)] $u$ is a solution to the energy minimization formulation, 
\begin{align}
    \underset{u \in K_\psi}{\text{ minimize: }}  &I(u) = \int_\Omega \frac{1}{2} \abs{\nabla u}^2 - fu 
  \end{align}
    \item[(b)] $u$ is a solution to the variational inequality formulation,
      \begin{equation}
    \int_\Omega \nabla u \cdot \nabla(v - u) \geq \int_\Omega f(v - u), \quad \text{ for all } v \in K_\psi.
  \end{equation} 
If also $u \in C(\overline{\Omega}) \cap C^2(\Omega)$ (is a classical solution), then (a) or (b) implies
    \item[(c)] $u$ is a solution to the complementarity problem formulation, for which the following hold over $\Omega$ a.e.
    \begin{subequations}
      \label{ncp}
      \begin{align}
        -\nabla^2 u - f \geq 0\label{ncp:1}\\
        u - \psi \geq 0\label{ncp:1}\\
        (-\nabla^2u - f)(u - \psi) = 0\label{ncp:1}
      \end{align}
    \end{subequations}
  \end{enumerate}
\end{ftheo}

\begin{proof}
  Suppose $u$ is a solution to (3.3) then it follows that for any $v \in K_\psi$, 
  \begin{align*}
    I(v) &= I(u + (v - u)),\\
    &= \frac{1}{2} \int_\Omega \abs{\nabla(u + (v - u))}^2 - f(u + (v - u)),\\
    &= \frac{1}{2} \int_\Omega \abs{\nabla u}^2+ \int_\Omega \nabla u \cdot \nabla(v - u) + \frac{1}{2} \int_\Omega \abs{\nabla(v - u)}^2 - \int_\Omega fu - \int_\Omega f(v - u),\\
    &= \underbrace{\frac{1}{2} \int_\Omega \abs{\nabla u}^2 - \int_\Omega fu}_{I(u)} + \underbrace{\frac{1}{2} \int_\Omega \abs{\nabla(v - u)}^2}_{\geq 0} + \underbrace{\int_\Omega \nabla u \cdot \nabla(v - u) - \int_\Omega f(v - u)}_{\geq 0},\\
    &\geq I(u). 
  \end{align*}
  Hence $u$ satisfies the energy minimization formulation (3.2). 

  Conversely suppose $u$ is a solution (3.2). Since $K_\psi$ is convex, for any test function $v \in K_\psi$ and $\epsilon \in [0, 1]$ then,
  \begin{equation}
      u + \epsilon(v - u) = \epsilon v + (1 - \epsilon)u \in K_\psi
  \end{equation} 
  By equation (3.2), 
  \begin{equation}
    I(u) \leq I(u + \epsilon (v - u)).
  \end{equation}
  So
  \begin{equation}
    0 \leq I(u + \epsilon (v - u)) - I(u).
  \end{equation}
  Note that 
  \begin{align}
    I(u + \epsilon(v - u)) &= \int_{\Omega} \frac{1}{2}\abs{\nabla (u + \epsilon(v - u))}^2 - \int_{\Omega}f(u + \epsilon(v - u)),\\
    &= \underbrace{\int_{\Omega} \frac{1}{2} \abs{\nabla u}^2}_{*} + \epsilon \int_\Omega \nabla u \cdot \nabla(v - u) + \epsilon^2 \int_\Omega \frac{1}{2} \abs{\nabla(v - u)}^2 \underbrace{-\int_\Omega fu}_{*} - \epsilon\int_\Omega f(v - u).
  \end{align}
  Expanding the right hand side of (3.7) we find that the $*$ terms cancel in (3.9) to get the following, 
  \begin{align*}
   0 &\leq  I(u + \epsilon (v - u)) - I(u),\\ 
   &=  \epsilon \int_\Omega \nabla u \cdot \nabla(v - u) + \epsilon^2 \int_\Omega \abs{\nabla(v - u)}^2 - \epsilon\int_\Omega f(v - u).
  \end{align*}
  Dividing by $\epsilon$ then we find (3.3) follows by sending $\epsilon \to 0$
  \begin{align*}
    0&\leq  \int_\Omega \nabla u \cdot \nabla(v - u) + \epsilon \int_\Omega \abs{\nabla(v - u)}^2 - \int_\Omega f(v - u),\\
   &= \int_\Omega \nabla u \cdot \nabla(v - u) - \int_\Omega f(v - u).
  \end{align*}
  So we conclude that $u$ is a solution to the variational inequality formulation. Now suppose that $u \in C(\overline{\Omega}) \cap C^2(\Omega)$ satisfies the variational inequality (3.3). We will show (3.4). Consider the inactive set $I_u$
  \begin{equation}
    I_u = \{x \in \Omega : u(x) > \psi(x)\}.
  \end{equation}
  Let $\phi \in C^\infty(\Omega)$ have compact support on $I_u$. Since $u$ is continuous there exists an $\epsilon > 0$ such that $v = u + t\phi \in K_\psi$ whenever $\abs{t} < \epsilon$. This is because the function $\phi$ has compact support in $I_u$ and is continuous. Now applying the variational inequality (3.3), and integration by parts, we get 
  \begin{align*}
    0 &\leq \int_\Omega \nabla u \cdot \nabla(v - u) - \int_\Omega f(v - u),\\
    &= t\left(\int_\Omega \nabla u \cdot \nabla\phi - \int_\Omega f\phi\right),\\
    &= t\left(-\int_\Omega \phi(\nabla^2u) + \int_{\partial \Omega} \phi(\nabla u \cdot n) - \int_\Omega f\phi \right),\\
    &= t\left(\int_\Omega (-\nabla^2u - f)\phi\right).
  \end{align*}
  Since this holds for any $\phi$ and $t$ of either sign it must follow that,
  \begin{equation}
    -\nabla^2u - f = 0 \text{ a.e. on } I_u.
  \end{equation}

  Now let $\phi \in C^\infty(\Omega)$ be nonnegative with compact support over the interior of the closed set $A_u := (I_u)^c$. For any $\epsilon > 0$ note $v = u + \epsilon\phi \in K_\psi$. We consider the variational inequality and apply integration by parts as before to conclude that
  \begin{equation}
   0 \leq \epsilon \left(\int_\Omega (-\nabla^2u - f)\phi\right).
  \end{equation}
  Since this holds for any such $\phi$ it must follow that
  \begin{equation}
    -\nabla^2u - f \geq 0 \text{ a.e. on } A_u.
  \end{equation}
  By definition of $A_u$, since $A_u \cup I_u = \Omega$ we have
  \begin{equation}
    -\nabla^2u - f \geq 0 \text{ a.e. on } \Omega.
  \end{equation}
  For emphasis, note that (3.12) holds both in $A_u$ and $I_u$. Since $u \in K_\psi$ by definition we have that $u - \psi \geq 0$. Furthermore since $u - \psi = 0$ on $A_u$ and since we have that $-\nabla^2u - f = 0$ a.e. on $I_u$, we can conclude the complementarity condition $(-\nabla^2 - f)(u - \psi) = 0$ a.e. Thus $u$ is a solution to (3.4). 
\end{proof}
The equivalence of these formulations is key to becoming familiar with free boundary problems. The energy minimization formulation is ideal for understanding the physical framework of the problem. The variational inequality formulation is necessary for defining the more broader class of free boundary problems whose solutions may not be derived from a potential function. The complementarity formulation reframes the problem in a way which can be solved using a variety of iterative methods.

\section{Reference Problem: Ball Obstacle Problem}
We will proceed by deriving an exact solution to a reference case for the obstacle problem \citep[Chapter 12]{bueler_petsc_2021}. This is a radially-symmetric problem. Consider $\Omega = [-2, 2]^2, f = 0$ and define the obstacle $\psi(r)$ as follows, 
\begin{equation}
  \psi(r) = \begin{cases}
    \sqrt{1 - r^2} & \text{ if } r \leq r_0\\
    \ell(r) & \text{ if } r > r_0
  \end{cases},
\end{equation}
where $r = \sqrt{x^2 + y^2}$, $r_0 = .9$, and $\ell(r) = \psi(r_0) + \psi'(r_0)(r - r_0)$. Note that $\psi$ is a hemisphere of radius 1 with a linear and continuous differentiable extension from $r = .9$ and onwards. The solution $u$ will satisfy the Poisson equation over the inactive set $I_u$. We consider the $\nabla^2$ operator in polar coordinates, 
\begin{equation}
  \nabla^2 = \frac{\partial^2 }{\partial r^2} + \frac{1}{r}\frac{\partial }{\partial r} + \frac{1}{r^2}\frac{\partial^2 }{\partial \theta^2}.
\end{equation}
Since our solution is radially-symmetric, the derivatives with respect to $\theta$ vanish, and clearing denominators the PDE simplifies to the following ODE, 
\begin{equation}
  ru''(r) + u'(r) = 0.
\end{equation}
Let $a$ be the radial distance from the origin to the free boundary. We will enforce $u(2) = 0$ as a radial homogenous Dirichlet boundary condition. Now our reference problem simplifies to
\begin{mdframed}
  Find $u(r)$ where for $r \in [0, 2]$,
  \begin{align}  \label{eq:Sphere Obstacle}
    ru''(r) + u'(r) = 0, &\text{ for } r \in (a, 2)\\
    u(2) = 0, \quad  &u(a) = \psi(a), \quad u'(a) = \psi'(a)
  \end{align}
\end{mdframed}
This ODE can be solved analytically and the solution has the following form, 
\begin{equation}
  u(r) = -A\log(r) + B \quad \text{ for } \quad a \leq r \leq 2.
\end{equation}
Using the boundary conditions (3.19) we have a system of equations to solve for $A$, $B$, and $a$. Doing so to high accuracy by Newton's method we get, 
 \begin{equation}
  a = 0.697965148223374, \quad A = 0.680259411891719, \quad B = 0.471519893402112.
 \end{equation}
 To find $g_D$ we use the solution values along $\partial \Omega$. The solution and a radial cross section of the solution are plotted below.
 \begin{figure}[H]
  \begin{center}
    \begin{tikzpicture}
      \begin{axis}[
        axis lines = left,
        xlabel = r,
        ylabel = {z},
        ylabel style={rotate = -90}, 
        xmin=0, xmax=2,
        ymin=0, ymax=1.5,
        no markers,
        samples=100,
        ]
        \addplot [domain=0.697965148223374:2, samples=100, smooth] {-0.680259411891719*ln(x) + 0.471519893402112};
  \addplot [domain=0:0.697965148223374, samples=100, smooth] {sqrt(1 - x^2)};
  \addplot [blue, domain=0:.9, samples=100, smooth] {sqrt(1 - x^2)}; 
  \addplot [red, domain=0:1.5] ({0.697965148223374},x); 
  \node at (axis cs:1,1) [red] {a}; 
  \node at (axis cs:.35,1.25) [black] {$A_u$}; 
  \node at (axis cs:1.5, .5) [black] {$I_u$}; 
  \node at (axis cs:1,.25) [blue] {$\psi(r)$}; 
\end{axis}
\end{tikzpicture}
  \caption{Radial cross section of $u(r)$.}
\end{center}
\end{figure}
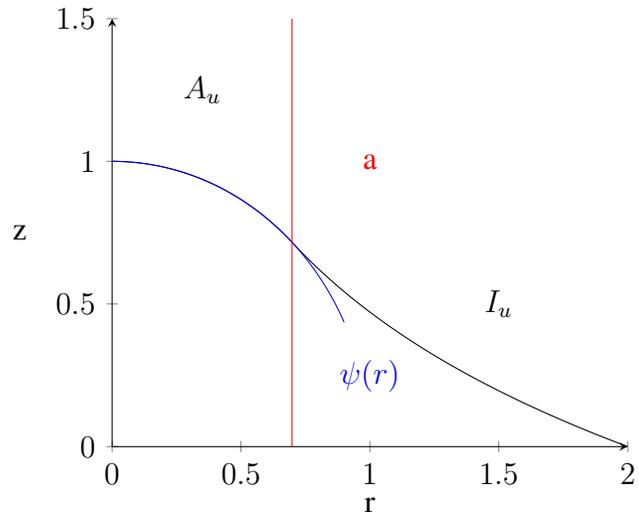
\begin{figure}[H]
  \includegraphics[width=\textwidth]{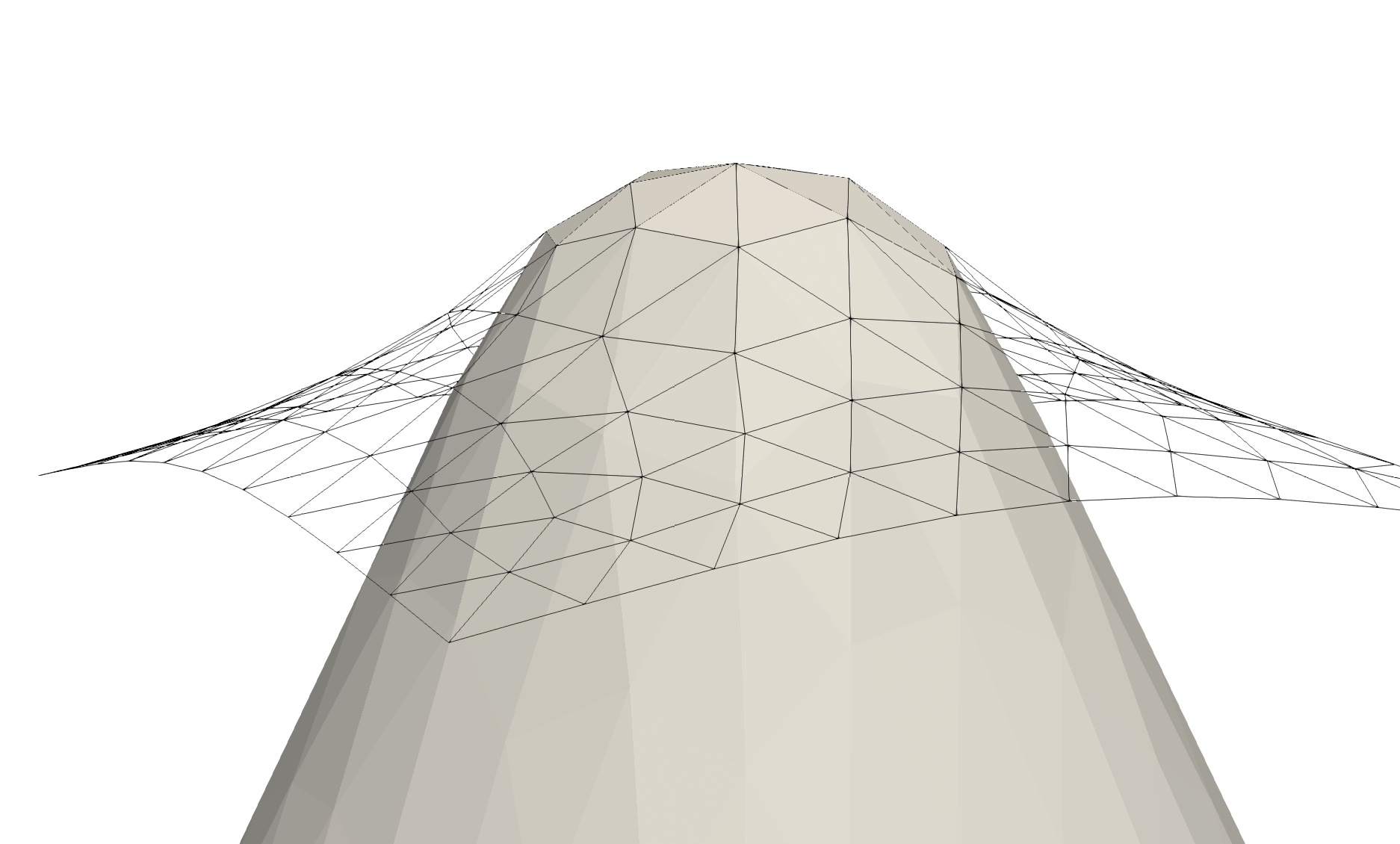}
  \caption{Reference Problem: $\psi$(gray) and $u$(black)}
\end{figure}

\section{Numerical Methods for Variational Inequalities}
\subsection{Active-Set Newton's Method}
We will be using the VINEWTONRSLS (VI-adapted Newton solver with reduced space line search) solver functionality in Firedrake, from the PETSc library \citep{balay_petsc_2015}, to solve our reference problem. We will only briefly outline the algorithm and key ideas. VINEWTONRSLS was introduced in \citet{benson_flexible_2003} and it solves finite dimensional nonlinear complementarity problems (NCP) of the form
\begin{equation}
  F(w) \geq 0, \quad w \geq 0, \quad F(w)w = 0, 
\end{equation}
where $w \in \RR^n$ and $F: \RR^n \to \RR^n$. Note that a finite-dimensional variational inequality can be formulated as a complementarity problem \citep{facchinei_finite-dimensional_2003}(Theorem 6). 

Here is how we will use the finite dimensional NCP form (3.20) to approximate the continuum CP (3.4). We will use the finite element method to discretize the obstacle $\psi$, the solution $u$ and the source term $f$. The $-\nabla^2$ operator will take the form of the stiffness matrix $A$ whose construction was described in Section 2.3. So we discretize via FEM and make the substitutions $w = u - \psi$ and $F(w) = -\nabla^2(w - \psi)-f$. Thus: 
\begin{equation}
  \underset{\text{Continuum CP}}{\begin{aligned}
    -\nabla^2u - f &\geq 0\\
    u - \psi &\geq 0\\
    (-\nabla^2u - f)(u - \psi) &= 0
  \end{aligned}}
  \quad
  \underset{\shortstack{\scriptsize \text{Discretize}\\ \scriptsize \text{\&}\\ \scriptsize \text{Substitute}}}{\implies}
  \quad
  \underset{\text{Finite Dimensional NCP}}{\begin{aligned}
    F(w) &\geq 0\\
    w &\geq 0\\
    F(w)w &= 0
  \end{aligned}}.
\end{equation}

Because of the constraints, this problem is inherently non-linear and we must use an iterative solver. The Newton algorithm below assumes that the solution $u$ is \emph{nondegenerate}, meaning that $-\nabla^2 u  - f > 0$ everywhere in the interior of the active set $A_u$. Since we assume nondegeneracy, the residual $F(w): \RR^n \to \RR^n$ will be positive on the active set. A revised residual function $\hat{F}: \RR^n \to \RR^n$ is defined which has the property that $\hat{F}(w) = 0$ when $w$ is the solution:
\begin{equation}
  \hat{F}_i(w) = \begin{cases}
    F_i(w) & \text{if } w_i > 0,\\
    \min\{F_i(w), 0\} & \text{if } w_i = 0.
  \end{cases}
\end{equation}
The algorithm proceeds similarly to any other Newton method with the exception of two key steps when computing the search direction and conducting the line search.
\begin{algorithm}
  \caption{VINEWTONRSLS}
  \begin{algorithmic}
    \Require $w^k \in \RR^n$.
    \State Define $A(w^k)$ ``certainly active" and $I(w^k)$ ``undecided" vertex sets.
    \State Compute $d$ via a reduced space Newton step on $I(w^k)$, (3.27) below. 
    \State Perform a line search along $\pi(w^k + \beta d)$ to stay inside $K = \{w \in \RR^n: w \geq 0\}$.
    \State Update $w^{k+1} = \pi(w^k + \beta d)$.
  \end{algorithmic}
  \end{algorithm}

Algorithm 1 begins by identifying index sets $A(w)$ and $I(w)$:
\begin{align}
  A(w^k) &= \{i \in \{1, ..., N\}| w_i^k = 0 \text{ and } F_i(w^k) > 0\},\\
  I(w^k) &= \{i \in \{1, ..., N\}| w_i^k > 0 \text{ or } F_i(w^k) \leq 0\}.
\end{align}

A search direction $d$ is computed by setting $d_{A_k} = 0$ and solving the following reduced system for the rest $d_{I_k}$,
\begin{equation}
  J(w^k)_{I^k, I^k}d_{I_k} = -F(w^k)_{I^k}.
\end{equation}
Here $J(w)$ is the Jacobian of $F(w)$ and the subscripts denote degrees of freedom relating to nodes in $I^k = I(w^k)$. 

The search direction is $d = [d_{A_k}, d_{I_k}]$. With a search direction $d$ it is possible that $w^k + d$ is inadmissable, for example an inactive node could be pushed underneath the obstacle. To remedy this we perform a projected line search. If $K = \{w \in \RR^n: w \geq 0\}$ is the admissable set of (3.20), and $w^k \in K$ we define $\pi: \RR^n \to K$ as follows, 
\begin{equation}
  \pi(w)_i = \begin{cases}
    w_i & \text{if } w_i > 0,\\
    0 & \text{if } w_i \leq 0.
  \end{cases}
\end{equation}
Then an Armijo line search \citep{armijo_minimization_1966} is conducted via a 2-norm merit function. The PETSc implementation of the line search uses parameters $\sigma = 10^{-4}$, $\beta = \frac{1}{2}$ and a minimum step size of $\gamma = 10^{-12}$ \citep[page 8]{benson_flexible_2003},
\begin{equation}
  \norm{\hat{F}(\pi(w^k + \beta d^k))}{2} \leq (1 - \sigma \beta)\norm{\hat{F}(w^k)}{2}.
\end{equation}

\subsection{Firedrake Implementation: Ball Obstacle}
We will implement the ball obstacle reference problem from Section 3.3 using Firedrake, and solve it using the VINEWTONRSLS solver. The code is very similar to a nonlinear PDE solver. 

First we import the Firedrake library, then use the mesh constructor to create a 2D mesh.
\begin{Verbatim}[frame = single,numbers = left,baselinestretch=0.75]
  from firedrake import *

  width, offset = 4.0, -2.0   # Omega = [-2,2]^2
  mesh = SquareMesh(20, 20, width, quadrilateral=False)
  mesh.coordinates.dat.data[:, :] += offset
\end{Verbatim}
Now we define the function space $V$, the obstacle $\psi$, the solution $u$, and the test function $v$. The obstacle $\psi$ is defined point-wise over the mesh using UFL. Note that $\psi$ is referenced as \texttt{lb} for ``lower bound". 
\begin{Verbatim}[frame = single,numbers = left,baselinestretch=0.75, firstnumber=last]
  V = FunctionSpace(mesh, "CG", 1)
  v = TestFunction(V)
  u = Function(V)

  # Define the obstacle function psi
  (x, y) = SpatialCoordinate(mesh)
  r = sqrt(x * x + y * y)
  r0 = 0.9
  psi0 = np.sqrt(1.0 - r0 * r0)
  dpsi0 = - r0 / psi0
  psi_ufl = conditional(le(r, r0), sqrt(1.0 - r * r),
                        psi0 + dpsi0 * (r - r0))
  lb = Function(V).interpolate(psi_ufl)
\end{Verbatim}
Recall that the exact solution is known and determines the Dirichlet boundary conditions. We define them both using UFL.
\begin{Verbatim}[frame = single,numbers = left,baselinestretch=0.75, firstnumber=last]
  afree = 0.697965148223374
  A = 0.680259411891719
  B = 0.471519893402112
  gbdry_ufl = conditional(le(r, afree), psi_ufl, - A * ln(r) + B)
  gbdry = Function(V).interpolate(gbdry_ufl)
  
  bdry_ids = (1, 2, 3, 4)   # boundary ids
  bcs = DirichletBC(V, gbdry, bdry_ids)
  
  uexact = gbdry.copy()
\end{Verbatim}
Finally we define the weak form. Previously we defined the left and right hand sides of our form separately and passed \texttt{a == L} to the solve command. To invoke the nonlinear solver we define the weak form by its residual form \texttt{F == 0} where $F = a - L$, but here the source term in our reference problem is zero.
\begin{Verbatim}[frame = single,numbers = left,baselinestretch=0.75, firstnumber=last]
  F = inner(grad(u), grad(v)) * dx
\end{Verbatim}

For the Poisson PDE problem we applied straightforward solver parameters to solve a linear Galerkin system using LU decomposition. The following example will be slightly more involved and will highlight the fact that Firedrake inherits solver parameters and logging functionality from PETSc. We set VINEWTONRSLS as the \texttt{snes\_type}. The zero tolerance for identifying the active and inactive sets is set with \texttt{snes\_vi\_zero\_tolerance}. The \texttt{sp} solver parameter dictionary allows the user to interface directly with PETSc options. We can define several other important parameters such as \texttt{snes\_atol},  \texttt{snes\_rtol}, and \texttt{snes\_stol} which are the absolute residual, relative iteration, and step tolerances respectively. Setting \texttt{snes\_vi\_monitor} prints $\abs{I(u^{k + 1})}/\abs{I(u^k)}$ at each iteration. The solver is then called with the bounds \texttt{lb} and \texttt{ub} set to the obstacle $\psi$ and infinity respectively.
\begin{Verbatim}[frame = single,numbers = left,baselinestretch=0.75, firstnumber=last]
  sp = {"snes_monitor": None,         
        "snes_type": "vinewtonrsls",
        "snes_converged_reason": None, 
        "snes_vi_monitor": None,      
        "snes_rtol": 1.0e-8,
        "snes_atol": 1.0e-12,
        "snes_stol": 1.0e-12,
        "snes_vi_zero_tolerance": 1.0e-12,
        # Newton step equations are solved by LU
        "ksp_type": "preonly",
        "pc_type": "lu",
        "pc_factor_mat_solver_type": "mumps"}

  problem = NonlinearVariationalProblem(F, u, bcs)
  solver = NonlinearVariationalSolver(
      problem, solver_parameters=sp)

  # it is necessary to define an upper obstacle
  ub = Function(V).interpolate(Constant(PETSc.INFINITY))  
  solver.solve(bounds=(lb, ub))  
\end{Verbatim}

Typical output from the solver looks like the following, 
{\small
\begin{Verbatim}[baselinestretch=0.75]
  0 SNES Function norm 2.046511757537e+00
  0 SNES VI Function norm 2.04651 
  Active lower constraints 17/29 
  upper constraints 0/0 
  Percent of total 0.125926 
  Percent of bounded 0.125926
  
  1 SNES Function norm 1.059615782871e-15
  1 SNES VI Function norm 1.05962e-15 
  Active lower constraints 17/17 
  upper constraints 0/0
  Percent of total 0.125926 
  Percent of bounded 0.125926
  
  Nonlinear  solve converged due to CONVERGED_FNORM_ABS iterations 1
\end{Verbatim}
}

\section{Convergence of Variational Inequalities}
As we will see, the convergence of numerical methods for variational inequalities is limited by the identification of the free boundary. In fact, the accuracy by which we locate the free boundary can dominate the other forms of numerical error. To demonstrate this consider the following one dimensional problem on $\Omega = [-1,1]$:
\begin{align}
  u \geq \psi\\
  u'' \geq 0\\
  (u - \psi)(u'') = 0\\
  u(-1) = u(1) &= 0\\
  \psi = 1/2 - &x^2
\end{align}
In this problem, the obstacle, domain, and Dirichlet conditions are even in $x$ and therefore the solution will be as well. Let $x = b$ be the right-most free boundary and note that since $u'' = 0$ on $(b, 1]$ it follows that $u$ is linear on $(b, 1]$:
\begin{equation}
  u(x) = (-\psi(b)/(1 - b))(x - 1)  \text{ on } (b, 1].
\end{equation}
Since $u'(b) = \psi'(b)$ we can solve the following equation for $b$, 
\begin{align}
  u'(b) &= \psi'(b),\\
  -\psi(b)/(1 - b) &= -2b,\\
  \frac{-(1/2 - b^2)}{1 - b} &= -2b,\\
  b &= \frac{2 \pm \sqrt{2}}{2}.
\end{align}
Note that only $b = \frac{2 - \sqrt{2}}{2} \in \Omega$, but by symmetry we know that $a = -b$ and therefore we have that the solution is piecewise defined as, 
\begin{equation}
  u = \begin{cases}
    \frac{(\frac{1}{2} - a^2)}{(a + 1)}(x + 1) & \text{ on } [-1, a]\\ 
    \frac{1}{2} - x^2 & \text{ on } (a, b)\\
    -\frac{(\frac{1}{2} - b^2)}{(1 - b)}(x - 1) & \text{ on } [b, 1]\\ 
  \end{cases}
\end{equation}

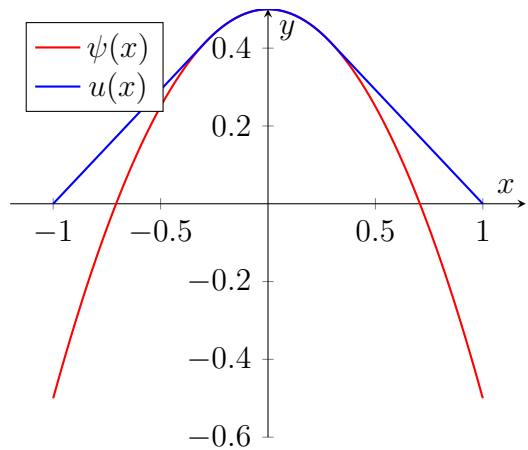
\begin{figure}[H]
  \centering
  \begin{tikzpicture}
    \begin{axis}[
      axis lines = middle,
      xlabel = {$x$},
      ylabel = {$y$},
      domain = -1:1,
      samples = 100,
      xmin = -1.2,
      xmax = 1.2,
      ymin = -0.6,
      ymax = 0.5,
      legend pos = north west
    ]
      \def\b{(2-sqrt(2))/2}
      \def\a{-\b}
      \addplot[
        thick,
        red
      ] {0.5 - x^2};
      \addlegendentry{$\psi(x)$}
  
      \addplot[
        thick,
        blue,
        domain=-1:\a
      ] {((.5 - \a*\a)/(\a + 1))*(x + 1)};

      \addplot[
        thick,
        blue,
        domain=\a:\b
      ] {0.5 - x^2};

      \addplot[
        thick,
        blue,
        domain=\b:1
      ] {-((.5 - \b*\b)/(1 - \b))*(x - 1)};
      \addlegendentry{$u(x)$}

    \end{axis}
  \end{tikzpicture}
\caption{The solution $u(x)$ and obstacle $\psi(x)$}
\end{figure}

Implemented in Firedrake, in Figure 3.4 we see that the error in the computed solution is actually proportional to the distance between the FE free boundary and the true free boundary; this distance is labeled as the ``gap" in the figure below. Therefore a sequence of uniform meshes would generate an $O(h)$ convergence rate with respect to the $L_2$ norm. Constructing an analagous PDE problem over $\Omega = [-1, 1]$:
\begin{align}
  u'' = 1\\
  u(-1) &= u(1) = 0
\end{align}
 Solving via FEM on the same sequence of meshes we achieve $O(h^2)$ convergence, as expected. 

\begin{figure}[H]
  \begin{center}
    \begin{subfigure}[b]{0.49\textwidth}
      \centering
      \includegraphics[width=\textwidth]{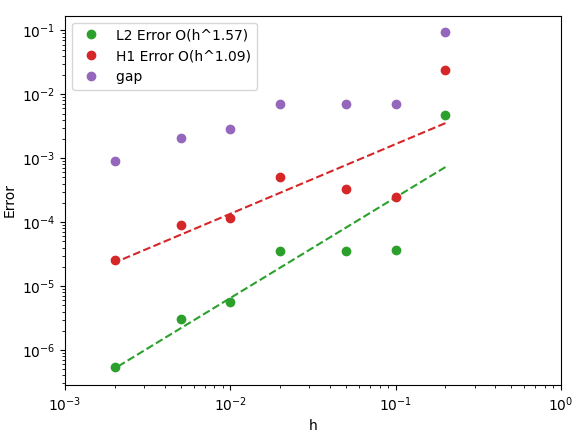}
      \label{fig:VI Reference Problem}
    \end{subfigure} 
    \begin{subfigure}[b]{0.49\textwidth}
      \centering
      \includegraphics[width=\textwidth]{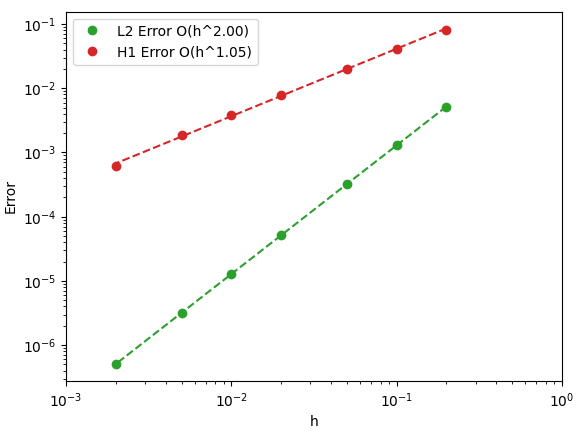}
      \label{fig:PDE Reference Problem}
    \end{subfigure}
    \caption{Convergence of comparable 1D VI (left) and PDE (right) problems.}
    \label{fig: Free Boundary Error}
  \end{center}
\end{figure}

Informally, the ``gap" in a VI problem is a geometric measure of the distance between the exact and FE free boundaries, or alternatively some measure for how well the FE free boundary approximates the exact free boundary. For the 1 dimensional VI problem described above the choice of measure was clear, $\abs{\Gamma_u - \Gamma_{\hat{u}}}$. For more complicated problems involving complex obstacles and higher dimensions there are several ways to measure the ``gap", which we will see in Section 6.1.

\chapter{Adaptive Mesh Refinement Techniques for PDE Problems}

To explain our methods for adaptive mesh refinement for VIs in the next chapter, we will first review certain methods for adaptive mesh refinement from the literature for PDEs.

\section{Tag and Refine Methods}

Tagging methods are a class of adaptive methods which assess the suitability of a mesh element-by-element in computing a Quantity of Interest (QoI) like $L_2$ error or another post computation functional like drag or lift \citep{bangerth_adaptive_2003}. The main idea behind most tagging methods is the refinement loop:

\begin{algorithm}
  \caption{Tagging Methods: Refinement Loop}
  \begin{algorithmic}
    \State Solve: Solve the PDE on the current mesh. 
    \State Estimate: Estimate the error in the QoI element-by-element.
    \State Tag: Tag elements for refinement based on the error estimate.
    \State Refine: Refine or coarsen the mesh maintaining the minimum angle criteria.
  \end{algorithmic}
\end{algorithm}

Throughout the literature (\citet{becker_feed-back_1996}, \citet{bangerth_adaptive_2003}, \citet{suttmeier_numerical_2008}) one finds a variety of ways to perform the ``Estimate" step. As we mentioned in the Introduction, a common way is by the Dual Weighted Residual (DWR) method introduced in \citet{becker_feed-back_1996}. For further details see \citet[Chapter 3]{bangerth_adaptive_2003}. A general approach for extending the DWR method to variational inequalities can be found in \citet{suttmeier_numerical_2008}. As we will see in Chapter 5 the methods proposed in this project are not based on the DWR method.

There are also several ways to perform the ``Tag" step. Consider the following fixed-rate strategy found in \citet[Chapter 4]{bangerth_adaptive_2003}. For fractions $X, Y$ with $1 - X > Y$ and a mesh with $N$ elements, refine the $X\cdot N$ elements with the largest error indicator and coarsen the $Y\cdot N$ elements with the smallest error indicator. For appropriate choices of $X$ and $Y$ this has the effect of keeping the degrees of freedom almost constant. There are other more exotic solutions which accomplish different goals. For example there is a fairly impractical ``error-balancing" strategy also described in \citet[Chapter 4]{bangerth_adaptive_2003} which seeks to equilibrate the error indicators across the mesh.

\section{Mesh Refinement}

The ``Refine" step in Algorithm 2 addresses the practical aspects of refinement of cells once they have been selected for refinement. This process involves two considerations: maintaining the minimum angle condition and managing hanging nodes. As illustrated in Figure 2.4, elements can be refined in such a way that the minimum angle condition is violated, leading to poor convergence properties. The second issue concerns hanging nodes, which are nodes that do not have a ``covering" relation to all neighboring elements; this concept will be elaborated on further in Chapter 5. The following Figure illustrates how a hanging node arises from refining a mesh. In the Figure, the refinement is shown in green and the hanging node is in red.
\begin{figure}[H]
  \centering
  \definecolor{ffqqqq}{rgb}{1.,0.,0.}
\definecolor{uuuuuu}{rgb}{0.26666666666666666,0.26666666666666666,0.26666666666666666}
\definecolor{xdxdff}{rgb}{0.49019607843137253,0.49019607843137253,1.}
\begin{tikzpicture}[line cap=round,line join=round,>=triangle 45,x=1.5cm,y=1.5cm]
\clip(2.5,-2.) rectangle (7.5,2.);
\draw [line width=2.pt] (3.,0.)-- (5.013325073747922,1.544475570087121);
\draw [line width=2.pt] (5.013325073747922,1.544475570087121)-- (5.013325073747922,-1.535428203717027);
\draw [line width=2.pt] (5.013325073747922,-1.535428203717027)-- (3.,0.);
\draw [line width=2.pt] (5.013325073747922,-1.535428203717027)-- (7.,0.);
\draw [line width=2.pt] (5.013325073747922,1.544475570087121)-- (7.,0.);

\draw [line width=2.pt, color = green] (4.006662536873961,0.7722377850435606)-- (4.006662536873961,-0.7677141018585135);
\draw [line width=2.pt, color = green] (5.013325073747922,0.004523683185047034)-- (4.006662536873961,-0.7677141018585135);
\draw [line width=2.pt, color = green] (4.006662536873961,0.7722377850435606)-- (5.013325073747922,0.004523683185047034);
\begin{scriptsize}
\draw [fill=xdxdff] (3.,0.) circle (2.5pt);
\draw [fill=xdxdff] (7.,0.) circle (2.5pt);
\draw [fill=xdxdff] (5.013325073747922,1.544475570087121) circle (2.5pt);
\draw [fill=xdxdff] (5.013325073747922,-1.535428203717027) circle (2.5pt);
\draw [fill=uuuuuu] (4.006662536873961,0.7722377850435606) circle (2.0pt);
\draw [fill=ffqqqq] (5.013325073747922,0.004523683185047034) circle (2.0pt);
\draw [fill=uuuuuu] (4.006662536873961,-0.7677141018585135) circle (2.0pt);
\end{scriptsize}
\end{tikzpicture}
\caption{Refinement (green) and hanging node(red).}
\end{figure}
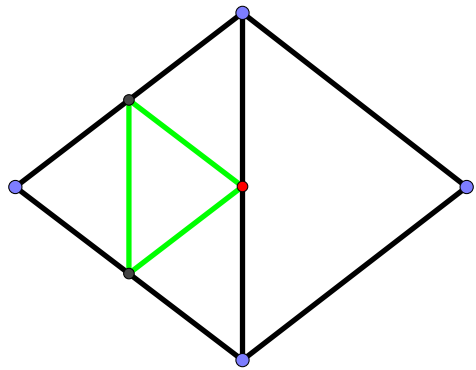
Hanging nodes do not increase the degrees of freedom of the resultant Galerkin system, but they do change the support of the basis functions. The following is a figure from \citet[Chapter 4]{bangerth_adaptive_2003} which illustrates how a hanging node affects the support of a basis function. 
\begin{figure}[H]
  \centering
  \includegraphics[width=0.8\textwidth]{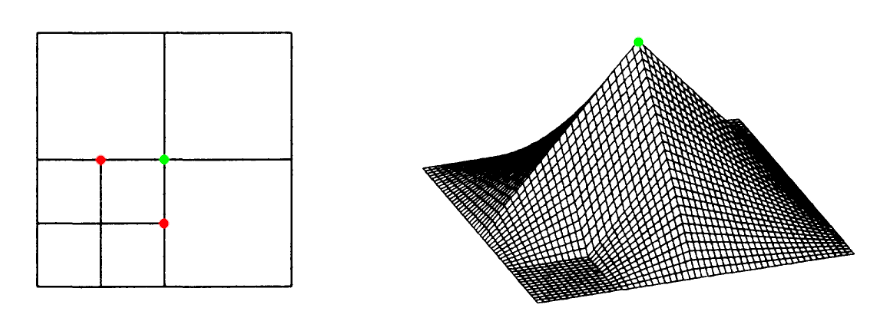}
  \caption{Hanging node (red) and a neighboring node (green) basis function.}
\end{figure}
There are several consequences of a refinement procedures which allow hanging nodes. Most immediately, they lead to greater complexity in the implementation of the solve step, as one needs to manage the different supports of the basis functions. Hanging nodes can be avoided altogether by using special transition cells to bridge edge lengths between refined and unrefined elements. This approach results in a refinement area that is larger than originally selected, but it is the approach we will use. 
\definecolor{wwccqq}{rgb}{0.4,0.8,0.}
\definecolor{ffqqqq}{rgb}{1.,0.,0.}
\definecolor{ffqqtt}{rgb}{1.,0.,0.2}
\definecolor{ududff}{rgb}{0.30196078431372547,0.30196078431372547,1.}
\captionsetup[subfigure]{justification=centering}
\begin{figure}[H]
  \centering
  
  \begin{subfigure}[b]{\textwidth}
    \centering
    \begin{tikzpicture}[scale = .90, line cap=round,line join=round,>=triangle 45,x=1.0cm,y=1.0cm]
      \clip(3.,0.53) rectangle (19.77,4.5);
      \fill[line width=2.pt,color=ffqqqq,fill=ffqqqq,fill opacity=0.10000000149011612] (7.380621075895451,3.9988352126590048) -- (9.071530166804537,1.0700947562243597) -- (10.76243925771362,3.9988352126589986) -- cycle;
      \fill[line width=2.pt,color=ffqqqq,fill=ffqqqq,fill opacity=0.10000000149011612] (9.077047122916897,-0.002898384860022918) -- (9.077047122916898,-2.9316388412946712) -- (10.767956213825988,-0.002898384860023917) -- cycle;
      \fill[line width=2.pt,color=wwccqq,fill=wwccqq,fill opacity=0.10000000149011612] (15.38328590011327,4.00142981467603) -- (13.692376809204186,1.0726893582413897) -- (17.07419499102237,1.0726893582413872) -- cycle;
      \fill[line width=2.pt,color=wwccqq,fill=wwccqq,fill opacity=0.10000000149011612] (17.07291427598191,-0.0036559327039541647) -- (16.22745973052736,-1.4680261609212581) -- (17.07291427598189,-2.9323963891385665) -- cycle;
      \draw [line width=2.pt] (12.000187003254672,-0.003655932703949632)-- (15.38200518507283,-0.003655932703949632);
      \draw [line width=2.pt] (12.000187003254672,-0.003655932703949632)-- (13.69109609416375,-2.932396389138574);
      \draw [line width=2.pt] (13.69109609416375,-2.932396389138574)-- (15.38200518507283,-0.003655932703949632);
      \draw [line width=2.pt] (17.07291427598189,-2.9323963891385665)-- (15.38200518507283,-0.003655932703949632);
      \draw [line width=2.pt] (13.69109609416375,-2.932396389138574)-- (17.07291427598189,-2.9323963891385665);
      \draw [line width=2.pt] (18.763823366890986,-0.0036559327039586975)-- (15.38200518507283,-0.003655932703949632);
      \draw [line width=2.pt] (17.07291427598189,-2.9323963891385665)-- (18.763823366890986,-0.0036559327039586975);
      \draw [line width=2.pt] (3.9988028940772704,3.9988352126590048)-- (7.380621075895451,3.9988352126590048);
      \draw [line width=2.pt] (3.9988028940772704,3.9988352126590048)-- (5.68971198498636,1.0700947562243581);
      \draw [line width=2.pt] (5.68971198498636,1.0700947562243581)-- (7.380621075895451,3.9988352126590048);
      \draw [line width=2.pt] (5.68971198498636,1.0700947562243581)-- (9.071530166804537,1.0700947562243597);
      \draw [line width=2.pt] (9.071530166804537,1.0700947562243597)-- (7.380621075895451,3.9988352126590048);
      \draw [line width=2.pt] (10.76243925771362,3.9988352126589986)-- (7.380621075895451,3.9988352126590048);
      \draw [line width=2.pt] (9.071530166804537,1.0700947562243597)-- (10.76243925771362,3.9988352126589986);
      \draw [line width=2.pt] (17.07291427598191,-0.0036559327039541647)-- (17.07291427598189,-2.9323963891385665);
      \draw [line width=2.pt] (17.07291427598191,-0.0036559327039541647)-- (16.22745973052736,-1.4680261609212581);
      \draw [line width=2.pt] (17.07291427598191,-0.0036559327039541647)-- (17.918368821436438,-1.4680261609212626);
      \draw [line width=2.pt] (15.38200518507283,-0.003655932703949632)-- (15.38200518507282,-2.93239638913857);
      \draw [line width=2.pt] (16.22745973052736,-1.4680261609212581)-- (15.38200518507282,-2.93239638913857);
      \draw [line width=2.pt] (12.001467718295102,4.00142981467603)-- (15.38328590011327,4.00142981467603);
      \draw [line width=2.pt] (12.001467718295102,4.00142981467603)-- (13.692376809204186,1.0726893582413897);
      \draw [line width=2.pt] (13.692376809204186,1.0726893582413897)-- (15.38328590011327,4.00142981467603);
      \draw [line width=2.pt] (17.07419499102237,1.0726893582413872)-- (15.38328590011327,4.00142981467603);
      \draw [line width=2.pt] (13.692376809204186,1.0726893582413897)-- (17.07419499102237,1.0726893582413872);
      \draw [line width=2.pt] (18.765104081931447,4.001429814676032)-- (15.38328590011327,4.00142981467603);
      \draw [line width=2.pt] (17.07419499102237,1.0726893582413872)-- (18.765104081931447,4.001429814676032);
      \draw [line width=2.pt,color=ffqqtt] (17.07419499102236,4.001429814676031)-- (17.07419499102237,1.0726893582413872);
      \draw [line width=2.pt,color=ffqqtt] (17.07419499102236,4.001429814676031)-- (16.22874044556782,2.5370595864587084);
      \draw [line width=2.pt,color=ffqqtt] (17.07419499102236,4.001429814676031)-- (17.91964953647691,2.5370595864587098);
      \draw [line width=2.pt,color=ffqqtt] (15.38328590011327,4.00142981467603)-- (15.383285900113279,1.0726893582413886);
      \draw [line width=2.pt,color=ffqqtt] (16.22874044556782,2.5370595864587084)-- (15.383285900113279,1.0726893582413886);
      \draw [line width=2.pt] (4.004319850189624,-0.0028983848600219186)-- (7.386138032007805,-0.0028983848600219186);
      \draw [line width=2.pt] (4.004319850189624,-0.0028983848600219186)-- (5.695228941098715,-2.931638841294669);
      \draw [line width=2.pt] (5.695228941098715,-2.931638841294669)-- (7.386138032007805,-0.0028983848600219186);
      \draw [line width=2.pt] (9.077047122916898,-2.9316388412946712)-- (7.386138032007805,-0.0028983848600219186);
      \draw [line width=2.pt] (5.695228941098715,-2.931638841294669)-- (9.077047122916898,-2.9316388412946712);
      \draw [line width=2.pt] (10.767956213825988,-0.002898384860023917)-- (7.386138032007805,-0.0028983848600219186);
      \draw [line width=2.pt] (9.077047122916898,-2.9316388412946712)-- (10.767956213825988,-0.002898384860023917);
      \draw [line width=2.pt] (9.077047122916897,-0.002898384860022918)-- (9.077047122916898,-2.9316388412946712);
      \draw [line width=2.pt] (9.077047122916897,-0.002898384860022918)-- (8.23159257746235,-1.4672686130773465);
      \draw [line width=2.pt] (9.077047122916897,-0.002898384860022918)-- (9.922501668371442,-1.4672686130773476);
      \draw [line width=2.pt] (7.386138032007805,-0.0028983848600219186)-- (7.386138032007807,-2.93163884129467);
      \draw [line width=2.pt] (8.23159257746235,-1.4672686130773465)-- (7.386138032007807,-2.93163884129467);
      \draw [line width=2.pt] (17.07291427598191,-0.0036559327039541647)-- (18.763823366890986,-0.0036559327039586975);
      \draw [line width=2.pt,color=ffqqtt] (17.918368821436445,-0.003655932703956431)-- (17.918368821436438,-1.4680261609212626);
      \draw [line width=2.pt] (17.918368821436438,-1.4680261609212626)-- (18.763823366890986,-0.0036559327039586975);
      \draw [line width=2.pt,color=ffqqtt] (17.495641548709173,-0.7358410468126084)-- (17.918368821436445,-0.003655932703956431);
      \draw [line width=2.pt,color=ffqqtt] (18.34109609416371,-0.7358410468126106)-- (17.918368821436445,-0.003655932703956431);
      \draw [line width=2.pt,color=ffqqtt] (17.0729142759819,-1.4680261609212604)-- (17.495641548709173,-0.7358410468126084);
      \draw [line width=2.pt,color=ffqqtt] (17.0729142759819,-1.4680261609212604)-- (17.918368821436438,-1.4680261609212626);
      \draw [line width=2.pt] (17.07291427598189,-2.9323963891385665)-- (17.918368821436438,-1.4680261609212626);
      \draw [line width=2.pt,color=ffqqtt] (17.0729142759819,-1.4680261609212604)-- (17.495641548709166,-2.2002112750299148);
      \draw [line width=2.pt,color=ffqqtt] (16.22745973052736,-1.4680261609212581)-- (17.0729142759819,-1.4680261609212604);
      \draw [line width=2.pt,color=ffqqqq] (7.380621075895451,3.9988352126590048)-- (9.071530166804537,1.0700947562243597);
      \draw [line width=2.pt,color=ffqqqq] (9.071530166804537,1.0700947562243597)-- (10.76243925771362,3.9988352126589986);
      \draw [line width=2.pt,color=ffqqqq] (10.76243925771362,3.9988352126589986)-- (7.380621075895451,3.9988352126590048);
      \draw [line width=2.pt,color=ffqqqq] (9.077047122916897,-0.002898384860022918)-- (9.077047122916898,-2.9316388412946712);
      \draw [line width=2.pt,color=ffqqqq] (9.077047122916898,-2.9316388412946712)-- (10.767956213825988,-0.002898384860023917);
      \draw [line width=2.pt,color=ffqqqq] (10.767956213825988,-0.002898384860023917)-- (9.077047122916897,-0.002898384860022918);
      \draw [line width=2.pt,color=wwccqq] (15.38328590011327,4.00142981467603)-- (13.692376809204186,1.0726893582413897);
      \draw [line width=2.pt,color=wwccqq] (13.692376809204186,1.0726893582413897)-- (17.07419499102237,1.0726893582413872);
      \draw [line width=2.pt,color=wwccqq] (17.07419499102237,1.0726893582413872)-- (15.38328590011327,4.00142981467603);
      \draw [line width=2.pt,color=wwccqq] (17.07291427598191,-0.0036559327039541647)-- (16.22745973052736,-1.4680261609212581);
      \draw [line width=2.pt,color=wwccqq] (16.22745973052736,-1.4680261609212581)-- (17.07291427598189,-2.9323963891385665);
      \draw [line width=2.pt,color=wwccqq] (17.07291427598189,-2.9323963891385665)-- (17.07291427598191,-0.0036559327039541647);
      \begin{scriptsize}
      \draw [fill=ududff] (12.000187003254672,-0.003655932703949632) circle (2.5pt);
      \draw [fill=ududff] (15.38200518507283,-0.003655932703949632) circle (2.5pt);
      \draw [fill=ududff] (13.69109609416375,-2.932396389138574) circle (2.5pt);
      \draw [fill=ududff] (13.69109609416375,-2.932396389138574) circle (2.5pt);
      \draw [fill=ududff] (15.38200518507283,-0.003655932703949632) circle (2.5pt);
      \draw [fill=ududff] (17.07291427598189,-2.9323963891385665) circle (2.5pt);
      \draw [fill=ududff] (17.07291427598189,-2.9323963891385665) circle (2.5pt);
      \draw [fill=ududff] (15.38200518507283,-0.003655932703949632) circle (2.5pt);
      \draw [fill=ududff] (18.763823366890986,-0.0036559327039586975) circle (2.5pt);
      \draw [fill=ududff] (3.9988028940772704,3.9988352126590048) circle (2.5pt);
      \draw [fill=ududff] (7.380621075895451,3.9988352126590048) circle (2.5pt);
      \draw [fill=ududff] (5.68971198498636,1.0700947562243581) circle (2.5pt);
      \draw [fill=ududff] (5.68971198498636,1.0700947562243581) circle (2.5pt);
      \draw [fill=ududff] (7.380621075895451,3.9988352126590048) circle (2.5pt);
      \draw [fill=ududff] (9.071530166804537,1.0700947562243597) circle (2.5pt);
      \draw [fill=ududff] (9.071530166804537,1.0700947562243597) circle (2.5pt);
      \draw [fill=ududff] (7.380621075895451,3.9988352126590048) circle (2.5pt);
      \draw [fill=ududff] (10.76243925771362,3.9988352126589986) circle (2.5pt);
      \draw [fill=ududff] (17.07291427598191,-0.0036559327039541647) circle (2.5pt);
      \draw [fill=ududff] (16.22745973052736,-1.4680261609212581) circle (2.5pt);
      \draw [fill=ududff] (17.918368821436438,-1.4680261609212626) circle (2.5pt);
      \draw [fill=ududff] (15.38200518507282,-2.93239638913857) circle (2.5pt);
      \draw [fill=ududff] (12.001467718295102,4.00142981467603) circle (2.5pt);
      \draw [fill=ududff] (15.38328590011327,4.00142981467603) circle (2.5pt);
      \draw [fill=ududff] (13.692376809204186,1.0726893582413897) circle (2.5pt);
      \draw [fill=ududff] (13.692376809204186,1.0726893582413897) circle (2.5pt);
      \draw [fill=ududff] (15.38328590011327,4.00142981467603) circle (2.5pt);
      \draw [fill=ududff] (17.07419499102237,1.0726893582413872) circle (2.5pt);
      \draw [fill=ududff] (17.07419499102237,1.0726893582413872) circle (2.5pt);
      \draw [fill=ududff] (15.38328590011327,4.00142981467603) circle (2.5pt);
      \draw [fill=ududff] (18.765104081931447,4.001429814676032) circle (2.5pt);
      \draw [fill=ududff] (17.07419499102236,4.001429814676031) circle (2.5pt);
      \draw [fill=ududff] (16.22874044556782,2.5370595864587084) circle (2.5pt);
      \draw [fill=ududff] (17.91964953647691,2.5370595864587098) circle (2.5pt);
      \draw [fill=ududff] (15.383285900113279,1.0726893582413886) circle (2.5pt);
      \draw [fill=ududff] (4.004319850189624,-0.0028983848600219186) circle (2.5pt);
      \draw [fill=ududff] (7.386138032007805,-0.0028983848600219186) circle (2.5pt);
      \draw [fill=ududff] (5.695228941098715,-2.931638841294669) circle (2.5pt);
      \draw [fill=ududff] (5.695228941098715,-2.931638841294669) circle (2.5pt);
      \draw [fill=ududff] (7.386138032007805,-0.0028983848600219186) circle (2.5pt);
      \draw [fill=ududff] (9.077047122916898,-2.9316388412946712) circle (2.5pt);
      \draw [fill=ududff] (9.077047122916898,-2.9316388412946712) circle (2.5pt);
      \draw [fill=ududff] (7.386138032007805,-0.0028983848600219186) circle (2.5pt);
      \draw [fill=ududff] (10.767956213825988,-0.002898384860023917) circle (2.5pt);
      \draw [fill=ududff] (9.077047122916897,-0.002898384860022918) circle (2.5pt);
      \draw [fill=ududff] (8.23159257746235,-1.4672686130773465) circle (2.5pt);
      \draw [fill=ududff] (9.922501668371442,-1.4672686130773476) circle (2.5pt);
      \draw [fill=ududff] (7.386138032007807,-2.93163884129467) circle (2.5pt);
      \draw [fill=ududff] (17.918368821436445,-0.003655932703956431) circle (2.5pt);
      \draw [fill=ududff] (17.495641548709173,-0.7358410468126084) circle (2.5pt);
      \draw [fill=ududff] (18.34109609416371,-0.7358410468126106) circle (2.5pt);
      \draw [fill=ududff] (17.0729142759819,-1.4680261609212604) circle (2.5pt);
      \draw [fill=ududff] (17.495641548709166,-2.2002112750299148) circle (2.5pt);
      \end{scriptsize}
      \end{tikzpicture}
      \caption{Example of transition element. Element marked for refinmenet in red. Transition element in green.}
      \label{fig:subfig1}
  \end{subfigure}

  \vspace{1em} 

  \begin{subfigure}[b]{\textwidth}
      \centering
      \begin{tikzpicture}[scale = .90, line cap=round,line join=round,>=triangle 45,x=1.0cm,y=1.0cm]
        \clip(3.,-3.5) rectangle (19.77,0.53);
        \fill[line width=2.pt,color=ffqqqq,fill=ffqqqq,fill opacity=0.10000000149011612] (7.380621075895451,3.9988352126590048) -- (9.071530166804537,1.0700947562243597) -- (10.76243925771362,3.9988352126589986) -- cycle;
        \fill[line width=2.pt,color=ffqqqq,fill=ffqqqq,fill opacity=0.10000000149011612] (9.077047122916897,-0.002898384860022918) -- (9.077047122916898,-2.9316388412946712) -- (10.767956213825988,-0.002898384860023917) -- cycle;
        \fill[line width=2.pt,color=wwccqq,fill=wwccqq,fill opacity=0.10000000149011612] (15.38328590011327,4.00142981467603) -- (13.692376809204186,1.0726893582413897) -- (17.07419499102237,1.0726893582413872) -- cycle;
        \fill[line width=2.pt,color=wwccqq,fill=wwccqq,fill opacity=0.10000000149011612] (17.07291427598191,-0.0036559327039541647) -- (16.22745973052736,-1.4680261609212581) -- (17.07291427598189,-2.9323963891385665) -- cycle;
        \draw [line width=2.pt] (12.000187003254672,-0.003655932703949632)-- (15.38200518507283,-0.003655932703949632);
        \draw [line width=2.pt] (12.000187003254672,-0.003655932703949632)-- (13.69109609416375,-2.932396389138574);
        \draw [line width=2.pt] (13.69109609416375,-2.932396389138574)-- (15.38200518507283,-0.003655932703949632);
        \draw [line width=2.pt] (17.07291427598189,-2.9323963891385665)-- (15.38200518507283,-0.003655932703949632);
        \draw [line width=2.pt] (13.69109609416375,-2.932396389138574)-- (17.07291427598189,-2.9323963891385665);
        \draw [line width=2.pt] (18.763823366890986,-0.0036559327039586975)-- (15.38200518507283,-0.003655932703949632);
        \draw [line width=2.pt] (17.07291427598189,-2.9323963891385665)-- (18.763823366890986,-0.0036559327039586975);
        \draw [line width=2.pt] (3.9988028940772704,3.9988352126590048)-- (7.380621075895451,3.9988352126590048);
        \draw [line width=2.pt] (3.9988028940772704,3.9988352126590048)-- (5.68971198498636,1.0700947562243581);
        \draw [line width=2.pt] (5.68971198498636,1.0700947562243581)-- (7.380621075895451,3.9988352126590048);
        \draw [line width=2.pt] (5.68971198498636,1.0700947562243581)-- (9.071530166804537,1.0700947562243597);
        \draw [line width=2.pt] (9.071530166804537,1.0700947562243597)-- (7.380621075895451,3.9988352126590048);
        \draw [line width=2.pt] (10.76243925771362,3.9988352126589986)-- (7.380621075895451,3.9988352126590048);
        \draw [line width=2.pt] (9.071530166804537,1.0700947562243597)-- (10.76243925771362,3.9988352126589986);
        \draw [line width=2.pt] (17.07291427598191,-0.0036559327039541647)-- (17.07291427598189,-2.9323963891385665);
        \draw [line width=2.pt] (17.07291427598191,-0.0036559327039541647)-- (16.22745973052736,-1.4680261609212581);
        \draw [line width=2.pt] (17.07291427598191,-0.0036559327039541647)-- (17.918368821436438,-1.4680261609212626);
        \draw [line width=2.pt] (15.38200518507283,-0.003655932703949632)-- (15.38200518507282,-2.93239638913857);
        \draw [line width=2.pt] (16.22745973052736,-1.4680261609212581)-- (15.38200518507282,-2.93239638913857);
        \draw [line width=2.pt] (12.001467718295102,4.00142981467603)-- (15.38328590011327,4.00142981467603);
        \draw [line width=2.pt] (12.001467718295102,4.00142981467603)-- (13.692376809204186,1.0726893582413897);
        \draw [line width=2.pt] (13.692376809204186,1.0726893582413897)-- (15.38328590011327,4.00142981467603);
        \draw [line width=2.pt] (17.07419499102237,1.0726893582413872)-- (15.38328590011327,4.00142981467603);
        \draw [line width=2.pt] (13.692376809204186,1.0726893582413897)-- (17.07419499102237,1.0726893582413872);
        \draw [line width=2.pt] (18.765104081931447,4.001429814676032)-- (15.38328590011327,4.00142981467603);
        \draw [line width=2.pt] (17.07419499102237,1.0726893582413872)-- (18.765104081931447,4.001429814676032);
        \draw [line width=2.pt,color=ffqqtt] (17.07419499102236,4.001429814676031)-- (17.07419499102237,1.0726893582413872);
        \draw [line width=2.pt,color=ffqqtt] (17.07419499102236,4.001429814676031)-- (16.22874044556782,2.5370595864587084);
        \draw [line width=2.pt,color=ffqqtt] (17.07419499102236,4.001429814676031)-- (17.91964953647691,2.5370595864587098);
        \draw [line width=2.pt,color=ffqqtt] (15.38328590011327,4.00142981467603)-- (15.383285900113279,1.0726893582413886);
        \draw [line width=2.pt,color=ffqqtt] (16.22874044556782,2.5370595864587084)-- (15.383285900113279,1.0726893582413886);
        \draw [line width=2.pt] (4.004319850189624,-0.0028983848600219186)-- (7.386138032007805,-0.0028983848600219186);
        \draw [line width=2.pt] (4.004319850189624,-0.0028983848600219186)-- (5.695228941098715,-2.931638841294669);
        \draw [line width=2.pt] (5.695228941098715,-2.931638841294669)-- (7.386138032007805,-0.0028983848600219186);
        \draw [line width=2.pt] (9.077047122916898,-2.9316388412946712)-- (7.386138032007805,-0.0028983848600219186);
        \draw [line width=2.pt] (5.695228941098715,-2.931638841294669)-- (9.077047122916898,-2.9316388412946712);
        \draw [line width=2.pt] (10.767956213825988,-0.002898384860023917)-- (7.386138032007805,-0.0028983848600219186);
        \draw [line width=2.pt] (9.077047122916898,-2.9316388412946712)-- (10.767956213825988,-0.002898384860023917);
        \draw [line width=2.pt] (9.077047122916897,-0.002898384860022918)-- (9.077047122916898,-2.9316388412946712);
        \draw [line width=2.pt] (9.077047122916897,-0.002898384860022918)-- (8.23159257746235,-1.4672686130773465);
        \draw [line width=2.pt] (9.077047122916897,-0.002898384860022918)-- (9.922501668371442,-1.4672686130773476);
        \draw [line width=2.pt] (7.386138032007805,-0.0028983848600219186)-- (7.386138032007807,-2.93163884129467);
        \draw [line width=2.pt] (8.23159257746235,-1.4672686130773465)-- (7.386138032007807,-2.93163884129467);
        \draw [line width=2.pt] (17.07291427598191,-0.0036559327039541647)-- (18.763823366890986,-0.0036559327039586975);
        \draw [line width=2.pt,color=ffqqtt] (17.918368821436445,-0.003655932703956431)-- (17.918368821436438,-1.4680261609212626);
        \draw [line width=2.pt] (17.918368821436438,-1.4680261609212626)-- (18.763823366890986,-0.0036559327039586975);
        \draw [line width=2.pt,color=ffqqtt] (17.495641548709173,-0.7358410468126084)-- (17.918368821436445,-0.003655932703956431);
        \draw [line width=2.pt,color=ffqqtt] (18.34109609416371,-0.7358410468126106)-- (17.918368821436445,-0.003655932703956431);
        \draw [line width=2.pt,color=ffqqtt] (17.0729142759819,-1.4680261609212604)-- (17.495641548709173,-0.7358410468126084);
        \draw [line width=2.pt,color=ffqqtt] (17.0729142759819,-1.4680261609212604)-- (17.918368821436438,-1.4680261609212626);
        \draw [line width=2.pt] (17.07291427598189,-2.9323963891385665)-- (17.918368821436438,-1.4680261609212626);
        \draw [line width=2.pt,color=ffqqtt] (17.0729142759819,-1.4680261609212604)-- (17.495641548709166,-2.2002112750299148);
        \draw [line width=2.pt,color=ffqqtt] (16.22745973052736,-1.4680261609212581)-- (17.0729142759819,-1.4680261609212604);
        \draw [line width=2.pt,color=ffqqqq] (7.380621075895451,3.9988352126590048)-- (9.071530166804537,1.0700947562243597);
        \draw [line width=2.pt,color=ffqqqq] (9.071530166804537,1.0700947562243597)-- (10.76243925771362,3.9988352126589986);
        \draw [line width=2.pt,color=ffqqqq] (10.76243925771362,3.9988352126589986)-- (7.380621075895451,3.9988352126590048);
        \draw [line width=2.pt,color=ffqqqq] (9.077047122916897,-0.002898384860022918)-- (9.077047122916898,-2.9316388412946712);
        \draw [line width=2.pt,color=ffqqqq] (9.077047122916898,-2.9316388412946712)-- (10.767956213825988,-0.002898384860023917);
        \draw [line width=2.pt,color=ffqqqq] (10.767956213825988,-0.002898384860023917)-- (9.077047122916897,-0.002898384860022918);
        \draw [line width=2.pt,color=wwccqq] (15.38328590011327,4.00142981467603)-- (13.692376809204186,1.0726893582413897);
        \draw [line width=2.pt,color=wwccqq] (13.692376809204186,1.0726893582413897)-- (17.07419499102237,1.0726893582413872);
        \draw [line width=2.pt,color=wwccqq] (17.07419499102237,1.0726893582413872)-- (15.38328590011327,4.00142981467603);
        \draw [line width=2.pt,color=wwccqq] (17.07291427598191,-0.0036559327039541647)-- (16.22745973052736,-1.4680261609212581);
        \draw [line width=2.pt,color=wwccqq] (16.22745973052736,-1.4680261609212581)-- (17.07291427598189,-2.9323963891385665);
        \draw [line width=2.pt,color=wwccqq] (17.07291427598189,-2.9323963891385665)-- (17.07291427598191,-0.0036559327039541647);
        \begin{scriptsize}
        \draw [fill=ududff] (12.000187003254672,-0.003655932703949632) circle (2.5pt);
        \draw [fill=ududff] (15.38200518507283,-0.003655932703949632) circle (2.5pt);
        \draw [fill=ududff] (13.69109609416375,-2.932396389138574) circle (2.5pt);
        \draw [fill=ududff] (13.69109609416375,-2.932396389138574) circle (2.5pt);
        \draw [fill=ududff] (15.38200518507283,-0.003655932703949632) circle (2.5pt);
        \draw [fill=ududff] (17.07291427598189,-2.9323963891385665) circle (2.5pt);
        \draw [fill=ududff] (17.07291427598189,-2.9323963891385665) circle (2.5pt);
        \draw [fill=ududff] (15.38200518507283,-0.003655932703949632) circle (2.5pt);
        \draw [fill=ududff] (18.763823366890986,-0.0036559327039586975) circle (2.5pt);
        \draw [fill=ududff] (3.9988028940772704,3.9988352126590048) circle (2.5pt);
        \draw [fill=ududff] (7.380621075895451,3.9988352126590048) circle (2.5pt);
        \draw [fill=ududff] (5.68971198498636,1.0700947562243581) circle (2.5pt);
        \draw [fill=ududff] (5.68971198498636,1.0700947562243581) circle (2.5pt);
        \draw [fill=ududff] (7.380621075895451,3.9988352126590048) circle (2.5pt);
        \draw [fill=ududff] (9.071530166804537,1.0700947562243597) circle (2.5pt);
        \draw [fill=ududff] (9.071530166804537,1.0700947562243597) circle (2.5pt);
        \draw [fill=ududff] (7.380621075895451,3.9988352126590048) circle (2.5pt);
        \draw [fill=ududff] (10.76243925771362,3.9988352126589986) circle (2.5pt);
        \draw [fill=ududff] (17.07291427598191,-0.0036559327039541647) circle (2.5pt);
        \draw [fill=ududff] (16.22745973052736,-1.4680261609212581) circle (2.5pt);
        \draw [fill=ududff] (17.918368821436438,-1.4680261609212626) circle (2.5pt);
        \draw [fill=ududff] (15.38200518507282,-2.93239638913857) circle (2.5pt);
        \draw [fill=ududff] (12.001467718295102,4.00142981467603) circle (2.5pt);
        \draw [fill=ududff] (15.38328590011327,4.00142981467603) circle (2.5pt);
        \draw [fill=ududff] (13.692376809204186,1.0726893582413897) circle (2.5pt);
        \draw [fill=ududff] (13.692376809204186,1.0726893582413897) circle (2.5pt);
        \draw [fill=ududff] (15.38328590011327,4.00142981467603) circle (2.5pt);
        \draw [fill=ududff] (17.07419499102237,1.0726893582413872) circle (2.5pt);
        \draw [fill=ududff] (17.07419499102237,1.0726893582413872) circle (2.5pt);
        \draw [fill=ududff] (15.38328590011327,4.00142981467603) circle (2.5pt);
        \draw [fill=ududff] (18.765104081931447,4.001429814676032) circle (2.5pt);
        \draw [fill=ududff] (17.07419499102236,4.001429814676031) circle (2.5pt);
        \draw [fill=ududff] (16.22874044556782,2.5370595864587084) circle (2.5pt);
        \draw [fill=ududff] (17.91964953647691,2.5370595864587098) circle (2.5pt);
        \draw [fill=ududff] (15.383285900113279,1.0726893582413886) circle (2.5pt);
        \draw [fill=ududff] (4.004319850189624,-0.0028983848600219186) circle (2.5pt);
        \draw [fill=ududff] (7.386138032007805,-0.0028983848600219186) circle (2.5pt);
        \draw [fill=ududff] (5.695228941098715,-2.931638841294669) circle (2.5pt);
        \draw [fill=ududff] (5.695228941098715,-2.931638841294669) circle (2.5pt);
        \draw [fill=ududff] (7.386138032007805,-0.0028983848600219186) circle (2.5pt);
        \draw [fill=ududff] (9.077047122916898,-2.9316388412946712) circle (2.5pt);
        \draw [fill=ududff] (9.077047122916898,-2.9316388412946712) circle (2.5pt);
        \draw [fill=ududff] (7.386138032007805,-0.0028983848600219186) circle (2.5pt);
        \draw [fill=ududff] (10.767956213825988,-0.002898384860023917) circle (2.5pt);
        \draw [fill=ududff] (9.077047122916897,-0.002898384860022918) circle (2.5pt);
        \draw [fill=ududff] (8.23159257746235,-1.4672686130773465) circle (2.5pt);
        \draw [fill=ududff] (9.922501668371442,-1.4672686130773476) circle (2.5pt);
        \draw [fill=ududff] (7.386138032007807,-2.93163884129467) circle (2.5pt);
        \draw [fill=ududff] (17.918368821436445,-0.003655932703956431) circle (2.5pt);
        \draw [fill=ududff] (17.495641548709173,-0.7358410468126084) circle (2.5pt);
        \draw [fill=ududff] (18.34109609416371,-0.7358410468126106) circle (2.5pt);
        \draw [fill=ududff] (17.0729142759819,-1.4680261609212604) circle (2.5pt);
        \draw [fill=ududff] (17.495641548709166,-2.2002112750299148) circle (2.5pt);
        \end{scriptsize}
        \end{tikzpicture}
      \caption{Further example showing how transition elements are refined.}
      \label{fig:subfig2}
  \end{subfigure}
  \vspace*{.25cm}
  \caption{Transition cell examples.}
  \label{fig:mainfig}
\end{figure}

\section{Firedrake-Netgen Integration}
A tagging strategy like the one above can be implemented in Firedrake in a variety of ways. In this project we used the Netgen/NGSolve integration \citep{zerbinati_ngspetsc_nodate}. This integration brings several new features to Firedrake. For our purposes the most important is the \texttt{.refine\_marked\_elements()} method. This method resides inside of a netgen mesh object and takes an indicator function over the domain representing which elements are marked for refinement. This method is capable of dealing with hanging nodes by use of transition cells, following the refinement pattern illustrated in Figure 4.3.

\section{Interpolation into DG0}
A DG0 finite element space is one in which the basis functions are constant with a value of 1 and are each supported over a single element.
\begin{figure}[H]
  \centering
  \includegraphics[width=0.6\textwidth]{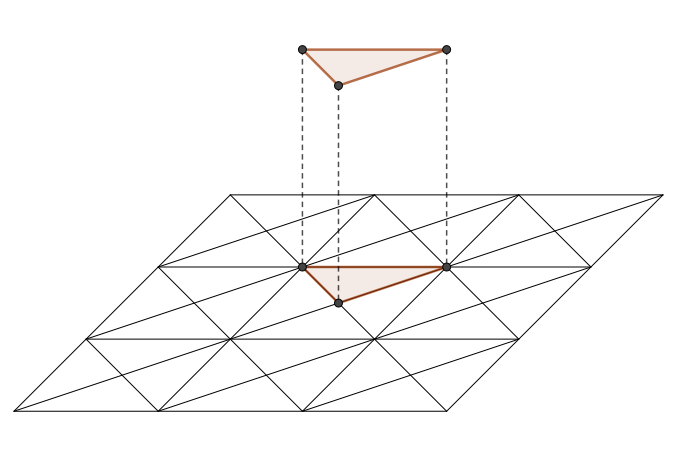}
  \caption{DG0 basis function.}
\end{figure}
The operation of interpolating a function into a DG0 space is central to our proposed methods for adaptive refinement. Conceptually this operation allows us to compute estimators (Estimate step of Algorithim 2) which take on a single value over vertices and then convert them, by averaging over an element into an estimator which takes on a single value over the element. 

Averaging a function over elements is equivalent to interpolation into a DG0 space. To see this, consider the definition of the interpolation operation given in the Firedrake user manual \citep{ham_firedrake_2023}. 
Let $u$ be a function over some domain $\Omega$ and $V$ be a finite element space defined on $\Omega$ with basis $\{\phi_i\}_{i = 1}^N$. Then $interpolate(u, V)$ is given by $\sum_{i = 1}^N v_i \phi_i$ where $v_i = \phi^*_i$ and $\phi^*_i$ is an element of the dual basis of $V$. (Recall that the dual basis of $V_\Omega$ is given by the vector $\{\phi^*_i\}$, such that $\phi^*_i(\phi_j) = \delta_{ij}$). The corresponding dual basis vector of a DG0 space is the following average
\begin{equation}
  \phi_j^*(f) = \frac{1}{area(\triangle_j)}\int_{\triangle_j} f\, dx.
\end{equation}
Note that this choice of functional has the property that $\phi^*_i(\phi_j) = \delta_{ij}$. Now $interpolate(u, V) = \sum_{i = 1}^n v_i\theta_i$ where
\begin{equation}
  v_i = \frac{1}{area(\triangle_i)}\int_{\triangle_i} u \, dx.
\end{equation}

\chapter{New Adaptive Refinement Schemes for VIs}

We propose two strategies for adaptive refinement of variational inequalities, which we call Variable Coefficient Elliptic Smoothing (VCES) and Unstructured Dilation Operator (UDO). These methods are designed to enhance mesh resolution at the free boundary found when solving the VI. These methods are the Estimate and Tag steps in Algorithm 2. They effectively control the error associated with approximating the free boundary, as observed in Section 3.3 and illustrated in Figure 3.4. Note that Firedrake/Netgen is used for the Solve and Refine steps in Algorithim 2.

\section{Variable Coefficient Elliptic Smoothing}
The first strategy is called Variable Coefficient Elliptic Smoothing (VCES). The idea is to use the residual $u^k - \psi$ and a positive tolerance to construct $s_0$, a node-wise indicator function for the active set. This indicator function is used as the initial condition of a variable-coefficient and time-dependent heat equation problem which has the effect of smoothing this indicator near the free boundary. This heat equation problem is solved for a single timestep via implicit Euler, the result of which we'll call $s_1$
\begin{equation}
  \frac{1}{\Delta t}(s_1 - s_0) = \nabla^2 s_1.
\end{equation}
This problem is a linear elliptic PDE, hence ``elliptic smoothing". Our choice of timestep is $\Delta t_i = \frac{1}{2}(\text{avg}(\text{diam}(\triangle_i)))^2$ where $\triangle_i$ is the set of all elements incident to vertex $i$. This choice depends on the element, thus ``variable coefficient". Varying the timestep based off of neighboring element diameter has the effect of applying the same amount diffusion across all elements regardless of size. The result is then interpolated into a DG0 space and thresholded to produce the refinement indicator, as in Algorithim 2. There are various parameters to consider with this technique. The choice of timestep and thresholding parameters will substantially affect the ``distance" about which the free boundary is resolved. 

\begin{algorithm}[H]
	\caption{Variable Coefficient Elliptic Smoothing Element Tagging for VIs}\label{alg:cap}
	\begin{algorithmic}[1]
		\Require $tol \in \mathbb{R}$, $u^k \in K_\psi, \psi \in V$, $W$ is DG0 FE space.
		\Require Threshold parameters $0\leq \alpha < \beta \leq 1$.
		\State Compute the nodal active set indicator function $s_0$
		  \begin{equation*}
			s_0 = \begin{cases}
			  1 & \text{ if } u^k - \psi < tol\\
			  0 & \text{ otherwise}
			\end{cases}
		  \end{equation*}
		\State Let $\Delta t_i = \frac{1}{2}(\text{avg}(\text{diam}(\triangle_i)))^2$, a CG1 field.
	  
		\State Solve $\frac{1}{\Delta t}(s_1 - s_0) = \nabla^2 s_1$ with $g_D = s_0|_{\partial\Omega}$ impliclty with Firedrake defaults settings. 
		\State Let $s_W = interpolate(s_1, W)$.
		\State Define the refinement indicator $I \in W$ as follows:
		\begin{equation*}
		  I(\triangle) = \begin{cases}
			1 & \text{ if } \alpha < s_W(\triangle) < \beta\\
			0 & \text{ otherwise}
		  \end{cases}
		\end{equation*}\\
		\Return $I$
	\end{algorithmic}
	\end{algorithm}

  Figure 5.1 illustrates the VCES algorithm applied to a problem similar to that in Section 3.3. 
\begin{figure}[H]
  \captionsetup[subfigure]{justification=centering}
  \centering
  \null\hfill
  \begin{subfigure}[b]{.35\textwidth}
    \begin{tikzpicture}[line cap=round,line join=round,>=triangle 45,scale=1.5]
      \clip(-2.25,-0.5) rectangle (2.25,1.25);
      \draw (0.4,0.9) -- (0.8,0.6);
      \draw (0.8,0.6) -- (1.2,0.4);
      \draw (1.2,0.4) -- (1.6,0.2);
      \draw (-0.4,0.9) -- (-0.8,0.6);
      \draw (-0.8,0.6) -- (-1.2,0.4);
      \draw (-1.2,0.4) -- (-1.6,0.2);
      \draw (-1.6,0.2) -- (-2,0);
      \draw (1.6,0.2) -- (2,0);
      \draw (-0.4,0.9) -- (0,1);
      \draw (0,1) -- (0.4,0.9);

      \draw[red] (1.6,-.2) -- (2,-.6);
      \draw[red] (1.2,0.2) -- (1.6,-.2);
      \draw[red] (0.8,0.6) -- (1.2,0.2);
      \draw[red] (0.4,0.9) -- (0.8,0.6);
      \draw[red] (0,1) -- (0.4,0.9);
      \draw[red] (-0.4,0.9) -- (0,1);
      \draw[red] (-0.4,0.9) -- (-0.8,0.6);
      \draw[red] (-0.8,0.6) -- (-1.2,.2);
      \draw[red] (-1.2, .2) -- (-1.6,-.2);
      \draw[red] (-1.6,-.2) -- (-2,-.6);

      \begin{scriptsize}
        \fill [color=black] (-0.4,0.9) circle (1.5pt);
        \fill [color=black] (0.4,0.9) circle (1.5pt);
        \fill [color=black] (0.8,0.6) circle (1.5pt);
        \fill [color=black] (1.2,0.4) circle (1.5pt);
        \fill [color=black] (1.6,0.2) circle (1.5pt);
        \fill [color=black] (-0.8,0.6) circle (1.5pt);
        \fill [color=black] (-1.2,0.4) circle (1.5pt);
        \fill [color=black] (-1.6,0.2) circle (1.5pt);
        \fill [color=black] (-2,0) circle (1.5pt);
        \fill [color=black] (2,0) circle (1.5pt);
        \fill [color=black] (0,1) circle (1.5pt);
      \end{scriptsize}

      \draw[-] (-2.5,0) -- (2.5,0); 
      \draw[-] (0,-0.5) -- (0,1.45); 
    \end{tikzpicture}
    \caption{Iterate $u^k$(black) and obstacle $\psi$(red).}
  \end{subfigure}
  \hfill
  \begin{subfigure}[b]{.35\textwidth}
    \begin{tikzpicture}[line cap=round,line join=round,>=triangle 45,scale=1.5]
      \clip(-2.25,-0.5) rectangle (2.25,1.25);
      \draw (1.6,0) -- (2,0);
      \draw (1.2,0) -- (1.6,0);
      \draw (0.8,1) -- (1.2,0);
      \draw (0.4,1) -- (0.8,1);
      \draw (0,1) -- (0.4,1);
      \draw (-0.4,1) -- (0,1);
      \draw (-0.4,1) -- (-0.8,1);
      \draw (-0.8,1) -- (-1.2,0);
      \draw (-1.2,0) -- (-1.6,0);
      \draw (-1.6,0) -- (-2,0);

      \begin{scriptsize}
        \fill [color=black] (1.6,0) circle (1.5pt);
        \fill [color=black] (2,0) circle (1.5pt);
        \fill [color=black] (1.2,0) circle (1.5pt);
        \fill [color=black] (0.8,1) circle (1.5pt);
        \fill [color=black] (1.2,0) circle (1.5pt);
        \fill [color=black] (0.4,1) circle (1.5pt);
        \fill [color=black] (0.8,1) circle (1.5pt);
        \fill [color=black] (0,1) circle (1.5pt);
        \fill [color=black] (0.4,1) circle (1.5pt);
        \fill [color=black] (-0.4,1) circle (1.5pt);
        \fill [color=black] (0,1) circle (1.5pt);
        \fill [color=black] (-0.4,1) circle (1.5pt);
        \fill [color=black] (-0.8,1) circle (1.5pt);
        \fill [color=black] (-1.2,0) circle (1.5pt);
        \fill [color=black] (-1.6,0) circle (1.5pt);
        \fill [color=black] (-2,0) circle (1.5pt);
      \end{scriptsize}
      
      \draw[-] (-2.5,0) -- (2.5,0); 
      \draw[-] (0,-0.5) -- (0,1.45); 
  
      \node at (-0.1,1.15) [left] {1};
      \draw[dashed] (-0.1,1) -- (0,1);
    \end{tikzpicture}
    \caption{Nodal active set indicator $s_0$.}
  \end{subfigure}
  \null\hfill
  \medskip

  \null\hfill
  \begin{subfigure}[b]{.35\textwidth}
    \centering
    \begin{tikzpicture}[line cap=round,line join=round,>=triangle 45,scale=1.5]
      \clip(-2.25,-0.5) rectangle (2.25,1.25);
      
      \draw (1.6, 0.1) -- (2,0);      
      \draw (1.2, 0.33) -- (1.6, 0.1);
      \draw (0.8, 0.66) -- (1.2, 0.33);
      \draw (0.4, 0.95) -- (0.8, 0.66);
      \draw (0,1) -- (0.4, 0.95);
      \draw (-0.4, 0.95) -- (0,1);
      \draw (-0.4, 0.95) -- (-0.8, 0.66);
      \draw (-0.8, 0.66) -- (-1.2, 0.33);
      \draw (-1.2, 0.33) -- (-1.6, 0.1);
      \draw (-1.6, 0.1) -- (-2,0);
  
      \begin{scriptsize}
        \fill [color=black] (1.6, 0.1) circle (1.5pt);
        \fill [color=black] (2,0) circle (1.5pt);
        \fill [color=black] (1.2, 0.33) circle (1.5pt);
        \fill [color=black] (0.8, 0.66) circle (1.5pt);
        \fill [color=black] (0.4, 0.95) circle (1.5pt);
        \fill [color=black] (0,1) circle (1.5pt);
        \fill [color=black] (-0.4, 0.95) circle (1.5pt);
        \fill [color=black] (-0.8, 0.66) circle (1.5pt);
        \fill [color=black] (-1.2, 0.33) circle (1.5pt);
        \fill [color=black] (-1.6, 0.1) circle (1.5pt);
        \fill [color=black] (-2,0) circle (1.5pt);
      \end{scriptsize}
  
      \draw[-] (-2.5,0) -- (2.5,0); 
      \draw[-] (0,-0.5) -- (0,1.45); 
  
      \node at (-0.1,1.15) [left] {1};
      \draw[dashed] (-0.1,1) -- (0,1);
    \end{tikzpicture}
    \caption{Smoothed $s_1$.}
  \end{subfigure}
  \hfill
  \begin{subfigure}[b]{.35\textwidth}
    \centering
    \begin{tikzpicture}[line cap=round,line join=round,>=triangle 45,scale=1.5]
      \clip(-2.25,-0.5) rectangle (2.25,1.25);

      \draw (-2, 0.05) -- (-1.6, 0.05);
      \draw (-1.6, 0.215) -- (-1.2, 0.215);
      \draw (-1.2, 0.495) -- (-0.8, 0.495);
      \draw (-0.8, 0.75) -- (-0.4, 0.75);
      \draw (-0.4, 0.9) -- (0, 0.9);
      \draw (0, 0.9) -- (0.4, 0.9);
      \draw (0.4, 0.75) -- (0.8, 0.75);
      \draw (0.8, 0.495) -- (1.2, 0.495);
      \draw (1.2, 0.215) -- (1.6, 0.215);
      \draw (1.6, 0.05) -- (2, 0.05);

      \begin{scriptsize}
        \fill [color=black] (-2, 0.05) circle (1pt);
        \fill [color=black] (-1.6, 0.05) circle (1pt);
        \fill [color=black] (-1.6, 0.215) circle (1pt);
        \fill [color=black] (-1.2, 0.215) circle (1pt);
        \fill [color=black] (-1.2, 0.495) circle (1pt);
        \fill [color=black] (-0.8, 0.495) circle (1pt);
        \fill [color=black] (-0.8, 0.75) circle (1pt);
        \fill [color=black] (-0.4, 0.75) circle (1pt);
        \fill [color=black] (-0.4, 0.9) circle (1pt);
        \fill [color=black] (0, 0.9) circle (1pt);
        \fill [color=black] (0, 0.9) circle (1pt);
        \fill [color=black] (0.4, 0.9) circle (1pt);
        \fill [color=black] (0.4, 0.75) circle (1pt);
        \fill [color=black] (0.8, 0.75) circle (1pt);
        \fill [color=black] (0.8, 0.495) circle (1pt);
        \fill [color=black] (1.2, 0.495) circle (1pt);
        \fill [color=black] (1.2, 0.215) circle (1pt);
        \fill [color=black] (1.6, 0.215) circle (1pt);
        \fill [color=black] (1.6, 0.05) circle (1pt);
        \fill [color=black] (2, 0.05) circle (1pt);
      \end{scriptsize}

      \draw[-] (-2.5,0) -- (2.5,0); 
      \draw[-] (0,-0.5) -- (0,1.25); 

      \node at (-0.1, 1) [left] {1};
      \draw[dashed] (-0.1, 1) -- (0, 1);

      \draw[red] (-2.5, 0.8) -- (2.5, 0.85); 
      \draw[red] (-2.5, 0.4) -- (2.5, 0.38);  
      \node[red] at (2, 1) {0.8};
      \node[red] at (2, 0.2) {0.4};
    \end{tikzpicture}
    \caption{$interpolate(s_1,W)$. Threshold values in red.}
  \end{subfigure}
  \null\hfill
  \medskip

  \begin{subfigure}[b]{0.9\textwidth}
    \centering
    \begin{tikzpicture}[line cap=round,line join=round,>=triangle 45,scale=1.5]
      \clip(-2.5,-0.5) rectangle (2.5,1.25);
      
      \draw (1.6,0) -- (2,0);      
      \draw (1.2,0) -- (1.6,0);
      \draw (0.8,1) -- (1.2,1);
      \draw (0.4,1) -- (0.8,1);
      \draw (0,0) -- (0.4,0);
      \draw (-0.4,0) -- (0,0);
      \draw (-0.4,1) -- (-0.8,1);
      \draw (-0.8,1) -- (-1.2,1);
      \draw (-1.2,0) -- (-1.6,0);
      \draw (-1.6,0) -- (-2,0);
  
      \draw[-] (-2.5,0) -- (2.5,0); 
      \draw[-] (0,-0.5) -- (0,1.45); 
  
      \node at (-0.1, 1) [left] {1};
      \draw (-0.1, 1) -- (0, 1);
  
      \begin{scriptsize}
      \fill [color=black] (1.6,0) circle (1.5pt);
      \fill [color=black] (2,0) circle (1.5pt);      
      \fill [color=black] (1.2,0) circle (1.5pt);
      \fill [color=black] (0.8,1) circle (1.5pt);
      \fill [color=black] (1.2,1) circle (1.5pt);
      \fill [color=black] (0.4,1) circle (1.5pt);
      \fill [color=black] (0.8,1) circle (1.5pt);
      \fill [color=black] (0,0) circle (1.5pt);
      \fill [color=black] (0.4,0) circle (1.5pt);
      \fill [color=black] (-0.4,0) circle (1.5pt);
      \fill [color=black] (0,0) circle (1.5pt);
      \fill [color=black] (-0.4,1) circle (1.5pt);
      \fill [color=black] (-0.8,1) circle (1.5pt);
      \fill [color=black] (-1.2,1) circle (1.5pt);
      \fill [color=black] (-1.2,0) circle (1.5pt);
      \fill [color=black] (-1.6,0) circle (1.5pt);
      \fill [color=black] (-2,0) circle (1.5pt);
      \end{scriptsize}
    \end{tikzpicture}
    \caption{Refinement indicator function $I$.}
  \end{subfigure}
  \vspace*{.25cm}
  \caption{Illustration of Variable Coefficient Elliptic Smoothing algorithm.}
\end{figure}
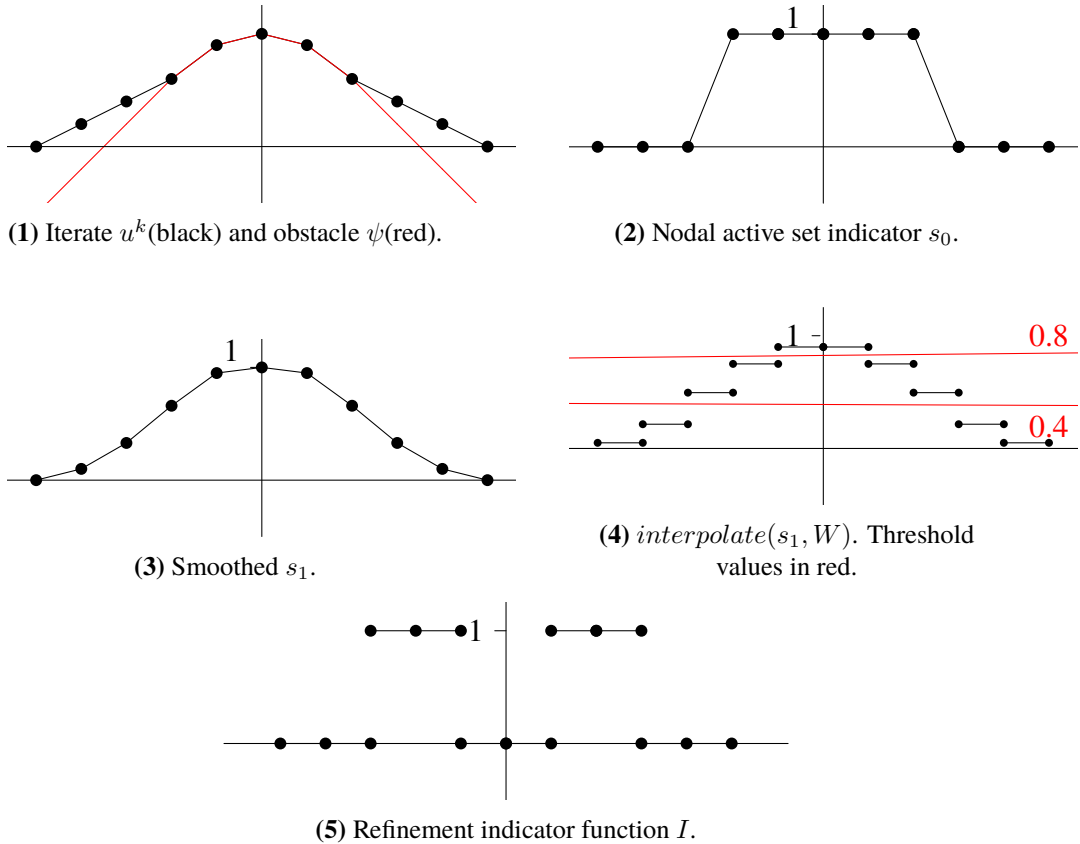

\section{Unstructured Dilation Operator}
Support for unstructured meshes in Firedrake comes from the DMPlex class in PETSc, as developed by \citet{lange_flexible_2015}. DMPlex is a data management object which can store the topology (connectivity of mesh entities) and geometry (coordinates) of a discretization. In the DMPlex object every mesh entity is assigned a unique index. The connectivity of a mesh is stored as a layered directed acyclic graph (DAG) in which a ``covering" relation specifies the edges of the graph. For example, for a tetrahedral element in a 3d mesh, a face is covered by 3 edges and a tetrahedral cell is covered by 4 faces. Each layer represents a class of mesh entity, i.e., vertices, edges, and so on, depending on the dimesion of $\Omega$. Below is an example of how a single tetrahedral cell is represented by DMPlex.
\begin{figure}[H]
  \null\hfill
  \begin{subfigure}[b]{.45\textwidth}
    \centering
    \definecolor{ffqqff}{rgb}{1.,0.,1.}
\definecolor{qqwuqq}{rgb}{0.,0.39215686274509803,0.}
\definecolor{ffqqqq}{rgb}{1.,0.,0.}
\definecolor{ududff}{rgb}{0.30196078431372547,0.30196078431372547,1.}
\begin{tikzpicture}[scale = 1.25, line cap=round,line join=round,>=triangle 45,x=1.0cm,y=1.0cm]
\clip(9.5,-2.5) rectangle (15.5,1.5);
\fill[line width=1.pt,color=qqwuqq,fill=qqwuqq,fill opacity=0.10000000149011612] (5.102909055151545,0.8598788713130618) -- (5.52,-2.96) -- (6.962909055151545,-1.060121128686938) -- cycle;
\fill[line width=2.pt,color=qqwuqq,fill=qqwuqq,fill opacity=0.10000000149011612] (2.6029090551515455,-1.8801211286869384) -- (5.52,-2.96) -- (5.102909055151545,0.8598788713130618) -- cycle;
\fill[line width=2.pt,color=qqwuqq,fill=qqwuqq,fill opacity=0.10000000149011612] (2.6029090551515455,-1.8801211286869384) -- (5.102909055151545,0.8598788713130618) -- (6.962909055151545,-1.060121128686938) -- cycle;
\fill[line width=2.pt,color=qqwuqq,fill=qqwuqq,fill opacity=0.10000000149011612] (2.6029090551515455,-1.8801211286869384) -- (6.962909055151545,-1.060121128686938) -- (5.52,-2.96) -- cycle;
\draw [line width=2.pt,color=ffqqqq] (2.6029090551515455,-1.8801211286869384)-- (5.52,-2.96);
\draw [line width=2.pt,color=ffqqqq] (5.52,-2.96)-- (6.962909055151545,-1.060121128686938);
\draw [line width=2.pt,color=ffqqqq] (2.6029090551515455,-1.8801211286869384)-- (6.962909055151545,-1.060121128686938);
\draw [line width=2.pt,color=ffqqqq] (5.52,-2.96)-- (5.102909055151545,0.8598788713130618);
\draw [line width=2.pt,color=ffqqqq] (2.6029090551515455,-1.8801211286869384)-- (5.102909055151545,0.8598788713130618);
\draw [line width=2.pt,color=ffqqqq] (6.962909055151545,-1.060121128686938)-- (5.102909055151545,0.8598788713130618);
\draw [->,>=stealth,line width=1pt] (11.,1.) -- (10.,0.);
\draw [->,>=stealth,line width=1pt] (11.,1.) -- (14.,0.);
\draw [->,>=stealth,line width=1pt] (11.,1.) -- (12.,0.);
\draw [->,>=stealth,line width=1pt] (12.,1.) -- (14.,0.);
\draw [->,>=stealth,line width=1pt] (12.,1.) -- (13.,0.);
\draw [->,>=stealth,line width=1pt] (12.,1.) -- (15.,0.);
\draw [->,>=stealth,line width=1pt] (13.,1.) -- (10.,0.);
\draw [->,>=stealth,line width=1pt] (13.,1.) -- (11.,0.);
\draw [->,>=stealth,line width=1pt] (13.,1.) -- (13.,0.);
\draw [->,>=stealth,line width=1pt] (14.,1.) -- (11.,0.);
\draw [->,>=stealth,line width=1pt] (14.,1.) -- (12.,0.);
\draw [->,>=stealth,line width=1pt] (14.,1.) -- (15.,0.);
\draw [->,>=stealth,line width=1pt] (10.,0.) -- (12.,-1.);
\draw [->,>=stealth,line width=1pt] (10.,0.) -- (14.,-1.);
\draw [->,>=stealth,line width=1pt] (11.,0.) -- (11.,-1.);
\draw [->,>=stealth,line width=1pt] (13.,0.) -- (11.,-1.);
\draw [->,>=stealth,line width=1pt] (15.,0.) -- (11.,-1.);
\draw [->,>=stealth,line width=1pt] (13.,0.) -- (12.,-1.);
\draw [->,>=stealth,line width=1pt] (14.,0.) -- (12.,-1.);
\draw [->,>=stealth,line width=1pt] (14.,0.) -- (13.,-1.);
\draw [->,>=stealth,line width=1pt] (15.,0.) -- (13.,-1.);
\draw [->,>=stealth,line width=1pt] (12.,0.) -- (13.,-1.);
\draw [->,>=stealth,line width=1pt] (12.,0.) -- (14.,-1.);
\draw [->,>=stealth,line width=1pt] (11.,0.) -- (14.,-1.);
\draw [->,>=stealth,line width=1pt] (11.,-1.) -- (12.5,-2.);
\draw [->,>=stealth,line width=1pt] (12.,-1.) -- (12.5,-2.);
\draw [->,>=stealth,line width=1pt] (13.,-1.) -- (12.5,-2.);
\draw [->,>=stealth,line width=1pt] (14.,-1.) -- (12.5,-2.);
\begin{scriptsize}
\draw [fill=ududff] (6.962909055151545,-1.060121128686938) circle (2.5pt);
\draw[color=ududff] (7.024972711500066,-0.8843347036153757) node {4};
\draw [fill=ududff] (5.52,-2.96) circle (2.5pt);
\draw[color=ududff] (5.64776858075123,-3.0023796771118603) node {3};
\draw [fill=ududff] (5.102909055151545,0.8598788713130618) circle (2.5pt);
\draw[color=ududff] (5.457809390303115,0.9772653627761532) node {2};
\draw [fill=ududff] (2.6029090551515455,-1.8801211286869384) circle (2.5pt);
\draw[color=ududff] (5.457809390303115,-1.6726653439750538) node {1};
\draw[color=ffqqqq] (4.042613421464657,-2.4799919033795437) node {5};
\draw[color=ffqqqq] (6.38860942349888,-2.0145918867816612) node {6};
\draw[color=ffqqqq] (4.840442021346742,-1.5301959511389676) node {7};
\draw[color=ffqqqq] (5.486303268870333,-0.950820420272216) node {8};
\draw[color=ffqqqq] (3.9951236238526286,-0.5234122417639568) node {9};
\draw[color=ffqqqq] (6.00869104260265,0.1984326819388808) node {10};
\draw[color=qqwuqq] (5.875719609288969,-1.473208194004533) node {11};
\draw[color=qqwuqq] (4.470021599972917,-1.2452571654667948) node {12};
\draw[color=qqwuqq] (4.95441753561561,-0.6088938774656087) node {13};
\draw[color=qqwuqq] (5.087388968929291,-1.8816204534679806) node {14};
\draw [fill=ffqqff] (5.014364521362902,-1.6680568902556219) circle (0.5pt);
\draw[color=ffqqff] (5.524295106959956,-1.5586898297061849) node {0};
\draw [fill=ffqqff] (12.5,-2.) circle (1.75pt);
\draw[color=ffqqff] (12.457805558316158,-2.1665592391401534) node {0};
\draw [fill=qqwuqq] (11.,-1.) circle (1.75pt);
\draw[color=qqwuqq] (10.748172844283122,-1.0268040964514622) node {11};
\draw [fill=qqwuqq] (12.,-1.) circle (1.75pt);
\draw[color=qqwuqq] (11.764454513180539,-1.0837918535858968) node {12};
\draw [fill=qqwuqq] (13.,-1.) circle (1.75pt);
\draw[color=qqwuqq] (12.685756586853898,-0.9698163393170275) node {13};
\draw [fill=qqwuqq] (14.,-1.) circle (1.75pt);
\draw[color=qqwuqq] (14.129446434259574,-0.8558408250481585) node {14};
\draw [fill=ffqqqq] (10.,0.) circle (1.75pt);
\draw[color=ffqqqq] (9.731891175385707,0.07495920814760594) node {5};
\draw [fill=ffqqqq] (11.,0.) circle (1.75pt);
\draw[color=ffqqqq] (10.824156520462367,0.09395512719241746) node {6};
\draw [fill=ffqqqq] (12.,0.) circle (1.75pt);
\draw[color=ffqqqq] (11.726462675090914,0.09395512719241746) node {7};
\draw [fill=ffqqqq] (13.,0.) circle (1.75pt);
\draw[color=ffqqqq] (13.17015252249659,0.10345308671482321) node {8};
\draw [fill=ffqqqq] (14.,0.) circle (1.75pt);
\draw[color=ffqqqq] (14.11045051521476,0.12244900575963474) node {9};
\draw [fill=ffqqqq] (15.,0.) circle (1.75pt);
\draw[color=ffqqqq] (15.145728103156989,0.11295104623722897) node {10};
\draw [fill=ududff] (11.,1.) circle (1.75pt);
\draw[color=ududff] (11.071103468044917,1.214714350836297) node {1};
\draw [fill=ududff] (12.,1.) circle (1.75pt);
\draw[color=ududff] (12.068389217897522,1.214714350836297) node {2};
\draw [fill=ududff] (13.,1.) circle (1.75pt);
\draw[color=ududff] (13.065674967750127,1.214714350836297) node {3};
\draw [fill=ududff] (14.,1.) circle (1.75pt);
\draw[color=ududff] (14.062960717602731,1.214714350836297) node {4};
\end{scriptsize}
\end{tikzpicture}
\caption{DAG representation of tetrahedral cell.}
\end{subfigure}
\hfill
\begin{subfigure}[b]{0.45\textwidth}
  \centering
  \definecolor{ffqqff}{rgb}{1.,0.,1.}
  \definecolor{qqwuqq}{rgb}{0.,0.39215686274509803,0.}
  \definecolor{ffqqqq}{rgb}{1.,0.,0.}
  \definecolor{ududff}{rgb}{0.30196078431372547,0.30196078431372547,1.}
  \begin{tikzpicture}[line cap=round,line join=round,>=triangle 45,x=1.0cm,y=1.0cm]
  \clip(2.,-3.25) rectangle (7.5,1.25);
  
  \fill[line width=2.pt,color=qqwuqq,fill=qqwuqq,fill opacity=0.10000000149011612] (5.102909055151545,0.8598788713130618) -- (5.52,-2.96) -- (6.962909055151545,-1.060121128686938) -- cycle;
  \fill[line width=2.pt,color=qqwuqq,fill=qqwuqq,fill opacity=0.10000000149011612] (2.6029090551515455,-1.8801211286869384) -- (5.52,-2.96) -- (5.102909055151545,0.8598788713130618) -- cycle;
  \fill[line width=2.pt,color=qqwuqq,fill=qqwuqq,fill opacity=0.10000000149011612] (2.6029090551515455,-1.8801211286869384) -- (5.102909055151545,0.8598788713130618) -- (6.962909055151545,-1.060121128686938) -- cycle;
  \fill[line width=2.pt,color=qqwuqq,fill=qqwuqq,fill opacity=0.10000000149011612] (2.6029090551515455,-1.8801211286869384) -- (6.962909055151545,-1.060121128686938) -- (5.52,-2.96) -- cycle;
  
  \draw [line width=1.pt,color=ffqqqq] (2.6029090551515455,-1.8801211286869384)-- (5.52,-2.96);
  \draw [line width=1.pt,color=ffqqqq] (5.52,-2.96)-- (6.962909055151545,-1.060121128686938);
  \draw [line width=1.pt, color=ffqqqq, dash pattern=on 3pt off 3pt] (2.6029090551515455,-1.8801211286869384) -- (6.962909055151545,-1.060121128686938);
  \draw [line width=1.pt,color=ffqqqq] (5.52,-2.96)-- (5.102909055151545,0.8598788713130618);
  \draw [line width=1.pt,color=ffqqqq] (2.6029090551515455,-1.8801211286869384)-- (5.102909055151545,0.8598788713130618);
  \draw [line width=1.pt,color=ffqqqq] (6.962909055151545,-1.060121128686938)-- (5.102909055151545,0.8598788713130618);
  
  \draw [->,line width=.5pt] (11.,1.) -- (10.,0.);
  \draw [->,line width=.5pt] (11.,1.) -- (14.,0.);
  \draw [->,line width=.5pt] (11.,1.) -- (12.,0.);
  \draw [->,line width=.5pt] (12.,1.) -- (14.,0.);
  \draw [->,line width=.5pt] (12.,1.) -- (13.,0.);
  \draw [->,line width=.5pt] (12.,1.) -- (15.,0.);
  \draw [->,line width=.5pt] (13.,1.) -- (10.,0.);
  \draw [->,line width=.5pt] (13.,1.) -- (11.,0.);
  \draw [->,line width=.5pt] (13.,1.) -- (13.,0.);
  \draw [->,line width=.5pt] (14.,1.) -- (11.,0.);
  \draw [->,line width=.5pt] (14.,1.) -- (12.,0.);
  \draw [->,line width=.5pt] (14.,1.) -- (15.,0.);
  \draw [->,line width=.5pt] (10.,0.) -- (12.,-1.);
  \draw [->,line width=.5pt] (10.,0.) -- (14.,-1.);
  \draw [->,line width=.5pt] (11.,0.) -- (11.,-1.);
  \draw [->,line width=.5pt] (13.,0.) -- (11.,-1.);
  \draw [->,line width=.5pt] (15.,0.) -- (11.,-1.);
  \draw [->,line width=.5pt] (13.,0.) -- (12.,-1.);
  \draw [->,line width=.5pt] (14.,0.) -- (12.,-1.);
  \draw [->,line width=.5pt] (14.,0.) -- (13.,-1.);
  \draw [->,line width=.5pt] (15.,0.) -- (13.,-1.);
  \draw [->,line width=.5pt] (12.,0.) -- (13.,-1.);
  \draw [->,line width=.5pt] (12.,0.) -- (14.,-1.);
  \draw [->,line width=.5pt] (11.,0.) -- (14.,-1.);
  \draw [->,line width=.5pt] (11.,-1.) -- (12.5,-2.);
  \draw [->,line width=.5pt] (12.,-1.) -- (12.5,-2.);
  \draw [->,line width=.5pt] (13.,-1.) -- (12.5,-2.);
  \draw [->,line width=.5pt] (14.,-1.) -- (12.5,-2.);
  \begin{scriptsize}
  \draw [fill=ududff] (6.962909055151545,-1.060121128686938) circle (2.5pt);
  \draw[color=ududff] (7.03,-0.87) node {4};
  \draw [fill=ududff] (5.52,-2.96) circle (2.5pt);
  \draw[color=ududff] (5.7,-3.0023796771118603) node {3};
  \draw [fill=ududff] (5.102909055151545,0.8598788713130618) circle (2.5pt);
  \draw[color=ududff] (5.15,1.1) node {2};
  \draw [fill=ududff] (2.6029090551515455,-1.8801211286869384) circle (2.5pt);
  \draw[color=ududff] (2.4,-1.7011592225422711) node {1};
  \draw[color=ffqqqq] (3.55,-2.413506186722703) node {5};
  \draw[color=ffqqqq] (6.388609423498879,-2.0145918867816612) node {6};
  \draw[color=ffqqqq] (4.21357669286796,-1.397224517825287) node {7};
  \draw[color=ffqqqq] (5.3,.1) node {8};
  \draw[color=ffqqqq] (4.128095057166308,0.14144492480444626) node {9};
  \draw[color=ffqqqq] (6.008691042602648,0.1984326819388808) node {10};
  \draw[color=qqwuqq] (5.875719609288968,-1.473208194004533) node {11};
  \draw[color=qqwuqq] (4.669478749943437,-1.7011592225422711) node {12};
  \draw[color=qqwuqq] (4.707470588033059,-0.637387756032826) node {13};
  \draw[color=qqwuqq] (5.115882847496508,-2.2330449557969936) node {14};
  \draw[color=ffqqff] (5.077891009406884,-1.5681877892285905) node {0};
  \draw [fill=ffqqff] (12.5,-2.) circle (1.75pt);
  \draw[color=ffqqff] (12.457805558316156,-2.1665592391401534) node {0};
  \draw [fill=qqwuqq] (11.,-1.) circle (1.75pt);
  \draw[color=qqwuqq] (10.74817284428312,-1.0268040964514622) node {11};
  \draw [fill=qqwuqq] (12.,-1.) circle (1.75pt);
  \draw[color=qqwuqq] (11.764454513180537,-1.0837918535858968) node {12};
  \draw [fill=qqwuqq] (13.,-1.) circle (1.75pt);
  \draw[color=qqwuqq] (12.685756586853895,-0.9698163393170275) node {13};
  \draw [fill=qqwuqq] (14.,-1.) circle (1.75pt);
  \draw[color=qqwuqq] (14.12944643425957,-0.8558408250481585) node {14};
  \draw [fill=ffqqqq] (10.,0.) circle (1.75pt);
  \draw[color=ffqqqq] (9.731891175385705,0.07495920814760594) node {5};
  \draw [fill=ffqqqq] (11.,0.) circle (1.75pt);
  \draw[color=ffqqqq] (10.824156520462365,0.09395512719241746) node {6};
  \draw [fill=ffqqqq] (12.,0.) circle (1.75pt);
  \draw[color=ffqqqq] (11.726462675090913,0.09395512719241746) node {7};
  \draw [fill=ffqqqq] (13.,0.) circle (1.75pt);
  \draw[color=ffqqqq] (13.170152522496588,0.10345308671482321) node {8};
  \draw [fill=ffqqqq] (14.,0.) circle (1.75pt);
  \draw[color=ffqqqq] (14.110450515214758,0.12244900575963474) node {9};
  \draw [fill=ffqqqq] (15.,0.) circle (1.75pt);
  \draw[color=ffqqqq] (14.480870936588584,0.0844571676700117) node {10};
  \draw [fill=ududff] (11.,1.) circle (1.75pt);
  \draw[color=ududff] (11.071103468044916,1.214714350836297) node {1};
  \draw [fill=ududff] (12.,1.) circle (1.75pt);
  \draw[color=ududff] (12.06838921789752,1.214714350836297) node {2};
  \draw [fill=ududff] (13.,1.) circle (1.75pt);
  \draw[color=ududff] (13.065674967750125,1.214714350836297) node {3};
  \draw [fill=ududff] (14.,1.) circle (1.75pt);
  \draw[color=ududff] (14.06296071760273,1.214714350836297) node {4};
  \end{scriptsize}
  \end{tikzpicture}
  \caption{Mesh entity numbering of tetrahedral cell. Vertices(blue), Edges(Red), Faces(green), Cells(pink)}
\end{subfigure}
\null\hfill
\caption{DMPlex representation of a tetrahedral cell \citep{lange_flexible_2015}.}
 \end{figure}
The DMPlex object has several methods which make querying the mesh topology and geometry simple \citep{lange_efficient_2016}. For example, let $p$ be an index assigned by DMPlex. Then $\emph{cone}(p)$ returns all the in-neighbors of $p$. In the example above $\emph{cone}(0) = \{11, 12, 13, 14\}$. The transitive closure of $\emph{cone}(p)$ is also available with $\emph{closure}(p)$. The dual of $\emph{cone}(p)$ is $\emph{support}(p)$ which returns all the out-neighbors of $p$. In the example above $\emph{support}(6) = \{6, 7, 10\}$. The transitive closure of $\emph{support}(p)$ is also available with $\emph{star}(p)$. 

The use of DMPlex queries is essential to our second strategy, the Unstructured Dilation Operator (UDO). We identify the set $B$ of elements that border the computed free boundary $\Gamma_{u^k}$ by interpolating a nodal active indicator function into DG0 and thresholding for values in the range (0, 1). We then use the \emph{closure} and \emph{star} methods to create vertex-to-cell and cell-to-vertex mappings. These mappings are then used to determine which elements are neighbors to the computed free boundary. We say an element neighbors another if it shares at least one vertex. The function $N(\triangle)$ returns a set of elements:
\begin{equation}
  N(\triangle) = \{\triangle_i \in T: \triangle \text{ shares at least 1 vertex with } \triangle_i\}.
\end{equation}
This process is then repeated $n$ times to create a set of elements that are within $n$ neighborhood levels of the border elements:
\begin{equation}
  N^n(\triangle) = \underbrace{N(...N(\triangle))}_{n \text{ times}}.
\end{equation}
As shown in Algorithm 3, for a border element set $B$, defined below in (5.5) we use breadth-first search to construct the set $N^n(B)$ and then assemble its corresponding indicator function. This process expand the support of the DG0 indicator function in a way that resembles the ``dilation" operation in image processing as seen in \citep{OpenCV} but it is applied over an unstructured mesh.

\begin{algorithm}[H]
  \caption{Unstructured Dilation Operator Element Tagging for VIs}
  \begin{algorithmic}[1]
    \Require $tol \in \mathbb{R}$, $u^k \in K, \psi \in V$, $W$ is a DGO FE space, mesh $T$.
    \Require Neighborhood Depth Parameter $n \in \mathbb{N}$.
    \State Compute the nodal active set indicator function $s$
    \begin{equation}
    s = \begin{cases}
      1 & \text{ if } u^k - \psi < tol\\
      0 & \text{ otherwise}
    \end{cases}.
    \end{equation}
  
    \State Let $s_W = interpolate(s, W)$ .
    \State Define the border element set $B$:
    \begin{equation}
    B = \{\triangle \in T: 0 < s_W(\triangle) < 1 \} .
    \end{equation}

    \State Use the \emph{closure} and \emph{star} methods to construct vertex-to-cell and cell-to-vertex mappings.

    \State Use breadth-first search to construct the set $N^n(B)$ as in (5.2) and (5.3). 
    \State Assemble the indicator function $I$ of the set $N^n(B)$. \\
    \Return $I$
  \end{algorithmic}
  \end{algorithm}

\section{Results: Dependence on Parameters}

Parameters in each method allow us vary the area of refinement about the free boundary. The following figures illustrate the effect of varying the thresholding parameters of the VCES method, and the neighborhood depth parameter of the UDO method. Figures 5.3 and 5.4 show results from 7 iterations of the refinment loop applied to the reference ball obstacle problem, where the initial mesh is generated using a uniform triangle height parameter set to $0.2$. 
\begin{figure}[H]
  \centering
  \begin{subfigure}[b]{0.45\textwidth}
      \centering
      \includegraphics[width=\textwidth]{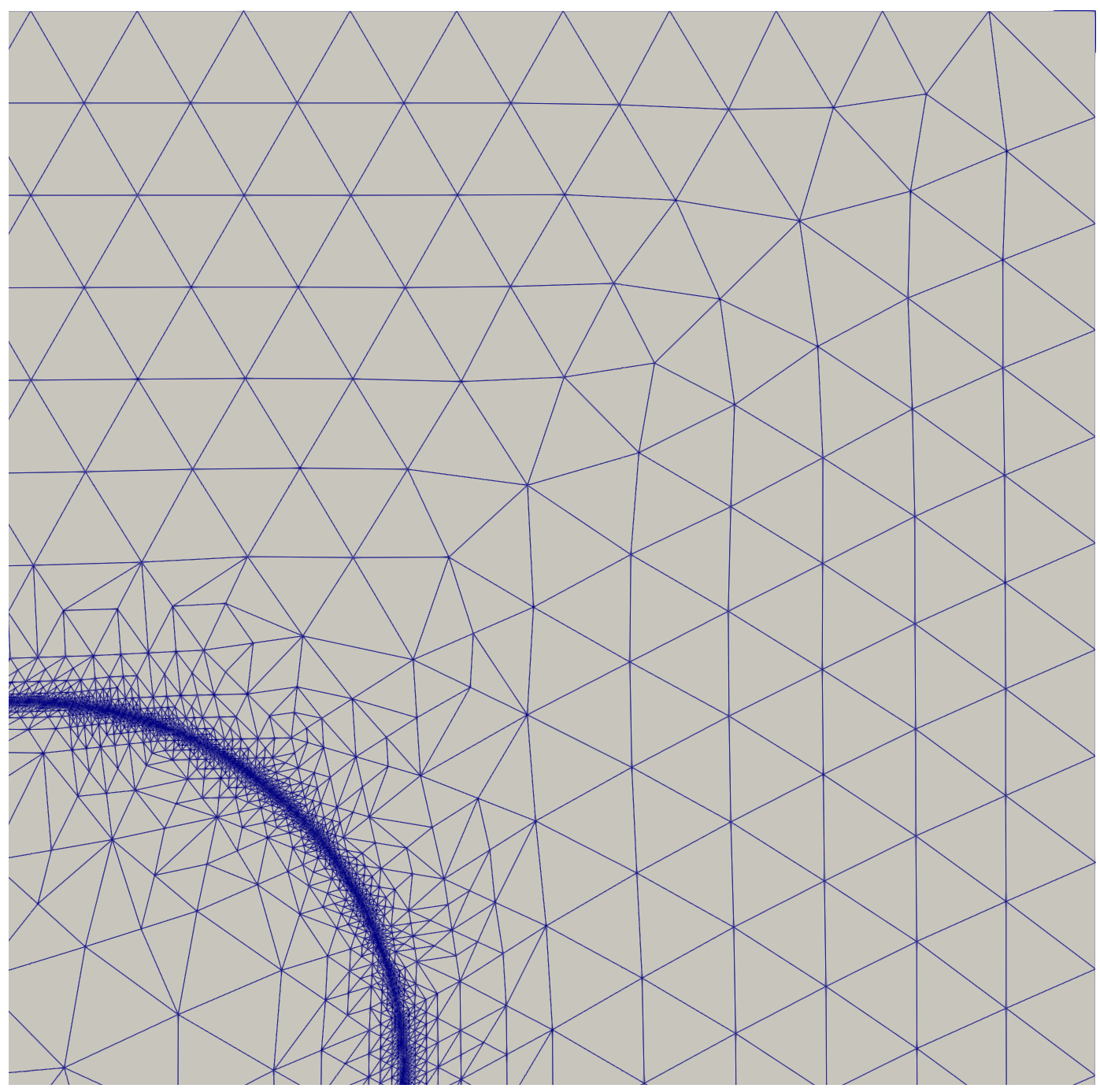}
      \caption{Thresholding parameters (.45, .65).}
      \label{fig:image1}
  \end{subfigure}
  \hfill
  \begin{subfigure}[b]{0.45\textwidth}
      \centering
      \includegraphics[width=\textwidth]{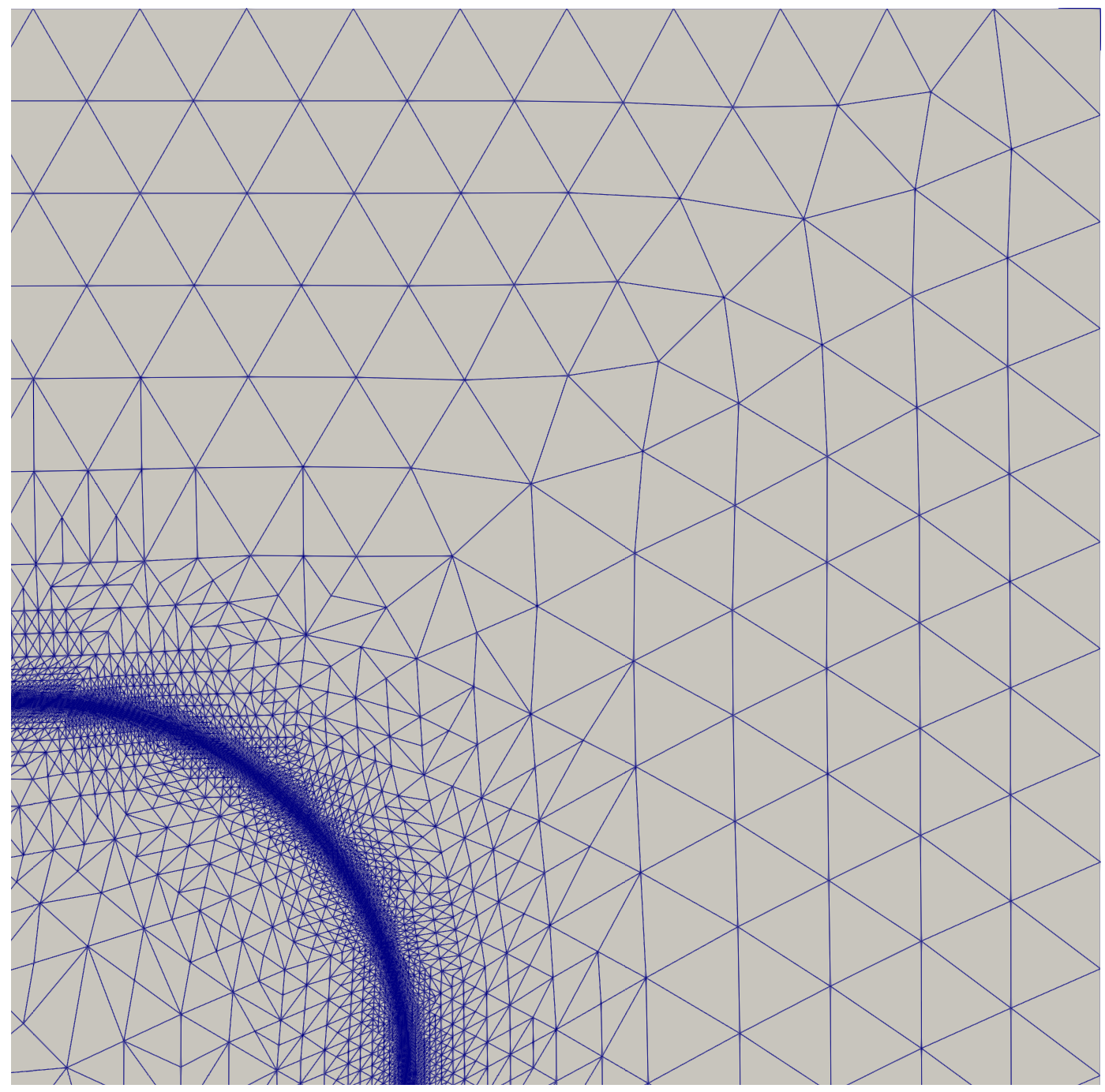}
      \caption{Thresholding parameters (.1, .9)}
      \label{fig:image2}
  \end{subfigure}
  \vspace*{.25cm}
  \caption{7 iterations of VCES method applied to the classical obstacle problem.}
  \label{fig:two_images}
\end{figure}

\begin{figure}[H]
  \centering
  \begin{subfigure}[b]{0.45\textwidth}
      \centering
      \includegraphics[width=\textwidth]{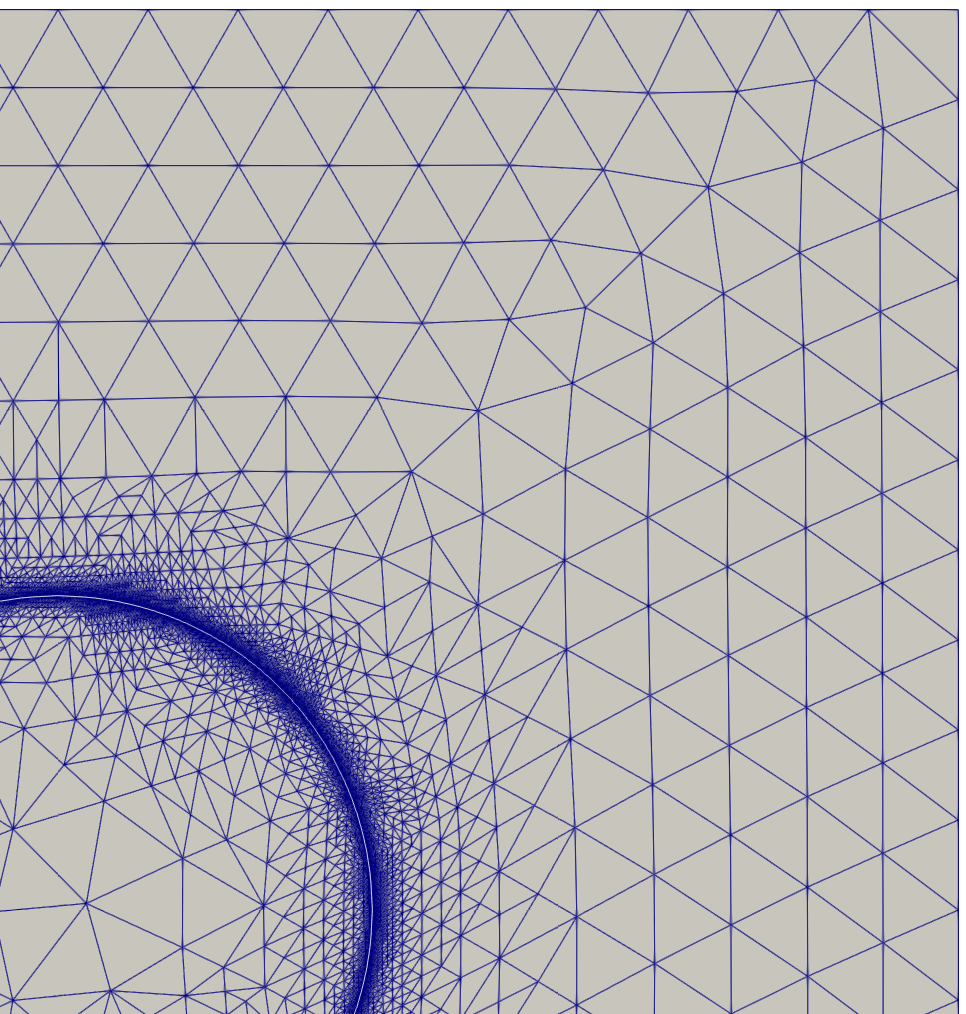}
      \caption{$n = 1$}
      \label{fig:image1}
  \end{subfigure}
  \hfill
  \begin{subfigure}[b]{0.45\textwidth}
      \centering
      \includegraphics[width=\textwidth]{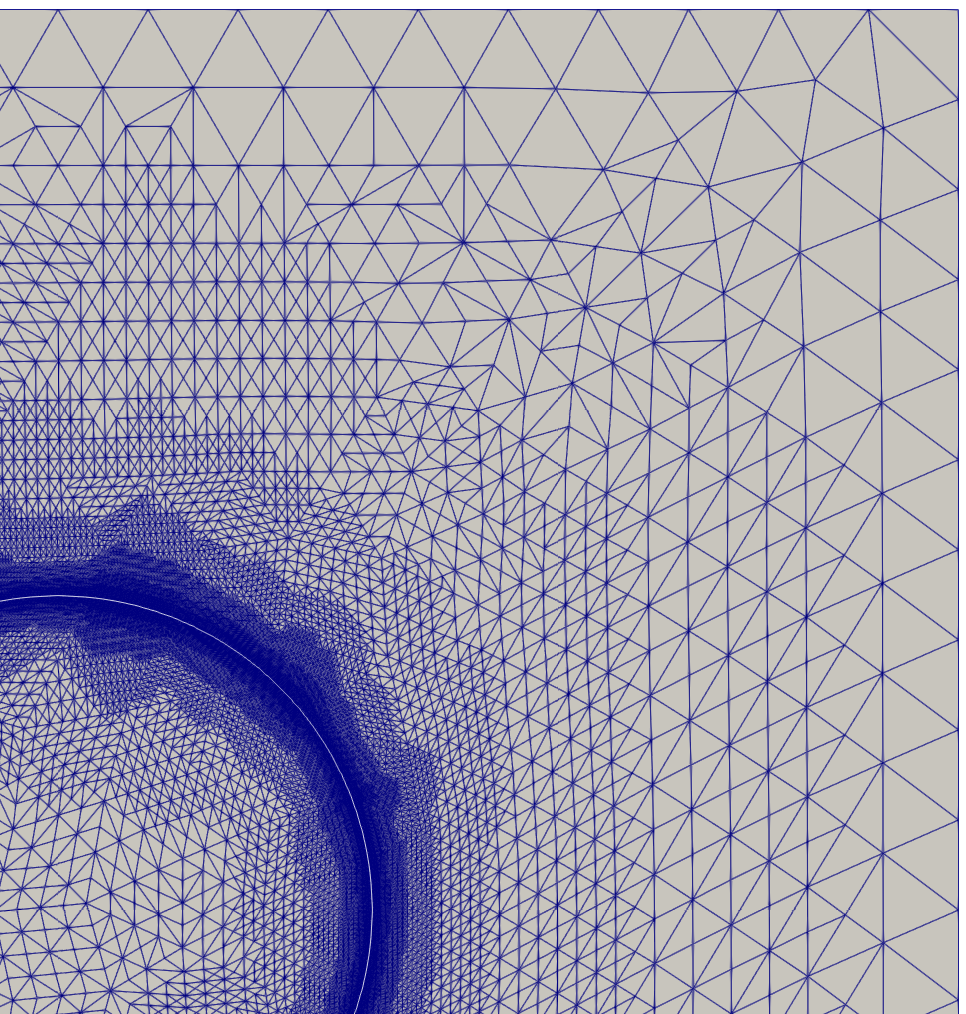}
      \caption{$n = 5$}
      \label{fig:image2}
  \end{subfigure}
  \vspace*{.25cm}
  \caption{7 iterations of UDO method applied to ball obstacle problem.}
  \label{fig:two_images}
\end{figure}
In the following experiment we ran 7 iterations of the refinement loop using both methods on the ball obstacle reference problem, using the same inital mesh as in the previous figures, while varying the thresholding and neighborhood depth parameters. 
\begin{figure}[H]
  \centering
  \begin{subfigure}[b]{.80\textwidth}
      \centering
      \includegraphics[width=\textwidth]{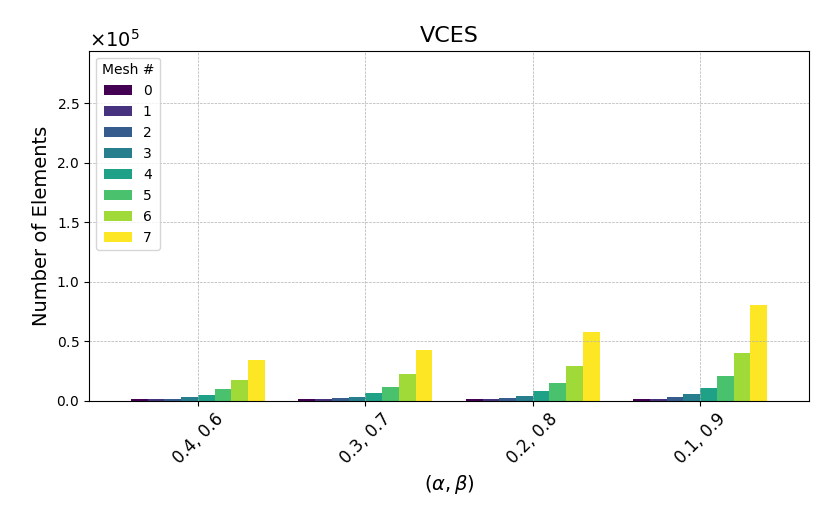}
      \caption{Varying threshold parameters for VCES.}
      \label{fig:top}
  \end{subfigure}
  
  \vspace{1em} 

  \begin{subfigure}[b]{.80\textwidth}
      \centering
      \includegraphics[width=\textwidth]{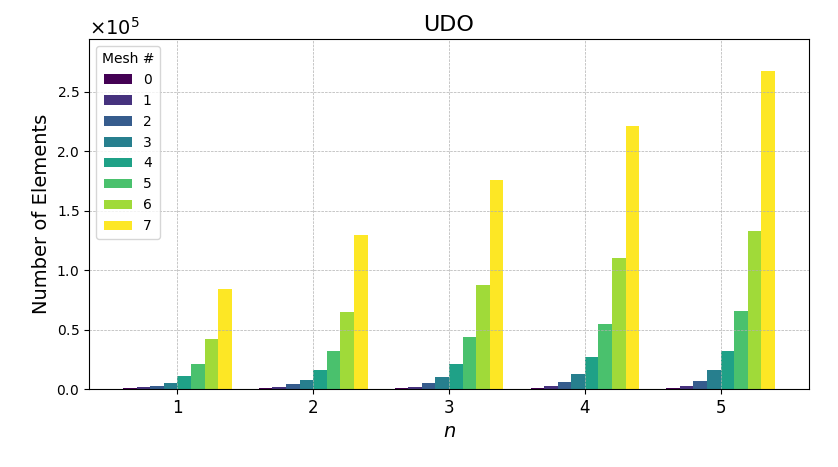}
      \caption{Varying neighborhood depth for UDO.}
      \label{fig:bottom}
  \end{subfigure}
  
  \caption{Comparison of parameters for VCES and UDO methods}
  \label{fig:combined}
\end{figure}

Both methods are also effective in the context of a more exotic obstacles. We consider a ``spiral" problem from \citet{graser_multigrid_2009}, which is a classical obstacle problem with obstacle function given in polar coordinates by 
\begin{equation}
  \psi(r, \theta) = \sin(2\pi/r + \pi/2 - \theta) + \frac{r(r + 1)}{r - 2} - 3r + 3.6,
\end{equation}
and $\psi(0, \theta) = 3.6$. This spiral produces a spiral active set as illustrated in Figure 5.6, computed on a high-resolution unifrom mesh. 
\begin{figure}[H]
  \centering
  \includegraphics[width=.45\textwidth]{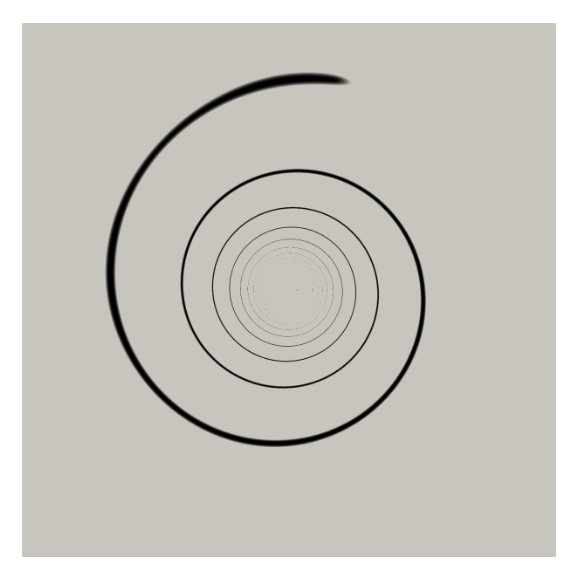}
  \caption{Active set for the spiral obstacle problem \citep{bueler_full_2023}.}
\end{figure} 

\begin{figure}[H]
  \centering
  \begin{subfigure}[b]{0.45\textwidth}
      \centering
      \includegraphics[width=\textwidth]{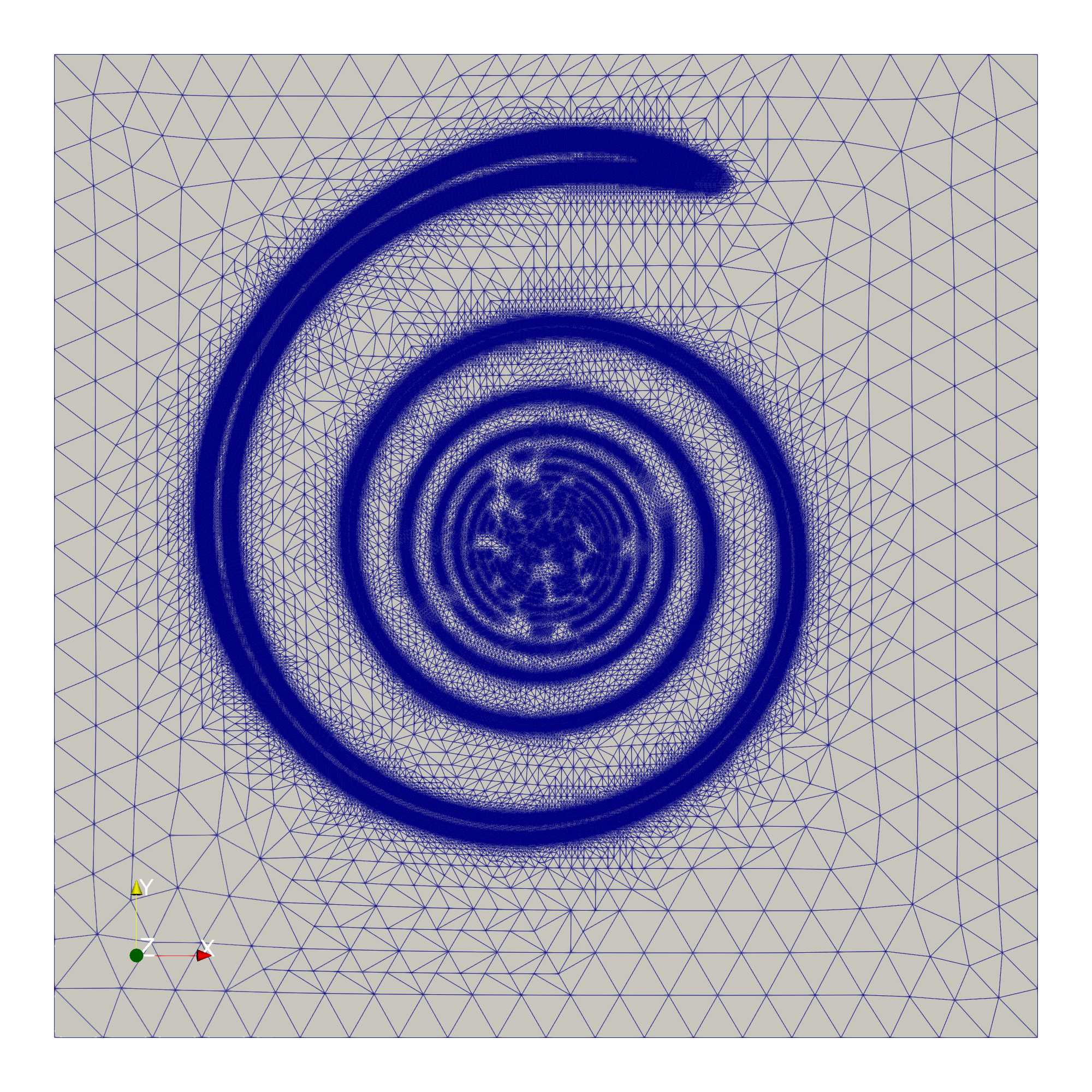}
      \caption{VCES: thresholding parameters (.1, .9).}
      \label{fig:image1}
  \end{subfigure}
  \hfill
  \begin{subfigure}[b]{0.45\textwidth}
      \centering
      \includegraphics[width=\textwidth]{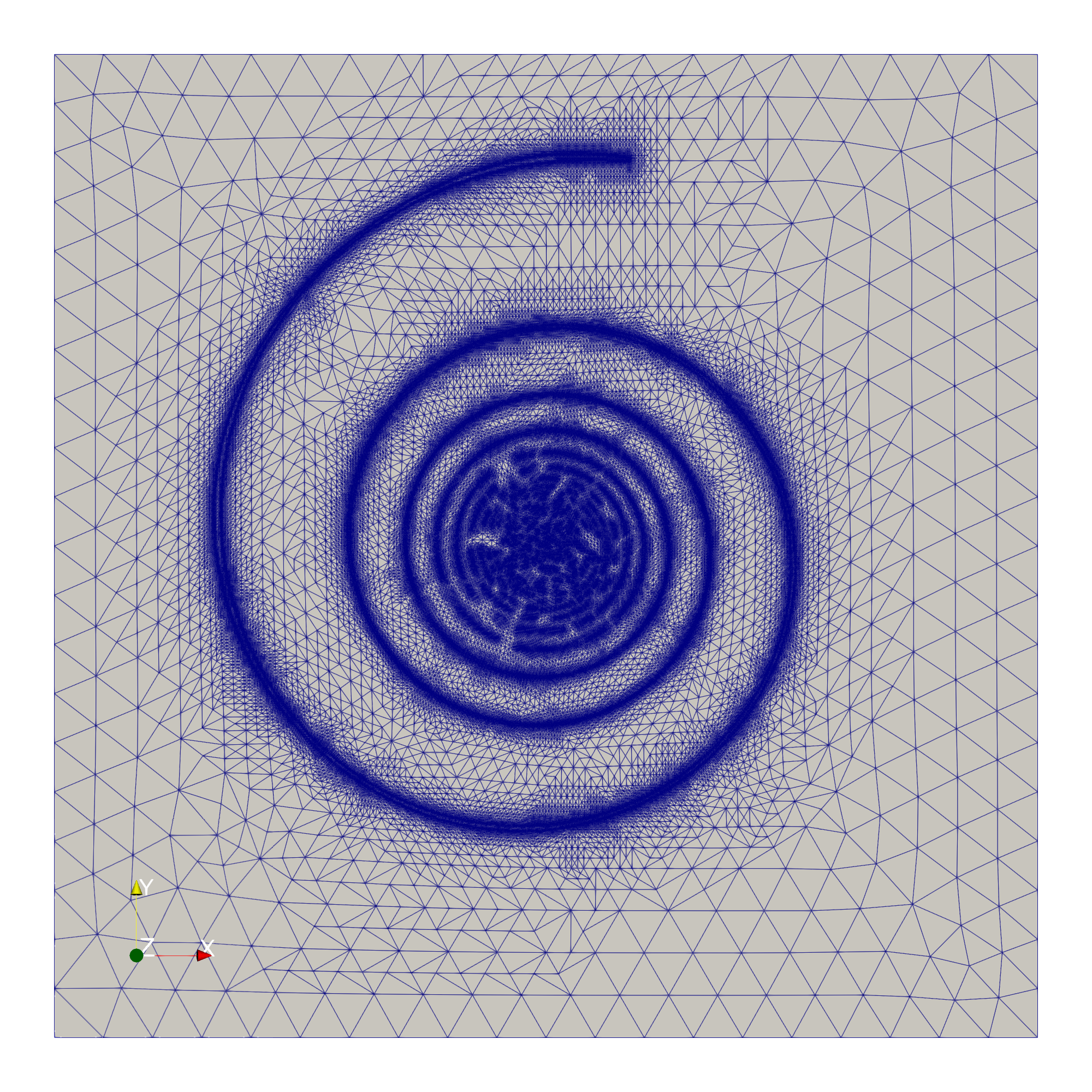}
      \caption{UDO: $n = 1$}
      \label{fig:image2}
  \end{subfigure}
  \vspace*{.25cm}
  \caption{7 iterations of both methods applied to the spiral problem.}
  \label{fig:two_images}
\end{figure}

As we will see in Chapter 6 these strategies can be used in conjunction with uniform refinement, as well as with other tagging methods, to create a more robust adaptive refinement strategy.

\chapter{Results}
\section{Metrics for Refinement Quality Near the Free Boundary}
The purpose of these methods is to reduce the error in the solution that arises from poorly approximating the free boundary. To test the effectiveness of these methods the typical metrics from FEM error analysis like the $L_2$ error are not sufficient. Instead we propose two heuristics which we can apply specifically to the obstacle problem.

The first is the Jaccard index \citep{kosub_note_2016}. If $A_u$ is the true active set and $A_{\hat{u}}$ is the computed active set on a given mesh, then the Jaccard index of $A_u$ and $A_{\hat{u}}$ is defined as
\begin{equation}
  J(A_u, A_{\hat{u}}) = \frac{|A_u \cap A_{\hat{u}}|}{|A_u \cup A_{\hat{u}}|}.
\end{equation}
Note that $J(\cdot, \cdot) \in [0, 1]$ where $J = 1$ indicates that the computed active set is identical to the true active set, and $J = 0$ indicates otherwise.   

The second heuristic we will use is Hausdorff distance. Let $(M, d)$ be a metric space, and consider a pair of non-empty subsets $X \subset M$ and $Y \subset M$. The Hausdorff distance between $X$ and $Y$ is defined as 
\begin{equation}
d_H(X, Y) = \max\left\{\sup_{x \in X} \inf_{y \in Y} d(x, y), \sup_{y \in Y} \inf_{x \in X} d(x, y)\right\}
\end{equation}
To accurately compute the Hausdorff distance, both curves are densely sampled with a substantial number of points. The Hausdorff distance is then calculated between these two sets of sampled points. For a more formal definition and a detailed discussion regarding the Hausdorff distance on finite and infinite sets, see \citet{jungeblut_complexity_2022}.

\section{Experiments}
Recall that the solution to the ball obstacle reference problem is known, so we have access to $u$, $A_u$, and $\Gamma_u$ to compute $L_2$ error, Jaccard index, and Hausdorff distance at each iteration. We will see in the results below that neither of our proposed methods will effectively address the discretization error over the inactive set. To address this issue we propose hybrid refinement strategies where our adaptive methods are used in conjunction with uniform refinement. For demonstration purposes we've devised a hybrid strategy which applies uniform refinement to the entire domain when $d_H(\Gamma_u, \Gamma_{\hat{u}}) < h^2$, where $h$ is the largest cell diameter in the inactive set, and $\Gamma_u$ and $\Gamma_{\hat{u}}$ are the true and computed free boundaries respectively. 

The following convergence results were obtained by applying VCES and UDO to the ball obstacle reference problem. For all cases the coarsest mesh was constructed using the default netgen mesh constructor with the triangle height parameter set to $.45$, resulting in a mesh with 113 vertices (188 elements). The mesh was then refined 7 times using uniform, adaptive, or hybrid refinement strategies. The adaptive strategy used either VCES or UDO at every refinement iteration. When VCES was applied, thresholding parameters where set to $(.2, .8)$ and when UDO was applied the neighborhood depth parameter was set to $3$. The results of the convergence experiments are in Figures 6.1-6.6.

\begin{figure}[H]
  \centering
  \includegraphics[width=\textwidth]{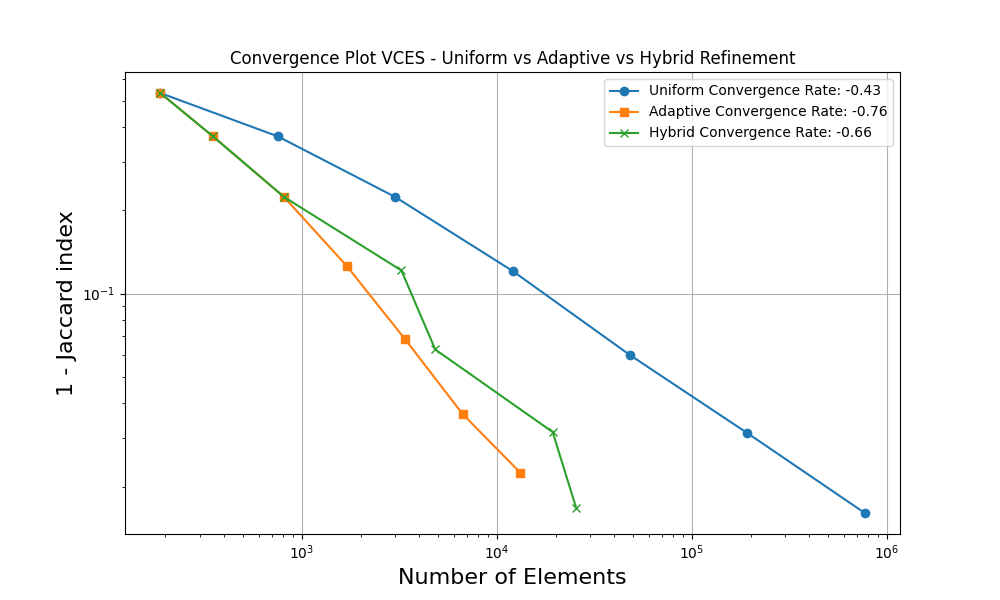}
  \caption{1 - Jaccard index convergence plot for VCES. Lower is better.} 
\end{figure}

\begin{figure}[H]
  \centering
  \includegraphics[width=\textwidth]{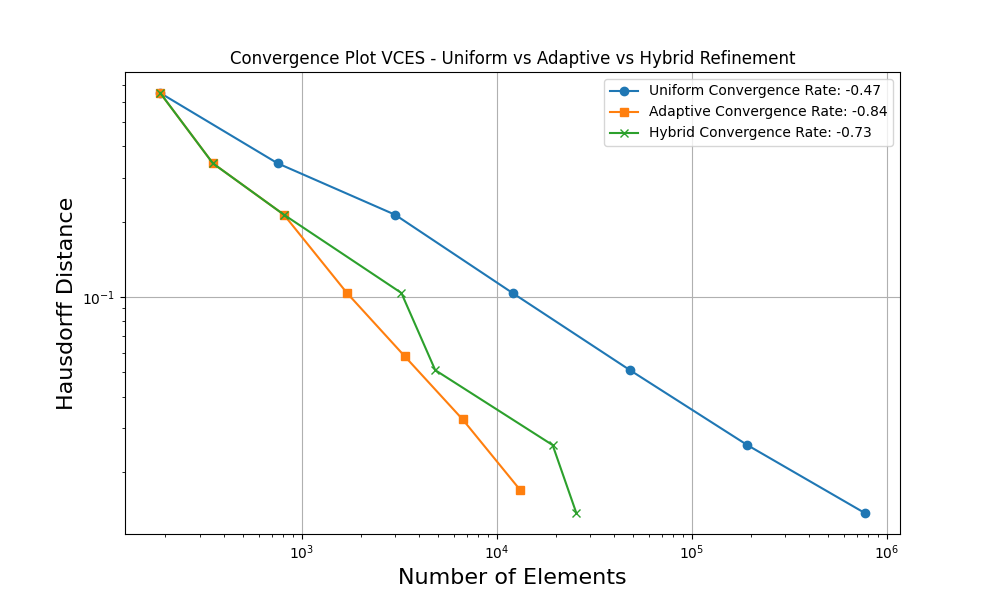}
  \caption{Hausdorff distance convergence plot for VCES. Lower is better.} 
\end{figure}

\begin{figure}[H]
  \centering
  \includegraphics[width=\textwidth]{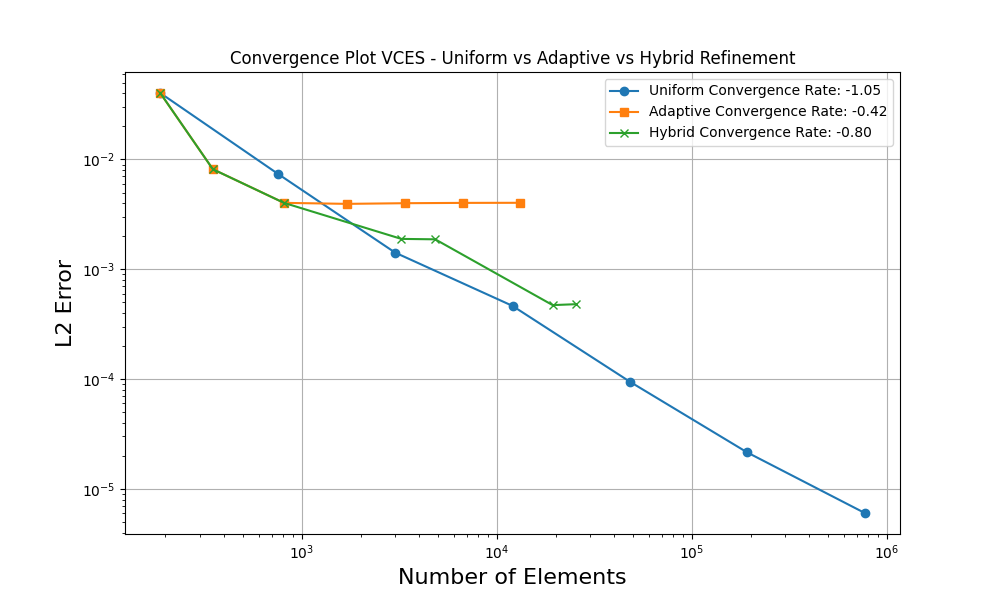}
  \caption{$L_2$ convergence plot for VCES. Lower is better.} 
\end{figure}

\begin{figure}[H]
  \centering
  \includegraphics[width=\textwidth]{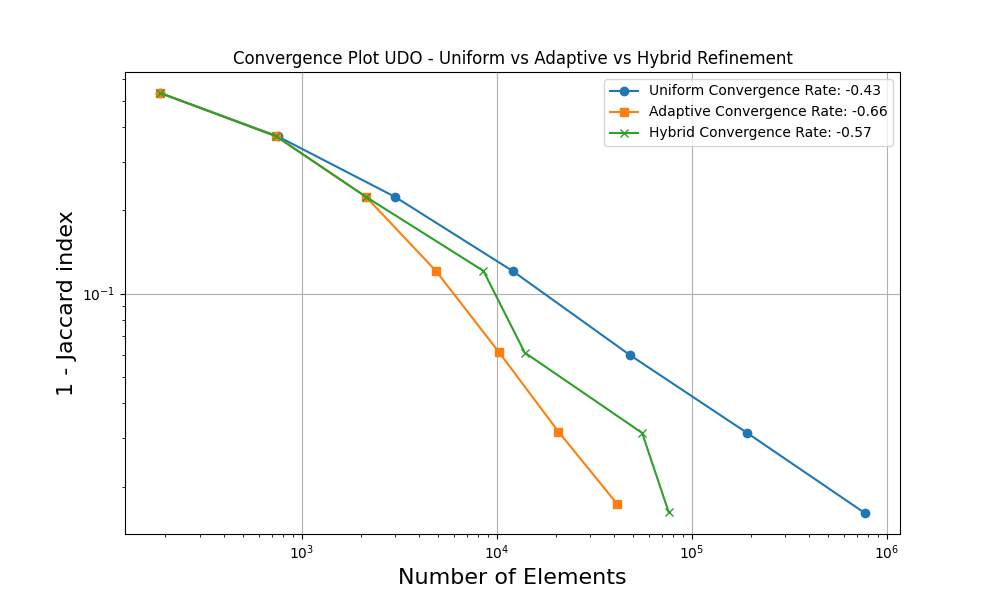}
  \caption{1 - Jaccard index convergence plot for UDO. Lower is better.} 
\end{figure}

\begin{figure}[H]
  \centering
  \includegraphics[width=\textwidth]{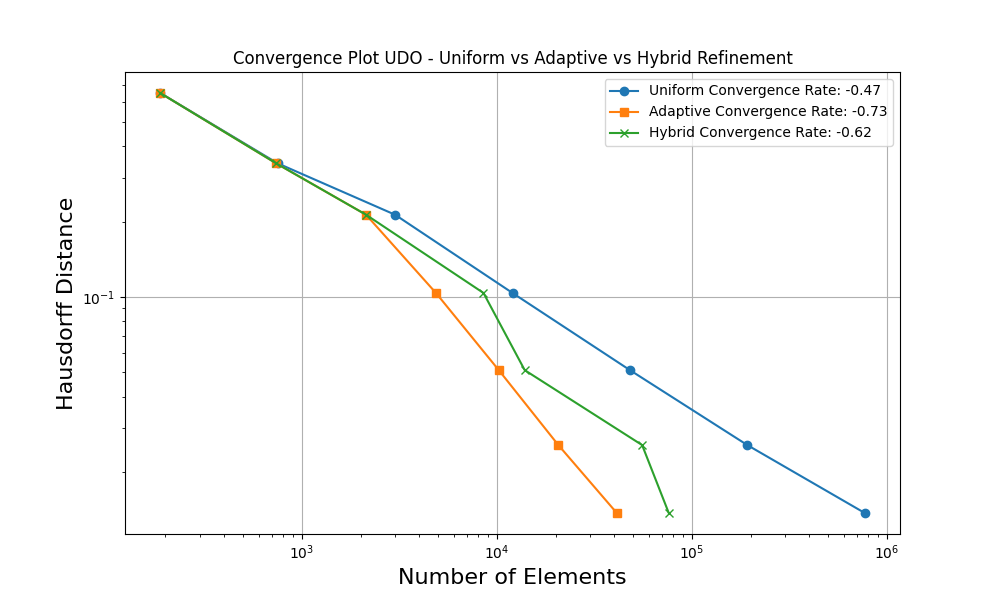}
  \caption{Hausdorff distance convergence plot for UDO. Lower is better.} 
\end{figure}

\begin{figure}[H]
  \centering
  \includegraphics[width=\textwidth]{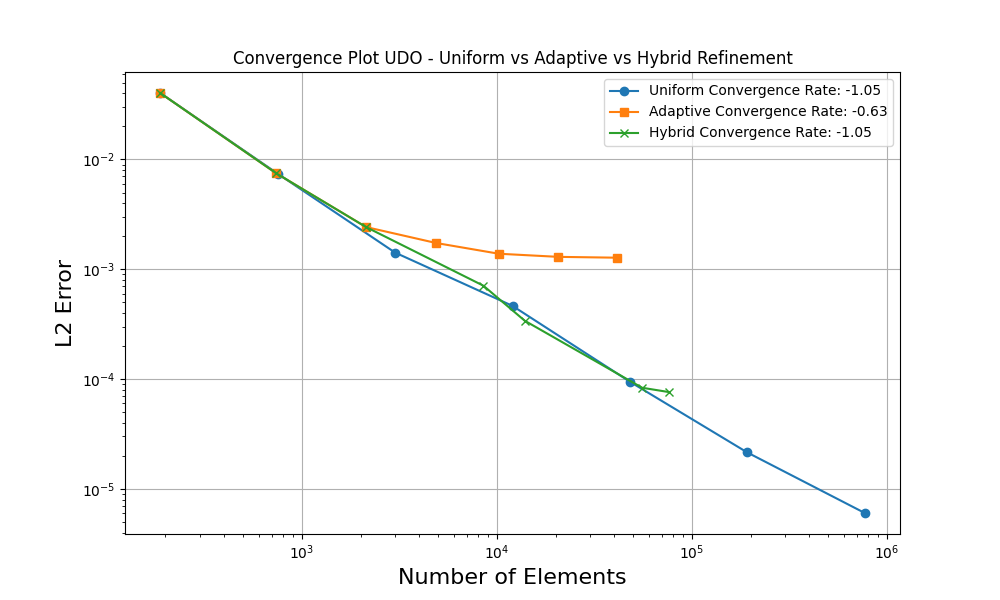}
  \caption{$L_2$ convergence plot for UDO. Lower is better.} 
\end{figure} 

For problems where the identification of the free boundary is paramount, these methods achieve the same performance as uniform refinement, as measured by Jaccard index or Hausdorff distance while using almost 2 orders of magnitude fewer elements. Our demonstration hybrid strategy was also effective in dynamically allocating resolution between approximating the free boundary and the general PDE solution. 

Both VCES and UDO begin with a the computation of a nodal active set. For problems where the active set is small e.g. the center of Figure 5.6, it is essential that the coarse mesh have enough resolution to capture the active set. Examples of where the mesh at an iteration failed to identify a small active set can be seen in the center of Figure 5.7. Choosing parameters which widen the spread of the refinement region and the use of a hybrid technique can help mitigate this issue.

\section{Parallel Considerations}
An important consideration for VI solver methods like VINEWTONRSLS is the imbalance of solver effort between mesh nodes in the active set of the converged solution where the solution is known and the corresponding inactive set where a PDE must be solved. In parallel, this imbalance can, in the worst case, cause a single process to contain most of the inactive set elements and do all the work. An indirect consequence of our proposed methods is balancing the parallel solver effort between the active and inactive set.

To demonstrate this inactive set balancing effect we consider the same ball obstacle problem as before. The inital mesh was constructed using  using the default netgen mesh constructor with the triangle height parameter set to $.45$ and then it was refined 5 times using VCES with thresholding parameters $(.2, .8)$ resulting in a final adaptive mesh with 2736 elements. The problem was then solved on a uniform mesh 3008 elements. Both meshes were then distributed across 5 processes and the number of active and inactive elements for each process was recorded. The ratio between inactive elements and total inactive elements was calculated for each process. The closer these ratios are to .20, which is the ratio of the exact active area over the inactive area, the more balanced the solver effort will be across all processes. As we can see from the table below, in the uniform mesh distribution, process 1 has significantly fewer inactive elements than other processes.

\begin{table}[H]
  \centering
  \label{tab:element_counts}
  \begin{tabular}{|c|c|c|c|c|c|c|} \hline
    \multicolumn{7}{|c|}{Element Counts and Ratios Across Processes} \\ \hline
    \multirow{2}{*}{Process} & \multicolumn{3}{c|}{VCES} & \multicolumn{3}{c|}{Uniform} \\ \cline{2-7}
    & Active & Inactive & Ratio & Active & Inactive & Ratio \\ \hline
    1 & 173 & 327 & 0.171 & 51 & 219 & 0.081 \\ \hline
    2 & 188 & 390 & 0.204 & 59 & 656 & 0.242 \\ \hline
    3 & 170 & 373 & 0.195 & 44 & 570 & 0.210 \\ \hline
    4 & 162 & 415 & 0.217 & 86 & 608 & 0.225 \\ \hline
    5 & 131 & 407 & 0.213 & 60 & 655 & 0.242 \\ \hline
  \end{tabular}
  \caption{Process element counts and ratios. VCES mesh has 2736 total elements, Uniform has 3008 total elements.}
\end{table}

\begin{figure}[H]
  \centering
  \begin{subfigure}[b]{0.45\textwidth}
      \centering
      \includegraphics[width=\textwidth]{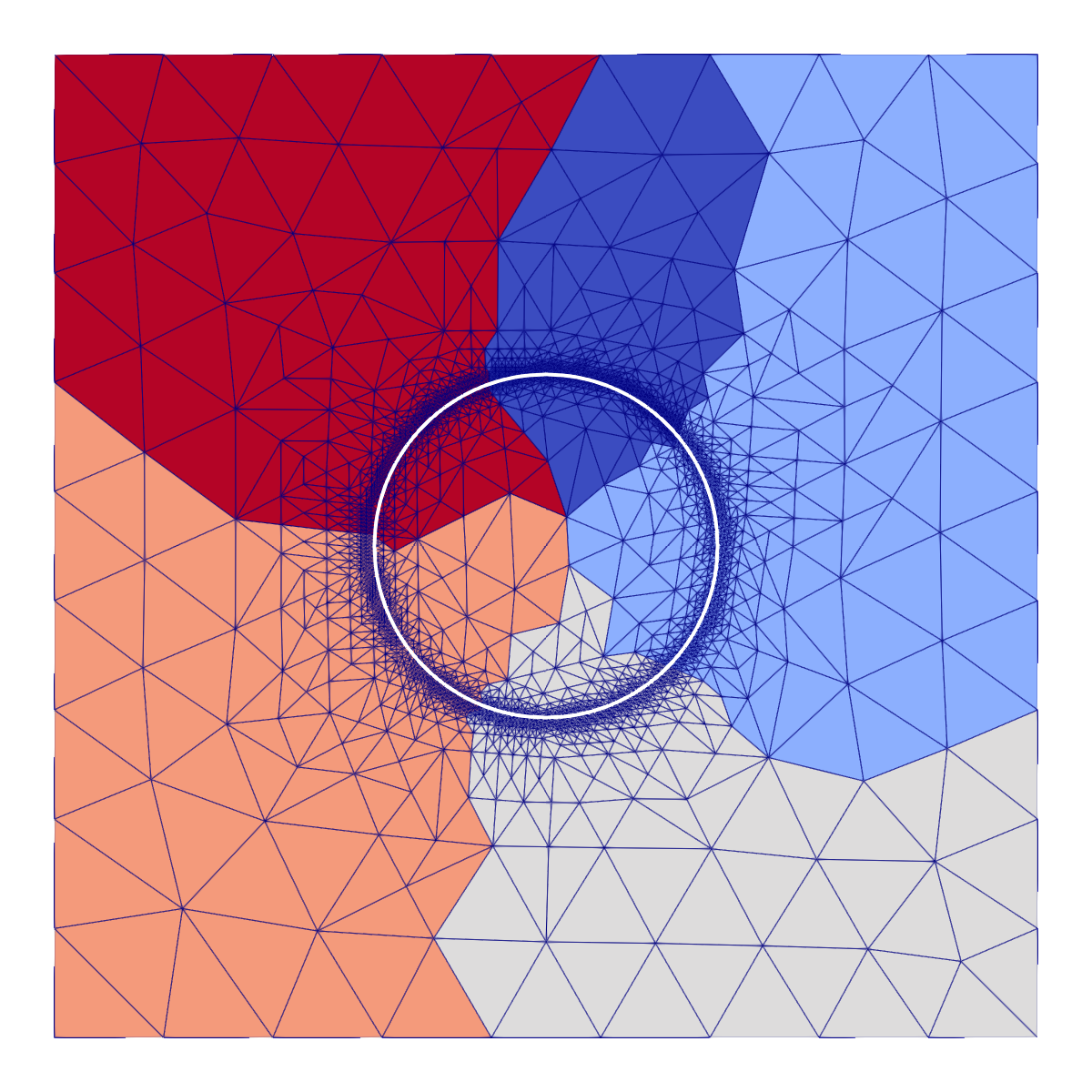}
      \label{fig:image1}
  \end{subfigure}
  \hfill
  \begin{subfigure}[b]{0.45\textwidth}
      \centering
      \includegraphics[width=\textwidth]{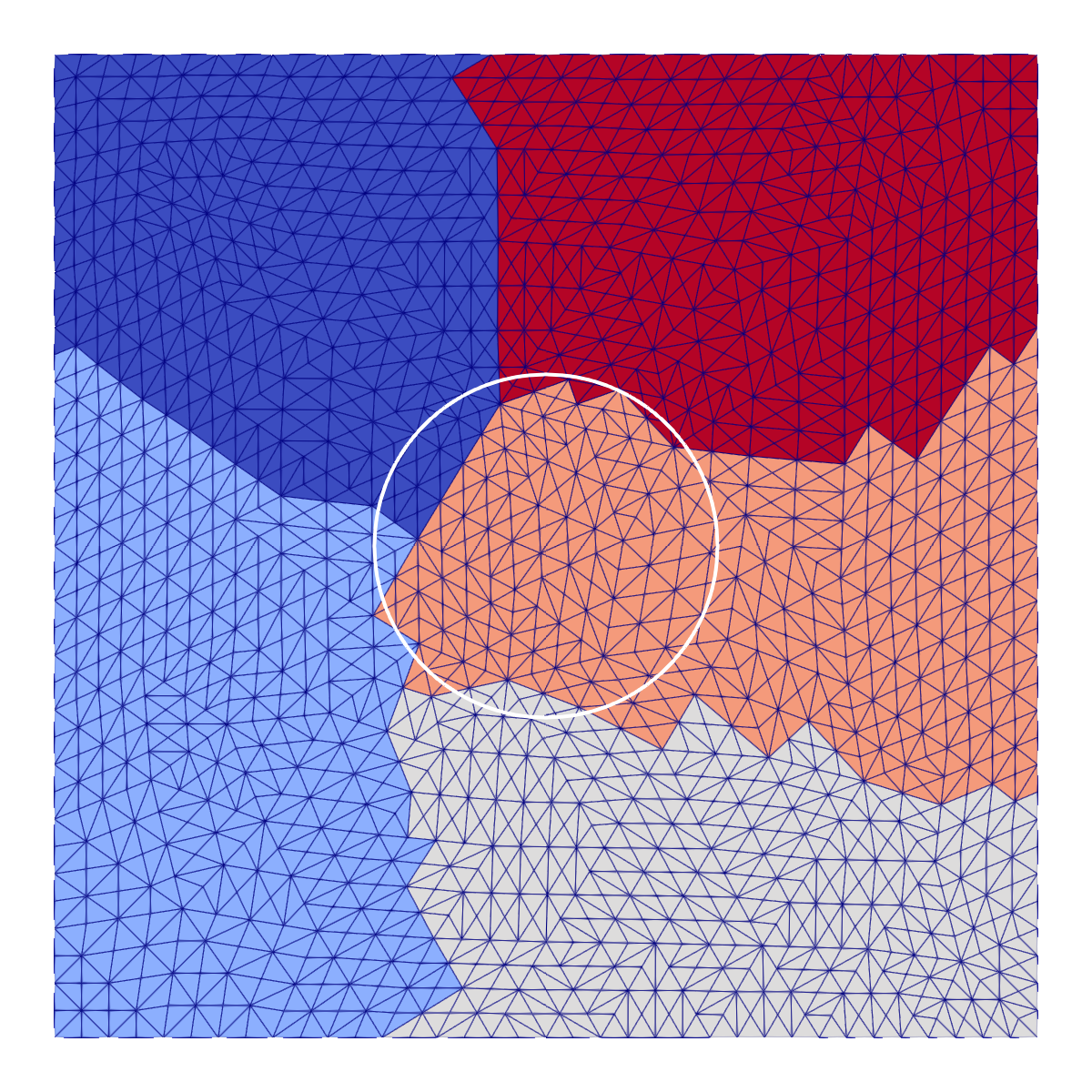}
      \label{fig:image2}
  \end{subfigure}
  \vspace*{.25cm}
  \caption{Distibuted meshes constructed via VCES(left) and uniform(right) refinement. Colors denote process, $\Gamma_u$ is displayed in white.}
  \label{fig:two_images}
\end{figure}

\chapter{Discussion}
The results in this project demonstrate the effectiveness of the proposed methods in refining about the free boundary, as well as providing load balancing benefits in parallel applications. As a consequence of the Firedrake-Netgen integration our implementations produce meshes with no hanging nodes, making our methods easily implementable with a variety of VI solvers and settings. Like other tagging methods, VCES and UDO expose parameters to the user which are capable of varying the spread of the refinement region. 

Further work will involve a deeper investigation into effective hybrid refinement schemes that are capable of balancing the error that arises from the computed free boundary versus the discretization error in the inactive set. Implementing the refinement step using PETSc's adaptation methods like DMAdaptLabel and DMAdaptMetric will allow us to investigate the effectiveness of these methods in a multigrid context as well as exploring the use of metric-based refinement methods.

\bibliography{References.bib}
\clearpage

\titleformat{\chapter}[hang]
  {\normalfont\normalsize\bfseries\centering}{Appendix \thechapter:}{3 pt}{\normalsize}

\appendix

\chapter{Supplemental Material}
This supplemental material is related to the methods and experiments discussed in this project. The referenced codes, which are integral to reproducing our results and implementing the described methods, are available in the following GitHub repository:

\begin{itemize} 
  \item \textbf{GitHub Repository}: \url{https://github.com/StefanoFochesatto/VI-AMR}
  \item \textbf{Contents}: The repository includes scripts for the numerical experiments demonstrated throughout this project. 
  \item \textbf{Usage Instructions}: Detailed instructions for running the code and reproducing the results can be found in the repository's \texttt{README.md} file. 
 \end{itemize}

\end{document}